\documentclass[bj,authoryear]{imsart}
\usepackage{subfiles}

\RequirePackage{amsthm,amsmath,amsfonts} 
\RequirePackage[colorlinks,citecolor=blue,urlcolor=blue]{hyperref}
\RequirePackage{graphicx}

\usepackage{multirow}   
\usepackage{verbatim}   
\usepackage[normalem]{ulem} 
\usepackage{xr}
\usepackage[section]{placeins}

\startlocaldefs
\newtheorem{theorem}{Theorem}[section]

\newtheorem{lemma}[theorem]{Lemma}

\newtheorem{proposition}[theorem]{Proposition}
\theoremstyle{remark}
\newtheorem{algorithm}[theorem]{Algorithm}

\newtheorem{remark}[theorem]{Remark}

\def\bs{\boldsymbol}
\newcommand{\ind}{\text{\bf 1}}

\newcommand{\re}{\mathbb{R}}
\newcommand{\ep}{\mathbb{P}}
\newcommand{\ol}[1]{\overline{#1}}

\newcommand\red[1]{#1}

\usepackage{longtable}
\usepackage{algorithm,algpseudocode}

\usepackage{multirow}   
\usepackage{verbatim}   
\usepackage[normalem]{ulem} 
\usepackage{xr}

\renewcommand{\thetable}{S\arabic{table}}

\endlocaldefs

\begin{document}

\begin{frontmatter}
\title{Inference for change-plane regression}
\runtitle{Inference for change-plane regression}

\begin{aug}
\author[A]{\fnms{Chaeryon} \snm{Kang}\ead[label=e1, mark]{crkang@pitt.edu}\orcid{0000-0001-9442-8515}}
\author[B]{\fnms{Hunyong} \snm{Cho}\ead[label=e2,mark]{hunycho@live.unc.edu}\orcid{0000-0002-0087-0079}}
\author[C]{\fnms{Rui} \snm{Song}\ead[label=e3,mark]{rsong@ncsu.edu}}
\author[D]{\fnms{Moulinath} \snm{Banerjee}\ead[label=e4,mark]{moulib@umich.edu}}
\author[E]{\fnms{Eric B.} \snm{Laber}\ead[label=e5,mark]{eric.laber@duke.edu}\orcid{0000-0003-2640-7696}}
\author[B]{\fnms{Michael R.} \snm{Kosorok}\ead[label=e6,mark]{kosorok@unc.edu}\orcid{0000-0002-6070-9738}}

\address[A]{University of Pittsburgh, \printead{e1}}
\address[B]{University of North Carolina at Chapel Hill, \printead{e2,e6}}
\address[C]{North Carolina State University, \printead{e3}}
\address[D]{University of Michigan, \printead{e4}}
\address[E]{Duke University, \printead{e5}}
\end{aug}

\begin{abstract}
A key challenge in analyzing the behavior of change-plane estimators is that the objective function has multiple minimizers. Two estimators are proposed to deal with this non-uniqueness. For each estimator, an $n$-rate of convergence is established, and the limiting distribution is derived.  Based on these results, we provide a parametric bootstrap procedure for inference. The validity of our theoretical results and the finite sample performance of the bootstrap are demonstrated through simulation experiments.  We illustrate the proposed methods to latent subgroup identification in precision medicine using the ACTG175 AIDS study data.
\end{abstract}


\begin{keyword}[class=MSC]
\kwd[Primary ]{62F12}  
\kwd[; secondary ]{62F40,62J05,62P10} 
\end{keyword}

\begin{keyword}
\kwd{Change-plane}
\kwd{Latent subgroup analysis}
\kwd{M-estimation}
\kwd{Precision medicine}
\end{keyword}

\end{frontmatter}



\section{Introduction}\label{sec:intro}

Change-plane regression models add flexibility by allowing distinct models on either side of a 
hyperplane \citep{kang2011new, wei2013latent, huang2021threshold}.  There is increasing interest in 
these models, especially in applications where domain knowledge suggests the input space is 
partitioned into (latent) subgroups; e.g., subjects with an enhanced treatment effect 
\citep{fan2017change, wang2022multi, wei2023subgroup, jin2023change}. For example, 
\cite{huang2021threshold} applied the change-plane logistic regression model to investigate  
potential heterogeneity in the efficacy of the human immunodeficiency virus (HIV) vaccine against 
HIV infection among participants in a phase 3 clinical trial study. Two subgroups were 
identified by a change-plane consisting of the four single nucleotide polymorphisms (SNPs, 
rs114945036, rs138747765, rs145835719, and rs75898867). A likelihood ratio test based on the change-plane
supported the existence of two subgroups: one with a higher vaccine efficacy (ranging from 
58.6-58.9\%) and another with lower efficacy (29.2-29.5\%). This result can be used to inform 
personalized recommendations HIV vaccination based on individual genotypes and demographic 
characteristics. \cite{li2021multithreshold} applied a multithreshold change-plane model to
identify latent subgroups among Scleroderma patients randomized into oral native collagen at a dose 
of 500 $\mu g$/day or a placebo. 
The change-plane regression model identified two subgroups where the association between the 
Modified Rodnan skin score and 11 demographic and biomarker covariates differed by a change-plane 
consisting of three variables: health assessment questionnaire,  patient self-assessment of disease 
progress, and lung performance measurement. The application identified heterogeneity by sex, which 
was missed by an ordinary regression model. 

While the preceding examples illustrate the clinical and scientific value of
change-plane models, their wider application is impeded by a lack of rigorous inferential methods. 
There are two major challenges to establishing theory for change-plane models: non-smoothness of 
the objective function and non-uniqueness of the estimated parameters.   The non-smoothness occurs 
because the model allows for an abrupt change in parameters at a boundary defined by a hyperplane; 
consequently, standard M-estimation theory does not apply.  The change in parameters is encoded as 
an indicator function, and because of the non-injective binary topology of indicator functions, 
multiple parameter values form a level set which maps to a single objective value. And, thus, the 
minimizer of the objective is not unique. This problem occurs not only in finite samples but also 
asymptotically.  Consequently, the argmin (or argmax) continuous mapping theorem is not applicable.

A consequence of the non-smoothness of the objective function in change-plane models is that it is possible to obtain 
convergence rates faster than $\sqrt n$.  In the special case of a change-\emph{point} model, multiple examples have been 
given in which estimated parameters converge at $n$-rate to the minimizer of a non-Gaussian process.  For example, 
\cite{kosorok2007inference} showed that a transformation model change-point estimator is $n$-consistent and weakly 
converges to the minimizer of a right-continuous jump process; the change-point problem in the Cox model 
\citep{pons2003estimation} is considered a special case of the change-point model in \cite{kosorok2007inference}, 
and similar asymptotic properties were presented \citep{pons2003estimation};  and \cite{song2016asymptotics} showed a least-squares
estimator (LSE) in a linear change-point model can attain $n$-rate convergence under a correctly specified 
model but that these rates can be slower under misspecification. Similar $n$-rate of convergence results of 
the maximum (profile) likelihood estimators (MLEs) or LSE in change-point models were 
presented earlier in the literature \citep{chernoff1956estimation, luo1996asymptotic, gijbels1999estimation}. The 
asymptotic theory of the change-point model is most well-developed in non-linear time series analysis. For example, 
\cite{chan1993consistency} showed that the conditional LSE of the change-point (also called a ``threshold") 
parameter has an $n$-consistency, and its limiting distribution is a functional of a compound Poisson process. 
\cite{hansen2000sample} established weak convergence of the change-point estimator through the inverse of the 
likelihood ratio statistic, with the convergence rate of $n^{1-2\alpha}$ for $0<\alpha<1$, assuming a diminishing threshold 
effect as sample size increases. Similar rate of convergence results for the change-point estimators in the time 
series analysis have been reported \citep{bai1998estimating}, where the convergence rate depends on the presence of a 
diminishing threshold effect assumption. Our study focuses on asymptotic theory for 
the change-\emph{plane} in spatially or 
temporally independent data in the presence of a non-zero constant threshold effect. 

The non-uniqueness issue in change-point estimators can be handled through a modification of the argmin (or argmax) continuous mapping 
theorem. For example, using such a modified theorem, the weak convergence of the smallest and the largest argmin can be 
obtained in the presence of multiple minimizers \citep{ferger2004continuous}. \cite{lan2009change} proposed a multistage 
adaptive procedure to estimate the change-point in a parametric regression model and investigated the asymptotic properties 
of the change-point estimator. They proved that the ``zoom-in" stage change-point estimator achieves $n$-consistency, and 
that the joint asymptotic distribution of the smallest and largest minimizers of the zoom-in estimators is given by the 
smallest and largest minimizers of a two-sided compound Poisson process. \cite{seijo2011continuous} proposed a version of 
the argmax continuous mapping theorem to deal with multiple maximizers in the change-point regression model defined in the 
Skorohod topology; the joint weak convergence of the smallest and largest argmax of the stochastic process was obtained 
when both the stochastic process and the associated pure jump process converged in the Skorohod topology with the smallest 
and largest argmax.

Among the few studies which have investigated the asymptotic behavior of change-plane models with more than two change-plane parameters, most replace the indicator function with a smooth surrogate or pose a diminishing threshold to ease the computational and/or theoretical burden. 
In this case, the convergence rate can be slower than $n$ depending on the degree of the smoothness and the rate of the diminishing threshold. For example, the change-plane estimator of \cite{seo2007smoothed}, which uses an integrated smooth kernel function, has a $\sqrt{n/\sigma_n}$ rate of convergence, attains a maximal rate of $n^{-3/4}$ under the choice of bandwidth $\sigma_n=\log n/\sqrt{n}$, and converges weakly to a normal distribution when $\log n/n \sigma_n^2 \rightarrow 0$. \cite{mukherjee2020asymptotic} showed that their kernel smoothed change-plane estimator can achieve a faster rate, close to $n^{-1}$, provided $n \sigma_n \rightarrow \infty$. \cite{lee2021factor} developed a change-plane estimator without modifying the change-plane indicator function. However, their weak convergence results rely on a diminishing threshold assumption; with a rate $0<\alpha<1/2$, their estimator converges to a non-Gaussian distribution at $n^{1-2\alpha}-$rate.

We establish the large-sample behavior of the prototypical linear change-plane estimator with an arbitrary number of invariant change-plane parameters without modifying the indicator function.  Most closely related to our work is the recent paper of \cite{deng2022maximum} in which they study asymptotic properties of the MLE of the change-plane Cox proportional hazards model. They established $n$-rate convergence and suggested the $m$-out-of-$n$ bootstrap for inference.  However, the non-uniqueness of the argmax was not fully addressed as the smallest argmax used in their approach is not well-defined in the case of multi-dimensional change-plane parameters, and thus they did not fully characterize the limiting distribution.  In addition, they required continuity of the change plane covariates and a bounded parameter space.

We propose M-estimators for the change-plane and establish their asymptotic properties under minimal assumptions; e.g., the change-plane covariate space is allowed to be at least partially discrete. Among many possible ways of summarizing the level set of the argmin, we propose two estimators based on two most intuitive approaches: mean and mode. We establish consistency and the rate of convergence of the proposed M-estimators using empirical process techniques.  The proposed estimators achieve a fast $n$-rate of convergence for the change-plane parameter. The most significant contribution of the article is to fully characterize the limiting distribution of the M-estimators, which are minimizers of sums of two compound jump processes.  Using this limiting distribution, we derive a parametric bootstrap procedure which is more efficient than the $m$-out-of-$n$ bootstrap suggested by \cite{deng2022maximum}. Simulation studies confirm the theoretical $n$-rate of convergence as well as the consistency of the parametric bootstrap distribution.

The remainder of this article is organized as follows. We describe the change-plane regression model in 
Section \ref{sec:dmod} and introduce two types of estimators in Section \ref{sec:est}. Consistency, rates of 
convergence, and the weak convergence of the proposed estimators are established in Sections 
\ref{sec:rateconv} and \ref{sec:WeakConvergence}. The parametric bootstrap procedure is described in Section  
\ref{sec:inference}. 
Section \ref{sec:numeric} provides numerical estimation procedures, and Section \ref{sec:sim} presents simulation studies of weak convergence and the 
validity of the parametric bootstrap.  We illustrate the application of the proposed methods using data from 
an AIDS study in Section \ref{sec:data_ACTG175}.  We conclude the paper with a discussion of 
future research directions in Section \ref{sec:conclusion}.

\section{Data and model}\label{sec:dmod}

We assume that the data consists of $n$ i.i.d. copies of the triple $(Y,Z,X)$, where $X\in\re^p$ 
are the change-plane covariates (to be defined shortly), $Z\in\re^d$ are the regression covariates, and 
$Y\in\re$ is the regression outcome.  In addition, we assume that $Y$ satisfies
\begin{eqnarray}
Y&=&\beta_0'Z \ind\{\omega_0'X-\gamma_0\leq 0\}     +\delta_0'Z \ind\{\omega_0'X-\gamma_0>0\}+\epsilon,\label{eq:dmodel}
\end{eqnarray}
where: $\epsilon$ is a continuous random variable which is independent of $(X,Z)$ and with mean zero and variance $0<\sigma^2<\infty$; $\beta_0,\delta_0\in\re^d$ are the regression parameters; $\omega_0\in S^{p-1}$ and $\gamma_0\in\re$ are the change-plane parameters, prime denotes transpose, and 
$S^{p-1}$ is the $p-1$ dimensional unit sphere embedded in $\re^p$ (i.e., $\left\lbrace \omega \in \re^p\,:\,\|\omega\|=1\right\rbrace$), for $p\geq 1$. We have used a subscript zero to denote the true parameter
values; we omit a subscript when referencing generic parameter values.
Let $S^0$ to be the set of points $\{-1,1\}$, and define the composite parameters $\zeta=(\beta,\delta)$, $\phi=(\omega,\gamma)$, and $\theta=(\zeta,\phi)$.  Define $\ol{\omega}_0$ to be the $p\times(p-1)$ matrix consisting of the orthonormal basis vectors
for the $p-1$-dimensional subspace in $\re^p$ which is orthogonal to the linear span of $\omega_0$; we
denote this subspace by  $\re^{p-1}_{\ol{\omega}_0}$.  

\begin{remark}\label{rem2}
The preceding model has the following interesting special case: Assume
\begin{eqnarray}
Y&=&\beta_0'\tilde{Z} \ind\{\tilde{X}\leq \gamma_0\}+\delta_0'\tilde{Z} \ind\{\tilde{X}>\gamma_0\}+\epsilon, \label{model2.e1}
\end{eqnarray}
where $\tilde{Z}=MX$, for some $(d+1)\times p$ matrix $M$, and $\tilde{X}=M_1'X$, where $M_1'$ is the first row of $M$, and where $M$ is unknown. This model posits that under an unknown linear rotation of $X$, there is a change-point (instead of a change-plane) regression model of the form given in~(\ref{model2.e1}). When the true value of $M_1$ is proportional to $\omega_0$, (\ref{model2.e1}) becomes a special case of (\ref{eq:dmodel}). Thus,  the change-plane regression model can be viewed as a natural extension of the change-point model when the choice of change-point variable is unknown.
\end{remark}

We will make use of the following conditions.
\begin{itemize}
\item[C1.] $\|X\|\leq k_1<\infty$ almost surely and its covariance is full rank.
\item[C2.] The random variable $U\equiv\omega_0'X-\gamma_0$ is continuous, and, in a neighborhood of zero, it has a continuous density $f$. Let $f_0\equiv f(0)$, then $0<f_0<\infty$. This implies the existence of $0<c_0,c_1<\infty$ such that  $f(u)\geq c_0$ for all $u\in[-c_1,c_1]$, and hence also implies that $P\{U\leq 0\}>0$ and $P\{U>0\}>0$.
\item[C2'.] Condition C2 is strengthened to require that the density $f$ be uniformly bounded on $\re$ and continuous at zero.
\item[C3.] Whenever $(\omega,\gamma)\neq(\omega_0,\gamma_0)$, $P\{(\omega'X-\gamma)(\omega_0'X-\gamma_0)<0\}>0$. Moreover, the joint distribution of $(Z,X)$ given $U=u$ converges to a probability measure $G$ as $u\rightarrow 0$; and, moreover, the covariance of $\overline{U}=\overline{\omega}_0'X$ is full rank under $G$.
\item[C3'.] Condition C3 is strengthened to require that $\overline{U}$ be continuous under $G$.
\item[C4.] $\|Z\|\leq k_2<\infty$ almost surely, both $E\left(ZZ'|\omega_0'X-\gamma_0\leq 0\right)$ and $E\left(ZZ'|\omega_0'X-\gamma_0> 0\right)$ are full rank with minimum eigenvalues bounded below by $c_2>0$.
\item[C5.] $\beta_0\neq\delta_0$, and $P\{(\beta_0-\delta_0)'Z=0\}=0$ under $G$.
\item[C6.] $\epsilon$ has uniformly bounded density $\xi$ on $\re$, where $\xi$ is uniformly equicontinuous, i.e., \[\lim_{\eta\downarrow 0}\;\sup_{s,t\in\re:\;|s-t|<\eta}\;|\xi(s)-\xi(t)|\;=\;0.\]
\end{itemize}
These assumptions are, overall, very mild. The most stringent are Condition C5 which 
implies existence of a change-plane, and Condition C2 which requires that the change-plane covariates ($X$) not
all be discrete. We further discuss each of the conditions.
Condition C1 implies that, without loss of generality, we can assume that $\gamma\in[l_0,u_0]$, where $l_0=-k_1-\rho$ and $u_0=k_1+\rho$, for any $\rho>0$, because $\ind\{\omega'X-\gamma\leq 0\}\geq\ind\{\omega'X-l_0\leq 0\}$ for all $\gamma\leq l_0$ and
$\ind\{\omega'X-\gamma\leq 0\}\leq\ind\{\omega'X-u_0\leq 0\}$ for all $\gamma\geq u_0$, almost surely, for any value of $\omega$. Accordingly, we fix $0<\rho<\infty$, $l_0=-k_1-\rho$ and $u_0=k_1+\rho$ hereafter. Condition C2 characterizes the limiting distribution of $U$ near zero and plays a crucial role in deriving the asymptotic distribution of the change-plane estimators described in the next section. This condition is easily satisfied if $X$ is continuous and full rank with density bounded above and below. It is also satisfied if $X$ includes both continuous and categorical variables, provided the weights in the vector $\omega_0$ are positive for at least one of the continuous variables. The strengthening of Condition C2 to C2' is used for the parametric bootstrap results of Section \ref{sec:inference}. C3 is needed for identifiability of $(\omega,\gamma)$ as well as stability in the limit when $U$ approaches zero. The strengthening of Condition C3 to C3' is used for weak convergence of one of the two estimators which will be introduced in Section 3. Condition C4 is needed for identifiability of $(\beta,\delta)$ as well as stability when $U$ is close to zero. Condition C5 is also crucial for model identifiability of $(\omega,\gamma)$, and violation of the $\beta_0\neq\delta_0$ condition is regarded as the ``absence of a change-plane." Condition C6, which is a strengthening of the stated properties of the residual $\epsilon$ above, is used to establish the parametric bootstrap results of Section \ref{sec:inference}.

\section{Estimation}\label{sec:est}

Our focus is on estimating $\theta=(\beta,\delta,\omega,\gamma)$. We will use the M-estimator obtained by minimizing
\begin{eqnarray}
\;\theta\mapsto M_n(\theta)&\equiv&\ep_n\left[\ind\{\omega'X-\gamma\leq 0\}(Y-\beta'Z)^2+\ind\{\omega'X-\gamma>0\}(Y-\delta'Z)^2\right],
\label{eq:mest}
\end{eqnarray}
where $\ep_n$ is the standard empirical measure. We have established previously that the range of $\theta$ is $K=K_1\times K_2$, where $K_1=\re^d\times\re^d$ and $K_2=S^{p-1}\times[l_0,u_0]$. Let $A_0(\phi)=\ind\{\omega'X-\gamma\leq 0\}$ and $A_1(\phi)=\ind\{\omega'X-\gamma>0\}$, and define $D_{0n}(\phi)=\ep_n[ZZ'A_0(\phi)]/[\ep_n A_0(\phi)]$ and $D_{1n}(\theta)=\ep_n[ZZ'A_1(\phi)]/[\ep_n A_1(\phi)]$. Fix a $c_3\in(0,c_2)$, where $c_2$ comes from assumption C4, and let $K_{2n}$ be the set of $\phi\in K_2$ such that the smallest eigenvalues of $D_{0n}(\phi)$ and $D_{1n}(\phi)$ are bounded below by $c_3$. Let $K_{2n}'\supset K_{2n}$ be similarly defined but with the weaker requirement that the minimum eigenvalues of $D_{0n}(\phi)$ and $D_{1n}(\phi)$ are both positive. Define $K_n=K_1\times K_{2n}$ and $K_n'=K_1\times K_{2n}'$, and note that membership in either $K_{2n}$ or $K_{2n}'$ does not involve $\zeta = (\beta,\delta)$.
Unfortunately, $K_n$ (and even $K_n'$) may be empty in small samples. To address this, 
we set $K_n^{\ast} = K_n$ if $K_n$ is non-empty, $K_n^{\ast}=K_n'$ if $K_n$ is empty but 
$K_n'$ is non-empty, and $K_n^{\ast} = K$ otherwise. 
Let $\tilde{\theta}_n$ be the $\arg\min$ over $K_n^{\ast}$ of $M_n(\theta)$. We will, however, need to
require that $K_n$ be non-empty for all sufficiently large $n$ to ensure identifiability of $(\beta,\delta)$. 
Combining Assumption C4 with the strong law of large numbers, it can be seen that $K_n$ will be non-empty and contain $\theta_0$ for all sufficiently large $n$, almost surely. Because $K_n\subset K_n'$, the same is true for $K_n'$.

An additional complication with $M_n(\theta)$ is the presence of certain level sets. To clarify this, fix $(\omega,\gamma)$ at a point $(\omega_1,\gamma_1)\in  S^{p-1}\times [l_0,u_0]$, and let $V_i=\ind\{\omega_1'X_i-\gamma_1>0\} -\ind\{\omega_1'X_i-\gamma_1\leq 0\}$, $1\leq i\leq n$. Let
\[\Phi_n=\left\{(\omega,\gamma)\in S^{p-1}\times[l_0,u_0]: \ind\{\omega'X_i-\gamma> 0\}-\ind\{\omega'X_i-\gamma\leq 0\}=V_i,\;1\leq i \leq n\right\},\]
and note that $(\omega,\gamma)\mapsto M_n(\beta,\delta,\omega,\gamma)$ is constant over $(\omega,\gamma)\in\Phi_n$ and that
$\Phi_n$ is a bounded set by construction. The sets of this form are defined by all possible realizations of $\bs V_n$ 
$\equiv(V_1,\ldots,V_n)\in\{-1,1\}^n$ that are obtained from partitioning $X_1,\ldots,X_n$ by hyperplanes in $\re^p$ into two non-empty groups. Suppose that $n$ is large enough to ensure that $K_{2n}'$ is non-empty and that $(\omega_1,\gamma_1)\in K_{2n}'$. Then $\Phi_n\subset K_{2n}'$ and there exists $j,k\in\{1,\ldots,n\}$ such that $V_j=-1$ and and $V_k=1$. 
Because the level set that minimizes (\ref{eq:mest}) is not necessarily a singleton, we define two types of midpoints of $\Phi_n$ which we can use to construct unique estimators. Define the maps $\omega\mapsto C_L^n(\omega)=\max_{1\leq i\leq n:V_i=-1}\omega'X_i$, $\omega\mapsto C_U^n(\omega)=\min_{1\leq i\leq n:V_i=1}\omega'X_i$, and $\omega\mapsto C_R^n(\omega)=C_U^n(\omega)-C_L^n(\omega)$. By construction, $\Phi_n=\{\phi\in K_2: C_L^n(\omega)-\gamma\leq 0<C_U^n(\omega)-\gamma\}$, and thus both $C_R^n(\omega)>0$ and $C_L^n(\omega)\leq\gamma<C_U^n(\omega)$ if and only if $(\omega,\gamma)\in\Phi_n$. 

We discuss two canonical midpoints---mean- and mode-midpoints. The $\omega$ part of the mean-midpoint is constructed as the weighted average of the $\omega$'s in the level set with their thickness ($C_R^n(\omega)$) as the weight, and that of the mode-midpoint is the $\omega$ element with the largest $C_R^n(\omega)$ value in the level set. Given the $\omega$ value of each midpoint approach, the $\gamma$ part is given as the median of the $\gamma$'s in the level set. A mode-midpoint corresponds to the hyperplane that maximizes the margin and hence, maximizes the distances between two latent classes. However, mode-midpoints are not necessarily unique, and as a result, a stronger, continuity condition (C3') is required for the weak convergence results than mean-midpoints. As mean-midpoints do not have such restriction, we use them as the baseline estimator.

For the first type of midpoint, we define the ``mean-midpoint'' of $\Phi_n$ to be $(\hat{\omega},\hat{\gamma})$, where 
\begin{align}
\label{omega.hat}
\hat{\omega}=\frac{\int_{R_n}\omega C_R^n(\omega)d\nu(\omega)}{\left\|\int_{R_n}\omega C_R^n(\omega)d\nu(\omega)\right\|},
\end{align}
$R_n=\{\omega\in S^{p-1}:C_R^n(\omega)>0\}$, $\nu$ is the uniform measure on $S^{p-1}$, and $\hat{\gamma}=\left[C_L^n(\hat{\omega})+C_U^n(\hat{\omega})\right]/2$. In Lemma~\ref{meanmidexistence} in Section \ref{asec:est} of Supplementary Material \citep{kang2024Inference}, we show that the mean-midpoint exists, is well-defined, and is contained in $\Phi_n$. Note that from a computational perspective, it is not difficult to approximate $(\widehat{\omega}, \widehat{\gamma})$ 
to any desired degree of accuracy using rejection sampling; moreover, it is easy to verify that any Monte Carlo estimate based on rejection sampling will be contained in $\Phi_n$. 

For the second type of midpoint, we define the ``mode-midpoint'' of the level set $\Phi_n$ to be $(\check{\omega},\check{\gamma})$, where $(\check{\omega},\check{\gamma})$ is the $\arg\max$ over $\Phi_n$ of
\begin{eqnarray}
&&(\omega,\gamma)\mapsto\lambda,\;\mbox{subject to}\nonumber\\
&&C_L^n(\omega)-\gamma\leq -\lambda\;\mbox{and}\:C_U^n(\omega)-\gamma\geq\lambda.\label{eq1.modemin}
\end{eqnarray}
We assume, as before, that $n$ is large enough so that $K_{2n}'$ is non-empty and that $(\omega_1,\gamma_1)\in K_{2n}'\cap\Phi_n$. 
In Lemma~\ref{modemidexistence} in Section \ref{asec:est} of Supplementary Material, we show that the mode-midpoint is unique, is contained in $\Phi_n$, and that
$$\check \gamma = \frac{C_L^n(\check\omega) + C_U^n(\check\omega)}{2}.$$

Now we define our two estimators, $\hat{\theta}_n$ and $\check{\theta}_n$, based on the two midpoint approaches. Let 
$\tilde{\theta}_n$ be as defined at the beginning of this section, and let $\Phi_n$ be the level set associated with 
$\tilde{\phi}_n$. If $n$ is small enough so that $K_n'$ is empty, then let 
$\hat{\theta}_n=\check{\theta}_n=\tilde{\theta}_n$. Otherwise, define 
$\hat{\theta}_n=
(\tilde{\beta}_n,\tilde{\delta}_n,\hat{\omega}_n,\hat{\gamma}_n)$, where $(\hat{\omega}_n,\hat{\gamma}_n)=(\hat{\omega},\hat{\gamma})$ is from (\ref{omega.hat}) for the level set $\Phi_n$ that minimizes (\ref{eq:mest}).
Similarly, define $\check{\theta}_n=(\tilde{\beta}_n,\tilde{\delta}_n,\check{\omega}_n,\check{\gamma}_n)$, where $(\check{\omega}_n,\check{\gamma}_n)=(\check{\omega},\check{\gamma})$ from (\ref{eq1.modemin}) for the level set $\Phi_n$ that minimizes (\ref{eq:mest}).


\section{Consistency and rates of convergence}\label{sec:rateconv}


In this section, we derive consistency and rates of convergence for the proposed M-estimators of the change-plane parameters and regression parameters. First, we establish consistency of $\hat{\theta}_n$ and $\check{\theta}_n$ in the following theorem, of which proof is in Appendix \ref{asec:con} of the Supplementary Material.

\begin{theorem}\label{thm.1}{(consistency)}
Let $\tilde{\theta}_n$ be any sequence in $K_n^{\ast}$ satisfying $M_n(\tilde{\theta}_n)\leq M_n(\theta_0)$ for all $n$ large enough almost surely.
Under Conditions {C1, C2, C3, C4 and C5}, $\tilde{\theta}_n\rightarrow\theta_0$ almost surely.
\end{theorem}
That both $\hat{\theta}_n$ and $\check{\theta}_n$ are sequences satisfying the conditions of Theorem~\ref{thm.1} implies their consistency. 

In the next theorem, we obtain the rates of convergence from the limiting behavior of the process $(M_n-M)(\theta)$; our line of argument follows Corollary 14.5 in \cite{kosorok2008introduction}.
Recall that $\zeta = (\beta, \delta)$ and $\phi=(\omega,\gamma)$. 

\begin{theorem}{(rates of convergence)}
\label{theorem:rateofconvergence1}
Assume Conditions {C1, C2, C3, C4 and C5} and let $M_n(\tilde{\theta}_n)=\min_{\theta\in K_n^{\ast}}M_n(\theta)$. Then $\sqrt{n} ~\|\tilde{\zeta}_n-\zeta_0\|=O_p(1)$ and ~$n (\|\tilde{\omega}_n-\omega_0\|+|\tilde{\gamma}_n-\gamma_0|)=O_p(1)$.
\end{theorem}
Because the sequences $\hat{\theta}_n$ and $\check{\theta}_n$ both satisfy the given conditions, Theorem~\ref{theorem:rateofconvergence1} yields that
$\sqrt{n}\|\hat{\zeta}_n-\zeta_0\|=O_P(1)$ and $n\|\hat{\phi}_n-\phi_0\|=O_P(1)$ and also $\sqrt{n}\|\check{\zeta}_n-\zeta_0\|=O_P(1)$ and $n\|\check{\phi}_n-\phi_0\|=O_P(1)$.

\begin{proof}[Proof of Theorem~\ref{theorem:rateofconvergence1}]
We verify the conditions of Corollary~14.5 of \cite{kosorok2008introduction} for the convergence rate $r_n=\sqrt{n}$ and
for convergence of a minimizer rather than a maximizer. This latter adjustment is straightforward because minimizers of an objective function are also maximizers of the negative of the objective function. Before giving the conditions, recall from expression~(\ref{eq:mest}) that $M_n(\theta)=\ep_n m_{\theta}(X,Y,Z)$, where
\[m_{\theta}(Y,Z,X)=\ind\{\omega'X-\gamma\leq 0\}(Y-\beta'Z)^2+\ind\{\omega'X-\gamma>0\}(Y-\delta'Z)^2,\]
and define $M(\theta)=Pm_{\theta}(X,Y,Z)$. 
Also, for any $\eta>0$, define $\Theta_{\eta}\equiv\{\theta\in K:\tilde{d}(\theta,\theta_0)\leq\eta\}$,
where $\tilde{d}^2(\theta, \theta_0)=\|\zeta-\zeta_0\|^2+\|\omega-\omega_0\|+|\gamma-\gamma_0|$. 
Here are the needed conditions for Corollary~14.5 of \cite{kosorok2008introduction}:
\begin{enumerate}
\item[R1] $\tilde{\theta}_n$ is consistent in outer probability for $\theta_0$ and satisfies $\ep_n m_{\tilde{\theta}_n}\leq \inf_{\theta\in\Theta_{\eta}}\ep_n m_{\theta}+O_p(n^{-1})$ for some $\eta>0$.
\item[R2] $M(\theta)-M(\theta_0)\geq c\tilde{d}^2(\theta,\theta_0)$ for all $\theta\in\Theta_{\eta}$ and some $0<c,\eta<\infty$.
\item[R3] There exists an increasing function $\varphi:[0,\infty)\mapsto[0,\infty)$ such that the following hold for some $0<c<\infty$:
\begin{enumerate}
\item $\eta\mapsto\varphi(\eta)/\eta^{\alpha}$ is decreasing for some $\alpha<2$,
\item $n\varphi(1/\sqrt{n})\leq c\sqrt{n}$ for some $0<c<\infty$, and
\item $E^{\ast}\|\mathbb{G}_n\|_{{\cal M}_\eta}\leq c\varphi(\eta)$, for all $0\leq\eta\leq\eta_1$ and some $0<c,\eta_1<\infty$ and all $n$ large enough, where
\[{\cal M}_{\eta}=\left\{m_{\theta}(Y,Z,X)-m_{\theta_0}(Y,Z,X)\;:\;\theta\in\Theta_{\eta}\right\}.\]
\end{enumerate}
\end{enumerate}
In the last condition, $E^{\ast}\|\mathbb{G}_n\|_{A}$ denotes the outer expectation of the uniform norm of the normalized empirical process over a set $A$, and the formal definition of each element can be found in Section \ref{asec:empirical_notation} of the Supplementary Materials.
Provided these conditions are satisfied, the corollary yields that $\sqrt{n}\tilde{d}(\tilde{\theta}_n,\theta_0)=O_p(1)$, and the conclusions of Theorem~\ref{theorem:rateofconvergence1} follow.

We will now verify each of the conditions R1--R3.
For R1, Conditions {C1, C2, C3, C4 and C5} imply that for $0<\eta$ small enough, $\Theta_\eta\subset K_n$ for all $n$ large enough almost surely. Combining this with
the fact that $M_n(\tilde{\theta}_n)\leq M_n(\theta_0)$ for all $n$ large enough almost surely, we obtain consistency of $\tilde{\theta}_n$ from Theorem~\ref{thm.1}. Combining these two results yields that Condition R1 holds.

The following proposition verifies that Condition R2 also holds:
\begin{proposition}\label{prop.rate1}
Under the given conditions, $M(\theta)-M(\theta_0) \geq c_4\tilde{d}^2(\theta, \theta_0)$, for all $\theta\in\Theta_{\eta_5}$, for some $0<c_4,\eta_5<\infty$.
\end{proposition}

For R3, we will use $\varphi(\eta)=\eta$. Part a) is satisfied for any $1<\alpha<2$. Part b) is satisfied because $n\varphi(1/\sqrt{n})=\sqrt{n}$. Finally, Part c) is satisfied through the following proposition.
\begin{proposition}\label{prop.rate2}
Under the given conditions, $E^{\ast}\|\mathbb{G}_n\|_{{\cal M}_\eta}\leq c_5\eta$, for some $0<c_5<\infty$ and all $\eta$ small enough. 
\end{proposition}

The proofs of Propositions \ref{prop.rate1} and \ref{prop.rate2} are in Section \ref{asec:rateconv} of the Supplementary Material. As all of the needed conditions are satisfied, the proof of Theorem \ref{theorem:rateofconvergence1} is complete.
\end{proof}

\section{Weak convergence}\label{sec:WeakConvergence}
We now derive the limiting distribution of $\left[\sqrt{n}(\hat{\zeta}_n-\zeta_0),\right.$ $\left. n(\hat{\phi}_n-\phi_0)\right]$ and $\left[\sqrt{n}(\check{\zeta}_n-\zeta_0),n(\check{\phi}_n-\phi_0)\right]$. The basic approach is to evaluate a suitably standardized minimized process and verify that it converges to a limiting process on compact sets in a manner which ensures that the mean-midpoint and mode-midpoint of the $\arg\min$s of the processes also converge on compact sets. Then we use this to derive the limiting distribution. However, there are a number of non-standard convergence issues which need to be addressed, including the fact that the uniform norm-based Argmax Theorem---see, e.g., Theorem 14.1 of \cite{kosorok2008introduction}---is not directly applicable because of the presence of discontinuous indicator functions in the limiting objective function. We highlight these issues as they arise. 

To begin, let $H_n=\re^d\times\re^d\times n(S^{p-1}-\omega_0)\times n([l_0,u_0]-\gamma_0)$, let $h=(h_1,h_2,h_3,h_4)\in H_n$ be an index, and define $h\mapsto Q_n(h)=$
\begin{align}
\label{eq:Qn}
Q_n(h_1,h_2,h_3,h_4)=n\left[M_n(\beta_0+h_1/\sqrt{n},\delta_0+h_2/\sqrt{n},\omega_0+h_3/n,\gamma_0+h_4/n)-M_n(\theta_0)\right].
\end{align}
Let $\tilde{h}_n\in\arg\min_{h\in H_n} Q_n(h)$, 
and note, by construction, that
\[\tilde{h}_n=\left(\sqrt{n}(\tilde{\beta}_n-\beta_0),\sqrt{n}(\tilde{\delta}_n-\delta_0),n(\tilde{\omega}_n-\omega_0),n(\tilde{\gamma}_n-\gamma_0)\right),\]
where the last two parameter estimates are contained in a level set $n(\Phi_n-\phi_0)$, which we will denote as $\Phi_n'$. By 
Theorem~\ref{theorem:rateofconvergence1}, $\tilde{h}_n$ and the entirety of $\Phi_n'$ are asymptotically bounded in probability. Also define $\hat{h}_n$ and $\check{h}_n$ by replacing $\tilde{\phi}_n$ in $\tilde{h}_n$ with $\hat{\phi}_n$ and $\check{\phi}_n$, respectively.

We now characterize the limiting distribution and work towards stating the main result. Let $W_1$ and $W_2$ be independent mean zero Gaussian random variables with respective covariances $\Sigma_1=\sigma^2\left(E\left[ZZ'\ind\{U\leq 0\}\right]\right)^{-1}$ and $\Sigma_2=\sigma^2\left(E\left[ZZ'\ind\{U> 0\}\right]\right)^{-1}$. Define a new random process indexed by $g=(g_1,g_2)\in\re^{p-1}\times\re$:
\begin{align}
\label{eq:Q02}
Q_{02}(g)=\sum_{j=1}^{\infty}\ind\{-g_1'\tilde{X}_j^{-}+g_2<-\tilde{U}_j^{-}\leq 0\}\tilde{E}_j^{-}+\ind\{0<\tilde{U}_j^{+}\leq-g_1'\tilde{X}_j^{+}+g_2\}\tilde{E}_j^{+},
\end{align}
where 
\begin{itemize}
\item $\tilde{E}_j^{-}=\left[(\beta_0-\delta_0)'Z_j^{-}\right]^2+2\epsilon_j^{-}(\beta_0-\delta_0)'Z_j^{-}$, $j\geq 1$; 
\item $(\epsilon_j^{-},\;j\geq 1)$ are i.i.d. realizations of the residual $\epsilon$ from the model~(\ref{eq:dmodel}); 
\item $\tilde{X}_j^{-}=\overline{\omega}_0'X_j^{-}$, $j\geq 1$ (recall that $\overline{\omega}_0$ is the $\omega_0-$orthonormal matrix defined in Section \ref{sec:dmod});
\item $(Z_j^{-},X_j^{-})\in\re^d\times\re^p$, $j\geq 1$, are i.i.d. draws from the distribution $G$ of Condition C3; 
\item $\tilde{U}_j^{-}=\sum_{k=1}^j M_k$, $j\geq 1$, where the $(M_k,k\geq 1)$ are i.i.d. exponential random variables with mean $f_0^{-1}$, and $f_0$ is from C2; 
\item $\tilde{E}_j^{+}=\left[(\beta_0-\delta_0)'Z_j^{+}\right]^2-2\epsilon_j^{+}(\beta_0-\delta_0)'Z_j^{+}$, $j\geq 1$; 
\item $(\epsilon_j^{+},\;j\geq 1)$ is an independent random replication of the sequence $(\epsilon_j^{-},\;j\geq 1)$; 
\item $\tilde{X}_j^{+}=\overline{\omega}_0'X_j^{+}$, $j\geq 1$;
\item $(Z_j^{+},X_j^{+},\;j\geq 1)$ is an independent random replication of $(Z_j^{-},X_j^{-},\;j\geq 1)$; and 
\item $(\tilde{U}_j^{+},\; j\geq 1)$ is an independent random realization of $(\tilde{U}_j^{-},\; j\geq 1)$.
\end{itemize}
We note that $\tilde{U}_1^{-}$, $\tilde{U}_2^{-}$, $\ldots$, are equivalent to the jump points of a stationary Poisson process $t\mapsto N^{-}(t)$ over $0< t<\infty$ with intensity $f_0$. 
In other words, for any $0<t<\infty$, $L(t)\equiv\sum_{j=1}^{\infty}\ind\{\tilde{U}_j^{-}\leq t\}$ has the same distribution as $N(t)$; and, conditional on $N(t)=m$, the collection of jump locations of $N(t)$ has the same distribution as $\{\tilde{U}_j^{-}:\; 1\leq j\leq m\}$ conditional on $L(t)=m$. 
Moreover, this collection of jump locations has the same distribution as an i.i.d. sample of $m$ uniform$[0,t]$ random variables. These are all direct consequences of the properties of a stationary Poisson process \citep{ross1995}. 
These same properties hold true for $\tilde{U}_1^{+}$, $\tilde{U}_2^{+}$, $\ldots$, for another independent Poisson process $t\mapsto N^{+}(t)$ having the same distribution as $N^{-}(t)$. 

Now let $\tilde{\Phi}_0=\arg\min_{g\in\re^{p-1}\times\re}Q_{02}(g)$, where we note that this need not be a single number but rather a level set. We will use $\tilde{g}$ to denote some element of $\tilde{\Phi}_0$. For any $g\in\re^{p-1}\times\re$, define $V_j^{-}(g)=\ind\{g_1'\tilde{X}_j^{-}-g_2-\tilde{U}_j^{-}>0\}-\ind\{g_1'\tilde{X}_j^{-}-g_2-\tilde{U}_j^{-}\leq 0\}$ and $V_j^{+}(g)=\ind\{g_1'\tilde{X}_j^{+}-g_2+\tilde{U}_j^{+}>0\}-\ind\{g_1'\tilde{X}_j^{+}-g_2+\tilde{U}_j^{+}\leq 0\}$, for all $j\geq 1$. The following theorem establishes some important properties regarding the elements of $\tilde{\Phi}_0$, the proof of which is given in Section \ref{asec:WeakConvergence_compact} of the Supplementary Material:
\begin{theorem}\label{argmax.compactness}
Under {Conditions C1, C2, C3, C4 and C5}, the following are true:
\begin{enumerate}
\item $T_0\equiv\sup_{g\in\tilde{\Phi}_0}\|g_1\|\vee |g_2|$ is a bounded random variable.
\item $V_j^{-}(g)=V_j^{-}(\tilde{g})$ and $V_j^{+}(g)=V_j^{+}(\tilde{g})$ for all $j\geq 1$ and for any $g,\tilde{g}\in\tilde{\Phi}_0$ almost surely. Accordingly, fix $\tilde{g}\in\tilde{\Phi}_0$, and define $\tilde{V}_j^{-}=V_j^{-}(\tilde{g})$ and $\tilde{V}_j^{+}=V_j^{+}(\tilde{g})$ for all $j\geq 1$.
\item For all $j\geq \tilde{J}_0^{-}\equiv\min\{j\geq 1:\tilde{U}_j^{-}>T_0(k_1+1)\}$, $\tilde{V}_j^{-}=-1$; and, for all $j\geq\tilde{J}_0^{+}\equiv\min\{j\geq 1:\tilde{U}_j^{+}>T_0(k_1+1)\}$, $\tilde{V}_j^{+}=1$. 
\end{enumerate}
\end{theorem}
Result 1 of Theorem~\ref{argmax.compactness} facilitates verifying compactness of our limiting distributions, while Results 2 and 3 allow us to characterize the level set $\tilde{\Phi}_0$ using only a finite portion of the $V_j^{-}$ and $V_j^{+}$ sequences. We now characterize the two main limiting distributions. First, we need to define the following maps from $\re^{p-1}$ to $\re$:
\begin{eqnarray*}
g_1\mapsto\tilde{C}_L(g_1)&=&\left(\max_{1\leq j\leq\tilde{J}_0^{-}:V_j^{-}=-1}g_1'\tilde{X}_j^{-}-\tilde{U}_j^{-}\right)\vee\left(\max_{1\leq j\leq\tilde{J}_0^{+}:V_j^{+}=-1}g_1'\tilde{X}_j^{+}+\tilde{U}_j^{+}\right),\\
g_1\mapsto\tilde{C}_U(g_1)&=&\left(\min_{1\leq j\leq\tilde{J}_0^{-}:V_j^{-}=1}g_1'\tilde{X}_j^{-}-\tilde{U}_j^{-}\right)\wedge\left(\min_{1\leq j\leq\tilde{J}_0^{+}:V_j^{+}=1}g_1'\tilde{X}_j^{+}+\tilde{U}_j^{+}\right),\;\mbox{and}\\
g_1\mapsto\tilde{C}_R(g_1)&=&\tilde{C}_U(g_1)-\tilde{C}_L(g_1),
\end{eqnarray*}
where we take the maximum of an empty set of numbers to be $-\infty$ and the minimum of an empty set to be $+\infty$. By the definitions of $\tilde{J}_0^{-}$ and $\tilde{J}_0^{+}$, we know that $\left(\min_{1\leq j\leq\tilde{J}_0^{-}}V_j^{-}\right)\wedge\left(\min_{1\leq j\leq\tilde{J}_0^{+}}V_j^{+}\right)=-1$ and $\left(\max_{1\leq j\leq\tilde{J}_0^{-}}V_j^{-}\right)\vee\left(\max_{1\leq j\leq\tilde{J}_0^{+}}V_j^{+}\right)=1$, and thus the above functions are well-defined. Let $\tilde{R}_0\equiv\{g_1\in\re^{p-1}:\tilde{C}_R(g_1)>0\}$. Analogously to Lemma~\ref{meanmidexistence}, we provide several useful properties of these limiting quantities, such as convexity and boundedness, of these functions and the set $\tilde{R}_0$ in Lemma \ref{limiting.properties} in the Supplementary Material.


Our main weak convergence results are contained in the next two theorems:
\begin{theorem}\label{convergence-meanmid}{(Weak convergence of the estimators using the mean-midpoint)}
Assume Conditions {C1, C2, C3, C4 and C5} hold. Then
\[\left(\begin{array}{c} \sqrt{n}(\hat{\beta}_n-\beta_0)\\ \sqrt{n}(\hat{\delta}_n-\delta_0)\\ n(\hat{\omega}_n-\omega_0)\\ n(\hat{\gamma}_n-\gamma_0)\end{array}\right)\leadsto\left(\begin{array}{c}W_1\\ W_2\\ \overline{\omega}_0\hat{g}_1\\ \hat{g}_2\end{array}\right),\]
where $W_1$, $W_2$, $\hat{g}=(\hat{g}_1,\hat{g}_2)$ are mutually independent, 
\[\hat{g}_1=\frac{\int_{\tilde{R}_0}g_1\tilde{C}_R(g_1)d\mu(g_1)}{\int_{\tilde{R}_0}\tilde{C}_R(g_1)d\mu(g_1)},\]
$\mu$ is Lebesgue measure on $\re^{p-1}$, and $\hat{g}_2=\left[\tilde{C}_L(\hat{g}_1)+\tilde{C}_U(\hat{g}_1)\right]/2$.
\end{theorem}
\begin{theorem}\label{convergence-modemid}{(Weak convergence of the estimators using the mode-midpoint)}
Assume Conditions {C1, C2, C3', C4 and C5} hold. Then
\[\left(\begin{array}{c} \sqrt{n}(\check{\beta}_n-\beta_0)\\ \sqrt{n}(\check{\delta}_n-\delta_0)\\ n(\check{\omega}_n-\omega_0)\\ n(\check{\gamma}_n-\gamma_0)\end{array}\right)\leadsto\left(\begin{array}{c}W_1\\ W_2\\ \overline{\omega}_0\check{g}_1\\ \check{g}_2\end{array}\right),\]
where $W_1$, $W_2$, $\check{g}=(\check{g}_1,\check{g}_2)$ are mutually independent, and where $\check{g}_1=\arg\max_{g_1\in\tilde{R}_0}\tilde{C}_R(g_1)$ and $\check{g}_2=\left[\tilde{C}_L(\check{g}_1)+\tilde{C}_U(\check{g}_1) \right]/2$.
\end{theorem}

As will be seen in the proofs, the reason the stronger C3' condition is needed in Theorem~\ref{convergence-modemid} is that without the condition, there is no guarantee that the $\arg\max$ over $\tilde{C}_R$ is unique. This non-uniqueness issue, however, is not a concern for Theorem~\ref{convergence-meanmid}. 

\subsection{Some preliminary results for proving Theorems \ref{convergence-meanmid} and \ref{convergence-modemid}}

The main idea of proving Theorems \ref{convergence-meanmid} and \ref{convergence-modemid} includes showing 
i) the uniform convergence of $Q_{2n}$ (a $\phi$-related part of $Q_n$ in (\ref{eq:Qn}) to be defined soon) to $Q_{02}$ in (\ref{eq:Q02})---specifically, the convergence of the boundary of level sets in the two compound Poisson processes, $Q_{2n}$ and $Q_{02}$---and
ii) the convergence of the remaining part, $Q_{1n}$, of $Q_{n}$ to a process of which the argmin is an independent Gaussian random vector.
As is the case for the standard Argmax Theorem, we show the convergence of $Q_{2n}$ to $Q_{02}$ within an arbitrary compact set first (in Theorems \ref{compact.convergence} and \ref{compact.convergence.mode}) and then remove the restriction at the end (Section \ref{sec:WeakConvergence_main}). 
One complication related to defining compact sets is that given a fixed $n$, some quantities may not be well-defined in a small set (although this is not an issue for large $n$)---e.g., argmax may not be well-defined in a summation over a null set. To circumvent this problem, we define a condition---$F_n(k)=1$ in (\ref{eq:Fn})---under which all quantities are unambiguously defined.
Also, because the support of $Q_{2n}$ and $Q_{02}$ are distinct---$n(\mathcal S^{p-1}-\omega_0)$ and  $\mathbb R^{p-1}$---in order for the Argmax approach to make sense, we need to come up with a bijective mapping from one to the other; see Lemma \ref{lemma.onetoone} below.
Thus, in this section, we first decompose the empirical process $Q_{n}$, introduce a one-to-one and onto mapping, and then we define suitable compact index sets for convergence.

Let's examine more closely $Q_n$ and establish an appropriate convergence on compact sets that can facilitate the desired weak convergence. Letting $U_i=\omega_0'X_i-\gamma_0$, we obtain after rearranging the constituent components
\begin{eqnarray}
\label{eq:Qn_decomp}
Q_n(h)&=&\sum_{i=1}^n \ind\{U_i\leq 0\}\left(\frac{(h_1'Z_i)^2}{n}-\frac{2\epsilon_i h_1'Z_i}{\sqrt{n}}\right)
+\sum_{i=1}^n \ind\{U_i>0\}\left(\frac{(h_2'Z_i)^2}{n}-\frac{2\epsilon_i h_2'Z_i}{\sqrt{n}}\right)\nonumber\\
&&+\sum_{i=1}^n\ind\{-h_3'X_i+h_4<n U_i\leq 0\}\nonumber\\
&&\times\left[\left\{\left(\beta_0-\delta_0-\frac{h_2}{\sqrt{n}}\right)'Z_i\right\}^2+2\epsilon_i\left(\beta_0-\delta_0+\frac{h_1-h_2}{\sqrt{n}}\right)'Z_i-\frac{(h_1'Z_i)^2}{n}\right]\nonumber\\
&&+\sum_{i=1}^n\ind\{0<n U_i\leq -h_3'X_i+h_4\}\nonumber\\
&&\times\left[\left\{\left(\delta_0-\beta_0-\frac{h_1}{\sqrt{n}}\right)'Z_i\right\}^2+2\epsilon_i\left(\delta_0-\beta_0+\frac{h_2-h_1}{\sqrt{n}}\right)'Z_i-\frac{(h_2'Z_i)^2}{n}\right]\nonumber\\
&=&Q_{1n}^{-}(h_1)+Q_{1n}^{+}(h_2)+Q_{2n}^{-}(h)+Q_{2n}^{+}(h),
\end{eqnarray}
where we also define 
\begin{eqnarray*}
E_{in}^{-}(h_1,h_2)&=&\left\{\left(\beta_0-\delta_0-\frac{h_2}{\sqrt{n}}\right)'Z_i\right\}^2+2\epsilon_i\left(\beta_0-\delta_0+\frac{h_1-h_2}{\sqrt{n}}\right)'Z_i-\frac{(h_1'Z_i)^2}{n},\;\mbox{and}\\
E_{in}^{+}(h_1,h_2)&=&\left\{\left(\delta_0-\beta_0-\frac{h_1}{\sqrt{n}}\right)'Z_i\right\}^2+2\epsilon_i\left(\delta_0-\beta_0+\frac{h_2-h_1}{\sqrt{n}}\right)'Z_i-\frac{(h_2'Z_i)^2}{n},
\end{eqnarray*}
for $i=1,\ldots,n$, so that $Q_{2n}^{-}(h)=\sum_{i=1}^n\ind\{-h_3'X_i+h_4<n U_i\leq 0\}E_{in}^{-}(h_1,h_2)$ and $Q_{2n}^{+}(h)=\sum_{i=1}^n\ind\{0<n U_i\leq -h_3'X_i+h_4\}E_{in}^{+}(h)$. Also define $Q_{1n}(h_1,h_2)=Q_{1n}^{-}(h_1)+Q_{1n}^{+}(h_2)$ and $Q_{2n}(h)=Q_{2n}^{-}(h)+Q_{2n}^{+}(h)$.

Note that $n(\tilde{\omega}_n-\omega_0)'\omega_0=o_P(1)$ by Lemma~\ref{lthm:cr2} combined with Theorem~\ref{theorem:rateofconvergence1}. Thus, $\tilde{h}_{3n}$ lives in $\re_{\overline{\omega}_0}^{p-1}$ in the limit as $n\rightarrow\infty$; and, moreover, $n(S^{p-1}-\omega_0)$ is a $p-1$ dimensional manifold for all $n\geq 1$. We leverage this structure to create a one-to-one map between $\re^{p-1}$ and $n(S^{p-1}-\omega_0)$ whenever $\|h_3\|\leq\sqrt{2}n$, which, by Theorem~\ref{theorem:rateofconvergence1}, is guaranteed to be almost surely true for all $n$ large enough. Accordingly, define the function $h_{\ast n}:\re^{p-1}\mapsto n(S^{p-1}-\omega_0)$ as
\[h_{\ast n}(g_1)=\overline{\omega}_0 g_1-n\left\{1-\sqrt{\left(1-\frac{\|g_1\|^2}{n^2}\right)_{+}}\right\}\omega_0.\]
The following lemma provides some useful properties of $h_{\ast n}$. The proof is in Section \ref{asec:WeakConvergence_prelim} of the Supplementary Material.
\begin{lemma}\label{lemma.onetoone}
Let $A_n=\{g_1\in\re^{p-1}:\|g_1\|\leq n\}$ and $B_n=\{h_3\in n(S^{p-1}-\omega_0):\|h_3\|\leq\sqrt{2}n\}$. Then $h_{\ast n}:A_n\mapsto B_n$ is continuous, one-to-one and onto, with continuous inverse $h_3\mapsto h_{\ast n}^{-1}(h_3)=\overline{\omega}_0'h_3$, which is also one to one and onto, for all $n\geq 1$.
\end{lemma}

We now use this to define some suitable compact index sets for convergence. For any $k\in(0,\infty)$, define $H_{20}^{\ast}(k)\equiv\{(g_1,g_2)\in\re^{p-1}\times\re:\:\|g_1\|\leq 0.9 k_1^{-1}k,\;|g_2|\leq k\}$ and $H_{2n}(k)\equiv\{(h_3,h_4)\in n(S^{p-1}-\omega_0)\times n([l,u]-\gamma_0):\;\|h_3\|^2\leq (0.9 k_1^{-1}k)^2+r_n^2(k),\;|h_4|\leq k\}$, where
\[r_n(k)=n\left[1-\sqrt{\left(1-\frac{(0.9 k_1^{-1}k)^2}{n^2}\right)_{+}}\right].\]
Note that Lemma~\ref{lemma.onetoone} yields that for a $k\in(0,\infty)$ and for all $n$ large enough, $(h_{\ast n}(g_1),g_2)\in H_{2n}(k)$ for any $(g_1,g_2)\in H_{20}^{\ast}(k)$ and also $(h_{\ast n}^{-1}(h_3),h_4)\in H_{20}^{\ast}(k)$ for any $(h_3,h_4)\in H_{2n}(k)$. Define $H_0^{\ast}(k)=\{(h_1,h_2,g_1,g_2)\in\re^d\times\re^d\times H_{20}^{\ast}(k):\; \|h_1\|\vee\|h_2\|\leq k\}$ and $H_n(k)=\{(h_1,h_2,h_3,h_4)\in\re^d\times\re^d\times H_{2n}(k):\|h_1\|\vee\|h_2\|\leq k\}$, and note that $\lim_{k\rightarrow\infty}H_0^{\ast}(k)=H_0^{\ast}\equiv\re^d\times\re^d\times\re^{p-1}\times\re$.

Let 
\begin{align}
\label{eq:tilde.hn.k}
\tilde{h}_n(k)=(\tilde{h}_{1n}(k),\tilde{h}_{2n}(k),\tilde{h}_{3n}(k),\tilde{h}_{4n}(k))\in\arg\min_{h\in H_n(k)} Q_n(h),
\end{align}
and let $\Phi_n'(k)$ be the level set containing $(\tilde{h}_{3n}(k),\tilde{h}_{4n}(k))$. Let $(\hat{h}_{1n}(k),\hat{h}_{2n}(k))=(\tilde{h}_{1n}(k),\tilde{h}_{2n}(k))$, and define $m_{2n}^{-}(k)=$ \newline $\sum_{i=1}^n \ind\{-4k \leq 
 nU_i<-2k\}$, $m_{2n}^{+}(k)=\sum_{i=1}^n\ind\{2k<nU_i\leq 4k\}$, and 
\begin{align}
     \label{eq:Fn}
    F_n(k)=\ind\{m_{2n}^{-}(k)\wedge m_{2n}^{+}(k)\geq 1\}.
\end{align}
As we will show soon, this indicator is needed to address the potential ill-definition of some quantities (e.g., taking an argmax over a null set); we will examine the convergence of our estimator conditionally on the event $F_n(k)=1$.
Now let $V_i(k)\equiv\ind\{\tilde{h}_{3n}'X_i-\tilde{h}_{4n}+nU_i>0\}-\ind\{\tilde{h}_{3n}'X_i-\tilde{h}_{4n}+nU_i\leq 0\}$, and note that $\max_{1\leq i\leq n}|V_i(k)-V_i|=0$ for $k$ large enough almost surely by Theorem~\ref{theorem:rateofconvergence1}, where the $V_i$'s are as defined in Section~\ref{sec:est}. 

Now, we need to define the following maps from $n(S^{p-1}-\omega_0)$ to $\re$, and an associated subset, for the setting where $F_n(k)=1$: 
\begin{eqnarray}
\label{eq:Cnk}
    h_3\mapsto C_L^{nk}(h_3)=\max_{1\leq i\leq n:\;V_i(k)=-1}h_3'X_i+nU_i, ~ ~  C_U^{nk}(h_3)=\min_{1\leq i\leq n:\;V_i(k)=1}h_3'X_i+nU_i,
\end{eqnarray}
$C_R^{nk}(h_3)= C_U^{nk}(h_3)-C_L^{nk}(h_3)$, and $h_3\mapsto R_n(k)=\{h_3\in n(S^{p-1}-\omega_0):\; C_R^{nk}(h_3)>0,\; \|h_3\|^2\leq (0.9 k_1^{-1}k)^2+r_n^2(k)\}$.
Note that by construction, when $F_n(k)=1$, $\tilde{h}_{3n}\in R_n(k)$ and thus $R_n(k)$ is nonempty. Note that, also by construction, any $h_3\in n(S^{p-1}-\omega_0)$ satisfying $\|h_3\|^2\leq (0.9 k_1^{-1}k)^2+r_n^2(k)$ also satisfies $\|h_3\|\leq k_1^{-1}k$ for all $n$ large enough as $\lim_{n\rightarrow\infty}r_n(k)=0$. Hence, for any $(h_3,h_4)\in H_{2n}$, $|h_3'X_i-h_4|\leq 2k$ for all $n$ large enough. Thus when $F_n(k)=1$ and $n$ is large enough, we know that there is a $nU_i$ for some $1\leq i\leq n$ such that $-4k\leq nU_i<-2k$ and hence $h_3'X_i-h_4+nU_i\leq 0\;\forall(h_3,h_4)\in H_{2n}(k)$ and there also exists another $nU_i$ for some $1\leq i\leq n$ for which $2k<nU_i\leq 4k$ and $h_3'X_i-h_4+nU_i\geq 0;\forall (h_3,h_4)\in H_{2n}(k)$. This guarantees that the maximums and minimums used in the construction of $C_L^{nk}$ and $C_U^{nk}$ are over nonempty sets and are thus well defined.

When $F_n(k)=1$ and $n$ is large enough for the conclusions of the previous paragraph to hold, define
\[\hat{h}_{3n}(k)=\frac{\int_{R_n(k)} h C_R^{nk}(h)d\nu_n(h)}{\int_{R_n(k)} C_R^{nk}(h)d\nu_n(h)},\]
where $\nu_n$ is the uniform measure on $n(S^{p-1}-\omega_0)$,
and $\hat{h}_{4n}(k)=\left[C_L^{nk}(\hat{h}_{3n}(k))+C_U^{nk}(\hat{h}_{3n}(k))\right]/2$. Note that $\hat{h}_{3n}(k)\in R_n(k)$ by arguments similar to those used in defining $\hat{\omega}_n$ in Section~\ref{sec:est}. Also, by construction, both $|C_L^{nk}(\hat{h}_{3n}(k))|\leq k$ and $|C_U^{nk}(\hat{h}_{3n}(k))|\leq k$, and thus $(\hat{h}_{3n}(k),\hat{h}_{4n}(k))\in H_{2n}(k)$. Hence we can define our mean-midpoint estimator $\hat{h}_n(k)\equiv(\hat{h}_{1n}(k),\hat{h}_{2n}(k),\hat{h}_{3n}(k),\hat{h}_{4n}(k))$. Moreover, by recycling previous arguments, it is easy to verify that $n(\hat{\theta}_n-\theta_0)=\hat{h}_n(k)$ for all $k$ large enough almost surely.

We now define the relevant restricted limiting estimators. Let $W_1$, $W_2$ and $Q_{02}$ be mutually independent and as defined above. Recall that $H_{20}^{\ast}(k)=\{(g_1,g_2)\in\re^{p-1}\times\re:\;\|g_1\|\leq 0.9k_1^{-1}k,\;|g_2|\leq k\}$ and $H_0^{\ast}(k)=\{(h_1,h_2,g_1,g_2)\in \re^{d}\times\re^{d}\times H_{20}^{\ast}(k):\;\|h_1\|\vee\|h_2\|\leq k\}$. 
Define $\tilde{h}_{j0}(k)=\arg\min_{h\in\re^d:\;\|h\|\leq k}$ $h'\Sigma_j h-2 h'\Sigma_j W_j$, $j=1,2$; 
$\tilde{g}(k)\in\arg\min_{g\in H_{20}^{\ast}(k)} Q_{02}(g)$; 
$(\tilde{g}_{1}(k),\tilde{g}_{2}(k))\equiv\tilde{g}(k)$; 
and $\tilde{h}_0(k)=(\tilde{h}_{10}(k),\tilde{h}_{20}(k),\overline{\omega}_0\tilde{g}_{1}(k),\tilde{g}_{2}(k))$. Let $(\hat{h}_{10}(k),\hat{h}_{20}(k))=(\tilde{h}_{10}(k),\tilde{h}_{20}(k))$, 
and define $m_{20}^{-}(k)=\sum_{j=1}^{\infty}$ $\ind\{2k<\tilde{U}_j^{-}\leq 4k\}$, $m_{20}^{+}(k)=\sum_{j=1}^{\infty}\ind\{2k<\tilde{U}_j^{+}\leq 4k\}$, and
\begin{align*}
    F_0(k)=\ind\{m_{20}^{-}(k)\wedge m_{20}^{+}(k)\geq 1\}.
\end{align*}

Let $\tilde{V}_{j}^{-}(k)\equiv V_j^{-}(\tilde{g}_0(k))$ and $\tilde{V}_{j}^{+}(k)\equiv V_j^{+}(\tilde{g}_0(k))$, for all $j\geq 1$, and note that $\max_{j\geq 1} |\tilde{V}_{j}^{-}(k)-V_j^{-}(\tilde{g})|=0$ and $\max_{j\geq 1} |\tilde{V}_{j}^{+}(k)-V_j^{+}(\tilde{g})|=0$ for $k$ large enough almost surely by Theorem~\ref{argmax.compactness}. Let $\tilde{\Phi}_0(k)$ be the level set in $\re^{p-1}\times\re$ which contains $\tilde{g}_0(k)$, i.e., $\tilde{\Phi}_0(k)$ is the set of all $g\in\re^{p-1}\times\re$ such that $V_{j}^{-}(g)=\tilde{V}_{j}^{-}(k)$ and $V_{j}^{+}(g)=\tilde{V}_{j}^{+}(k)$, for all $j\geq 1$. Let $\tilde{J}_0^{-}(k)=\min\{j\geq 1:\; \tilde{U}_j^{-}\geq 4k\}$ and $\tilde{J}_0^{+}(k)=\min\{j\geq 1:\; \tilde{U}_j^{+}\geq 4k\}$, and note that both $\tilde{J}_0^{-}(k)$ and $\tilde{J}_0^{+}(k)$ are finite and well-defined. We now need to define the following maps from $\re^{p-1}\mapsto\re$, for the setting where $F_0(k)=1$:
\begin{eqnarray*}
g_1\mapsto\tilde{C}_L^k(g_1)&=&\left(\max_{1\leq j\leq\tilde{J}_0^{-}(k):\tilde{V}_{j}^{-}(k)=-1}g_1'\tilde{X}_j^{-}-\tilde{U}_j^{-}\right)\vee\left(\max_{1\leq j\leq\tilde{J}_0^{+}(k):\tilde{V}_{j}^{+}(k)=-1}g_1'\tilde{X}_j^{+}+\tilde{U}_j^{+}\right),\\
g_1\mapsto\tilde{C}_U^k(g_1)&=&\left(\min_{1\leq j\leq\tilde{J}_0^{-}(k):\tilde{V}_{j}^{-}(k)=1}g_1'\tilde{X}_j^{-}-\tilde{U}_j^{-}\right)\wedge\left(\min_{1\leq j\leq\tilde{J}_0^{+}(k):\tilde{V}_j^{+}(k)=1}g_1'\tilde{X}_j^{+}+\tilde{U}_j^{+}\right),
\end{eqnarray*}
and $g_1\mapsto\tilde{C}_R^k(g_1)=\tilde{C}_U^k(g_1)-\tilde{C}_L^k(g_1)$. Let $\tilde{R}_0(k)=\{g_1\in\re^{p-1}:\;\tilde{C}_R^k(g_1)>0,\;\|g_1\|\leq 0.9 k_1^{-1}k\}$. Note that when $F_0(k)=1$, $\tilde{g}_1(k)\in\tilde{R}_0(k)$ and thus $\tilde{R}_0(k)$ is nonempty. Note also that for any $(g_1,g_2)\in H_{20}^{\ast}(k)$, $|g_1'\tilde{X}_j^{-}-g_2|\leq 2k$ for all $1\leq j\leq \tilde{J}_0^{-}(k)$ and $|g_1'\tilde{X}_j^{+}-g_2|\leq 2k$ for all $1\leq j\leq \tilde{J}_0^{+}(k)$. Thus when $F_0(k)=1$, we know that there is a $\tilde{U}_j^{-}\in(2k,4k]$, for some $1\leq j\leq \tilde{J}_0^{-}(k)$, and thus also $g_1'\tilde{X}_j^{-}-g_2-\tilde{U}_j^{-}\leq 0$ for all $(g_1,g_2)\in H_{20}^{\ast}(k)$. We also know that there is a $\tilde{U}_j^{+}\in(2k,4k]$, for some $1\leq j\leq \tilde{J}_0^{+}(k)$, and hence also $g_1'\tilde{X}_j^{+}-g_2+\tilde{U}_j^{+}> 0$ for all $(g_2,g_2)\in H_{20}^{\ast}(k)$. This means that the minimums and maximums used in the construction of $\tilde{C}_L^k$ and $\tilde{C}_U^k$ are well defined and finite.

When $F_0(k)=1$, define
\[\hat{g}_1(k)=\frac{\int_{\tilde{R}_0(k)} g \tilde{C}_R^{k}(g)d\mu(g)}{\int_{\tilde{R}_0(k)} \tilde{C}_R^{k}(g)d\mu(g)},\]
and $\hat{g}_2(k)=\left[\tilde{C}_L^{k}(\hat{g}_{1}(k))+\tilde{C}_U^{k}(\hat{g}_{1}(k))\right]/2$. Note that $\hat{g}_1(k)\in \tilde{R}_0(k)$ by arguments similar to those used previously. Also, by construction, both $|\tilde{C}_L^{k}(\hat{g}_{1}(k))|\leq k$ and $|\tilde{C}_U^{k}(\hat{g}_{1}(k))|\leq k$, and thus $(\hat{g}_{1}(k),\hat{g}_{2})\in H_{20}^\ast(k)$. Hence we can define our mean-midpoint limiting estimator as $\hat{h}_0(k)\equiv (\tilde{h}_{10}(k),\tilde{h}_{20}(k),$ $\overline{\omega}_0\hat{g}_1(k),\hat{g}_{2}(k))$. Moreover, by recycling previous arguments, it is easy to verify that $\hat{h}_0(k)=(W_1,W_2,\overline{\omega}_0\hat{g}_1,\hat{g}_2)$ for all $k$ large enough almost surely. 

Now we are ready for the following convergence theorem for restrictions over compact sets:
\begin{theorem}\label{compact.convergence}
Assume conditions {C1, C2, C3, C4 and C5}. Then, for every $k\in (0,\infty)$,
\[\left(\begin{array}{c}F_n(k)\hat{h}_n(k)\\ F_n(k)\end{array}\right)\leadsto \left(\begin{array}{c}F_0(k)\hat{h}_0(k)\\ F_0(k)\end{array}\right).\]
\end{theorem}
With this theorem, we have enough results to prove Theorem~\ref{convergence-meanmid}. We postpone the relatively lengthy proof of Theorem \ref{compact.convergence} to the later Section \ref{asec:WeakConvergence_compact_conv}. In Theorem \ref{compact.convergence.mode} in Section \ref{asec:compact.convergence.mode} of the Supplementary Material, we provide a result for the compact convergence of the mode-midpoint estimator in parallel with Theorem \ref{compact.convergence}, which again gives sufficient grounds for proving Theorem \ref{convergence-modemid}.

As a side note, the decomposition (\ref{eq:Qn_decomp}) is helpful in understanding the estimation quality of the change-plane parameters; while $E_{in}^-(h_1, h_2)$ (or $E_{in}^+(h_1, h_2)$) determines the jump sizes of the compound Poisson process of $Q_{2n}^-$ (or $Q_{2n}^+$), the $\left\{\left(\beta_0-\delta_0\right)'Z_i\right\}^2$ part without $h_2$ in the first term of $E_{in}^-(h_1, h_2)$ provides the needed signal for the estimator, while $h_2$ in the first term and second term---a mean zero noise term---makes the estimation problem harder. Thus, one can first see that when the difference in the two regression coefficients is large at the boundary (i.e., $\Xi(\beta_0,\delta_0)\equiv \lim_{u\searrow0}E\Big[|\beta_0 - \delta_0|'Z\Big| U=u\Big]$ is large), because first quadratic term would dominate the noise term, the estimation quality is improved. Second, when the overall estimation error for $\beta_0$ and $\delta_0$ (represented by $h_1$ and $h_2$) is small, its contribution in the first term is minimized (note that when the sample sizes are unbalanced, the MSE for one of $\beta$ and $\delta$ could be lower, but the overall MSE for $\zeta\equiv(\beta,\delta)$ is not optimal). One can speculate that the optimality is achieved when the sample size is balanced between the two groups; or when $\rho\equiv\sqrt{\Pr(\omega_0'X - \gamma_0)\{1-\Pr(\omega_0'X - \gamma_0)\}}$ is maximized. These quantities, $\Xi(\beta_0,\delta_0)$ and $\rho$, affect the estimation quality in finite sample problems.
This is analogous to "the energy of the jump," $\Xi\sqrt {t^*(1-t^*)}$, in the change-point in the mean problem \citep{verzelen2023optimal}, where $t^*$ is the change point in a unit interval, and $\Xi$ is the change in mean. However, in the limit, as Theorem \ref{theorem:rateofconvergence1} implies, the $h$ terms in $E_{in}^-(h_1, h_2)$ vanishes, and the estimation depends primarily on the size of $\Xi(\beta_0,\delta_0)$ relative to $\text{var}(\epsilon_i)$.

\subsection{Proof of Theorems \ref{convergence-meanmid} and \ref{convergence-modemid}}
\label{sec:WeakConvergence_main}

{\bf Proof of Theorem~\ref{convergence-meanmid}}.
Let $F\subset\re^{2d+p+1}$ be closed. Suppose we can show that  for any $\eta>0$, there exists a $0<k<\infty$ such that $\lim\inf_{n\rightarrow\infty}P(F_n(k)=1)\geq 1-\eta$, $P(F_0(k)=1)\geq 1-\eta$, $\lim\inf_{n\rightarrow\infty}P(F_n(k)\hat{h}_n(k)=D_n(\hat{\theta}_n-\theta_0))\geq 1-\eta$, and $P(F_0(k)\hat{h}_0(k)=\hat{h}_0)\geq 1-\eta$, where $D_n$ is a diagonal matrix made of a vector $(\sqrt n \boldsymbol 1_{2d}^\top, n\boldsymbol 1_{p+1}^\top)^\top$. Thus, fixing an $\eta>0$  and finding a corresponding $k<\infty$ which simultaneously satisfy these criteria, we have
\begin{eqnarray*}
\lefteqn{\lim\sup_{n\rightarrow\infty}P(D_n(\hat{\theta}_n-\theta_0)\in F)}&&\\
&\leq&\lim\sup_{n\rightarrow\infty}P(D_n(\hat{\theta}_n-\theta_0)\in F,\;F_n(k)=1,\;F_n(k)\hat{h}_n(k)=D_n(\hat{\theta}_n-\theta_0))\\
&&+\lim\sup_{n\rightarrow\infty}P(F_n(k)=0)+\lim\sup_{n\rightarrow\infty}P(F_n(k)\hat{h}_n(k)\neq D_n(\hat{\theta}_n-\theta_0))\\
&\leq&\lim\sup_{n\rightarrow\infty}P(F_n(k)\hat{h}_n(k)\in F)+2\eta\\
&\leq&P(F_0(k)\hat{h}_0(k)\in F)+2\eta\\
&\leq&P(\hat{h}_0\in F,\;F_0(k)=1,\;F_0(k)\hat{h}_0(k)=\hat{h}_0) +P(F_0(k)=0)+P(F_0(k)\hat{h}_0(k)\neq \hat{h}_0)+2\eta\\
&\leq&P(\hat{h}_0\in F)+4\eta,
\end{eqnarray*}
which implies that $\lim\sup_{n\rightarrow\infty}P(D_n(\hat{\theta}_n-\theta_0)\in F)\;\leq\;P(\hat{h}_0\in F)$,
because $\eta>0$ is arbitrary. As the closed set $F$ was also arbitrary, the desired conclusion follows from the Portmanteau Theorem for weak convergence (see, e.g., Theorem 7.6 of \cite{kosorok2008introduction}). We next need to verify the required probability bounds for all $\eta>0$. 

We begin by showing that for any $\eta>0$, there exists a $k<\infty$ for which $P(F_0(k)=1)\geq 1-\eta$. Recall that $F_0(k)=\ind\{m_{20}^{-}(k)\wedge m_{20}^{+}(k)\geq 1\}$, where $m_{20}^{-}(k)$ 
is the number of $\tilde{U}_j^{-}$ values, among all $j\geq 1$, in the interval $(2k,4k]$. Because the $\tilde{U}_j^{-}$ are progressive sums of i.i.d. exponentials with mean $f_0^{-1}$, $m_{20}^{-}$ is Poisson distributed with parameter $2kf_0$. We can similarly reason that $m_{20}^{+}$ is an independent copy of the same distribution. Thus $\lim\inf_{k\rightarrow\infty}P(F_0(k)=1)=1$, and we have the desired result for this step. Because Theorem~\ref{compact.convergence} yields that $F_n(k)\rightarrow F_0(k)$ in probability, as $n\rightarrow\infty$, we also have that for any $\eta>0$, there exists a $k<\infty$ for which $\lim\inf_{n\rightarrow\infty}P(F_n(k)=1)\geq 1-\eta$. 

Now let $\Phi_n'$ be the level set containing $n(\hat{\phi}_n-\phi_0)$. As argued above, the entirety of $\Phi_n'$ is simultaneously asymptotically bounded in probability as a result of Theorem~\ref{theorem:rateofconvergence1}. Hence for any $\eta>0$, there exists a $k<\infty$ such that both $\lim\inf_{n\rightarrow\infty}P(F_n(k)=1)\geq 1-\eta/2$ and $\lim\inf_{n\rightarrow\infty}P(D_n(\hat{\theta}_n-\theta_0)\in H_n(k))\geq 1-\eta/2$,
and thus
$\lim\inf_{n\rightarrow\infty}P(F_n(k)\hat{h}_n(k)=D_n(\hat{\theta}_n-\theta_0))\geq 1-\eta$. We can apply similar arguments for $\hat{h}_0$ via Theorem~\ref{argmax.compactness} to obtain that for any $\eta>0$ there exists a $k<\infty$ such that $P(F_0(k)\hat{h}_0(k)=\hat{h}_0)\geq 1-\eta$, and our proof is therefore complete.$\Box$

The proof of Theorem~\ref{convergence-modemid} proceeds along very similar arguments as the proof of Theorem~\ref{convergence-meanmid} and will borrow from that proof, especially as all the results of that theorem hold in this proof because Condition C3' is stronger than Condition C3. A detailed proof is given in Section \ref{asec:convergence-modemid} of the Supplementary Material.

\subsection{Convergence on compact sets---Proof of Theorems \ref{compact.convergence}}
\label{asec:WeakConvergence_compact_conv}
We prove Theorem~\ref{compact.convergence} by first defining a certain ensemble process restricted to a compact set---soon to be defined as $\mathcal P_n(k)$---and sufficient to construct $Q_{n}$ within the corresponding compact set. Then, we show in Lemma \ref{compact.convergence.concatenated} that it converges uniformly weakly to another ensemble process, $\mathcal P_0(k)$, which is also sufficient to generate $Q_{02}$. We will introduce a suitable uniform metric, $d_{**}$, needed to define the uniform convergence, and take advantage of the almost sure representation of weak convergence to establish weak convergence of the argmax of $Q_{n}$ to $Q_0$.

{\bf Proof of Theorem~\ref{compact.convergence}}. Define $J_n^{-}(k)=\{1\leq i\leq n:\;-4k\leq nU_i\leq 0\}$, $J_n^{+}(k)=\{1\leq i\leq n:\;0< nU_i\leq 4k\}$, $m_{1n}^{-}(k)=\sum_{i=1}^n\ind\{-2k\leq nU_i\leq 0\}$, $m_{1n}^{+}(k)=\sum_{i=1}^n\ind\{0<nU_i\leq 2k\}$, $m_n^{-}(k)=m_{1n}^{-}(k)+m_{2n}^{-}(k)$ and $m_n^{+}(k)=m_{1n}^{+}(k)+m_{2n}^{+}(k)$, where we recall that  $m_{2n}^{-}(k)=$ $\sum_{i=1}^n \ind\{-4k \leq 
 nU_i<-2k\}$, and 4$m_{2n}^{+}(k)=\sum_{i=1}^n\ind\{2k<nU_i\leq 4k\}$.
 Define also $\hat{\Sigma}_{1n}=\{n^{-1}\sum_{i=1}^n\ind\{U_i\leq 0\} Z_iZ_i'\}^{-1}$, $\hat{\Sigma}_{2n}=\{n^{-1}\sum_{i=1}^n \ind\{U_i>0\} Z_iZ_i'\}^{-1}$, $W_{1n}=n^{-1/2}\sum_{i=1}^n\ind\{U_i\leq 0\}\epsilon_iZ_i$, and $W_{2n}=n^{-1/2}\sum_{i=1}^n\ind\{U_i>0\}\epsilon_iZ_i$. Now we define the ensemble process 
\begin{eqnarray*}
{\cal P}_n(k)&=&\left\{m_n^{-}(k),(X_i,Z_i,\epsilon_i,-nU_i):\;i\in J_n^{-}(k), m_{1n}^{-}(k), m_{2n}^{-}(k), \hat{\Sigma}_{1n}, W_{1n};\;\right.\\
&& \left.m_n^{+}(k),(X_i,Z_i,\epsilon_i,nU_i):\;i\in J_n^{+}(k), m_{1n}^{+}(k), m_{2n}^{+}(k), \hat{\Sigma}_{2n}, W_{2n}\right\},
\end{eqnarray*}
where if either $J_n^{-}(k)$ or $J_n^{+}(k)$ are empty, the corresponding list of constituent variables is taken as a null set. Also, when not null, we order the elements $(X_i,Z_i,\epsilon_i,-nU_i)$ for $i\in J_n^{-}(k)$ in ascending order of $|nU_i|$; and we likewise order $(X_i,Z_i,\epsilon_i,nU_i)$ for $i\in J_n^{+}(k)$ in ascending order of $|nU_i|$. Also define $m_{10}^{-}(k)=\sum_{j=1}^{\infty}\ind\{0\leq \tilde{U}_j^{-}\leq 2k\}$, $m_{10}^{+}(k)=\sum_{j=1}^{\infty}\ind\{0<\tilde{U}_j^{+}\leq 2k\}$, $m_0^{-}(k)=m_{10}^{-}(k)+m_{20}^{-}(k)$, and $m_0^{+}(k)=m_{10}^{+}(k)+m_{20}^{+}(k)$, where $m_{20}^{-}(k)$ and $m_{20}^{+}(k)$ are as defined above. We now define the limiting ensemble process
\begin{eqnarray*}
{\cal P}_0(k)&=&\left\{m_0^{-}(k),(X_j^{-},Z_j^{-},\epsilon_j^{-},\tilde{U}_j^{-}):\;1\leq j\leq m_0^{-}(k), m_{10}^{-}(k), m_{20}^{-}(k), \Sigma_1, W_1;\;\right.\\
&& \left.m_0^{+}(k),(X_j^{+},Z_j^{+},\epsilon_j^{+},\tilde{U}_j^{+}):\;1\leq j\leq m_0^{+}(k), m_{10}^{+}(k), m_{20}^{+}(k), \Sigma_2, W_2\right\},
\end{eqnarray*}
where we use null sets as needed if either $m_0^{-}(k)$ or $m_0^{+}(k)$ are zero, and other terms are as defined previously.

Our current goal is to show that ${\cal P}_n(k)\leadsto{\cal P}_0(k)$ with respect to a suitable uniform metric. Let $\mathbb{D}_{\infty}^q$ be the space of infinite sequences $x_1,x_2,\ldots$, with $x_j\in\re^q$ for all $j\geq 1$, and let $\mathbb{Z}^{0+}$ be the set of non-negative integers. Define the space $\mathbb{E}^q$ consisting of the set of elements $(x_0,\{x_j:\;j\geq 1\})\in\mathbb{Z}^{0+}\times\mathbb{D}_{\infty}^q$ such that $x_j=0$ for all $j>x_0$ (and thus consists of all zeros when $x_0=0$). For $x,y\in\mathbb{E}^q$, define the metric $d_{\ast}(x,y)=|x_0-y_0|+\max_{1\leq j\leq x_0\wedge y_0}\|x_j-y_j\|$, where $d_{\ast}(x,y)=0$ when $x_0=y_0=0$. It is easy to verify that $d_{\ast}$ satisfies the triangle inequality and is otherwise a valid metric making $(\mathbb{E}^q,d_{\ast})$ into a complete metric space. The actual space and metric for the processes ${\cal P}_n(k)$ and ${\cal P}_0(k)$ is $\mathbb{E}^q\times\re^{2+d^2+d}\times\mathbb{E}^q\times\re^{2+d^2+d}$, where $q=p+d+2$, and the metric consisting of the appropriate concatenation of two copies of $d_{\ast}$ and the other needed uniform metrics, which concatenated metric we will denote as $d_{\ast\ast}$ and which will be the default uniform metric on this space. We now have the following lemma the proof of which is in Section \ref{asec:proof.compact.convergence.concatenated} of the Supplemental Material:
\begin{lemma}\label{compact.convergence.concatenated}
Under assumptions {C1, C2, C3, C4 and C5}, and for any $0<k<\infty$, ${\cal P}_n(k)\leadsto{\cal P}_0(k)$ uniformly, as $n\rightarrow\infty$.
\end{lemma}

It is easy to verify that the data in ${\cal P}_n(k)$ is 
sufficient to fully generate the processes $Q_{1n}^{-}(h_1)$ and $Q_{1n}^{+}(h_2)$ for all $h_j\in\re^d:\; \|h_j\|\leq k$, $j=1,2$, and also to generate $F_n(k)$. Also, when $(h_3,h_4)\in H_{2n}(k)$, $|-h_3'X_i+h_4|\leq 2k$ for all $n$ large enough and for any $1\leq i\leq n$, and thus $Q_{2n}^{-}(h)$ and $Q_{2n}^{+}(h)$ is also fully generated by the data in ${\cal P}_n(k)$ for all $h\in H_n(k)$, because for all observations $1\leq i\leq n$ such that $i\not\in J_n^{-}(k)\cup J_n^{+}(k)$ we have that $\ind\{-h_3'X_i+h_4<nU_i\leq 0\}\vee\ind\{0<nU_i\leq -h_3'X_i+h_4\}=0$. We can also readily verify that the data in ${\cal P}_0(k)$ is sufficient to fully generate the processes $h_j\mapsto h_j'\Sigma_j h_j-2h_j'\Sigma_j W_j$, for $h_j\in\re^{d}:\;\|h_j\|\leq k$, and $j=1,2$, and also $F_0(k)$ and $Q_{02}(g)$, for all $g\in H_{20}^{\ast}(k)$. Accordingly, provided $F_n(k)=1$, the data in ${\cal P}_n(k)$ generates $\hat{h}_n(k)$ for all $n$ large enough almost surely; and, similarly, provided $F_0(k)=1$, the data in ${\cal P}_0(k)$ generates $\hat{h}_0(k)$.

By Lemma~\ref{compact.convergence.concatenated}, we know by the almost sure representation theorem (see, e.g., Theorem~7.26 of \cite{kosorok2008introduction}) that there exists a new probability space where the marginal distributions of ${\cal P}_n(k)$ and ${\cal P}_0(k)$ are unchanged but that ${\cal P}_n(k)\rightarrow{\cal P}_0(k)$ uniformly outer almost surely, as $n\rightarrow\infty$. Note that for clarity we are not changing notation to reflect the new probability space. If we can show that this outer almost sure convergence implies that 
\[\left(\begin{array}{c}F_n(k)\hat{h}_n(k)\\ F_n(k)\end{array}\right)\rightarrow\left(\begin{array}{c}F_0(k)\hat{h}_0(k)\\ F_0(k)\end{array}\right)\]
outer almost surely, then our proof of Theorem~\ref{compact.convergence} will be complete. Accordingly, assume going forward that we are in this new probability space where ${\cal P}_n(k)\rightarrow{\cal P}_0(k)$ uniformly outer almost surely. It is now easy to verify that $\|\hat{h}_{jn}(k)-\hat{h}_{j0}(k)\|\rightarrow 0$, as $n\rightarrow\infty$, for $j=1,2$. It is also easy to verify that if $x_n\rightarrow x$ on the space of integers, then $x_n=x$ for all $n$ large enough. Accordingly, we have that  $(F_n(k),m_n^{-}(k),m_{1n}^{-}(k),m_{2n}^{-}(k),m_n^{+}(k),m_{1n}^{+}(k),$ $m_{2n}^{+}(k))=(F_0(k),m_0^{-}(k),m_{10}^{-}(k),m_{20}^{-}(k),m_0^{+}(k),$ $m_{10}^{+}(k),m_{20}^{+}(k))$ for all $n$ large enough almost surely. When $F_0(k)=0$, then the results above verify that $0=F_0(k)\hat{h}_0(k)=F_n(k)\hat{h}_n(k)$ for all $n$ large enough, and thus the proof of the theorem is trivial in this case. Accordingly, assume going forward that $F_0(k)=1$. This then assures us that $F_n(k)=1$ for all $n$ large enough. 

Next, we show that the reordered process is sufficient to generate the argmin of $Q_{n}(h)$ and the mean-midpoint $\hat h_{3n}(k)$. Then, once we establish that the level sets in a finite regime converge to the limiting level sets almost surely, we are almost done. Specifically, we show that
\begin{eqnarray*}
\lim_{n\rightarrow\infty}\sup_{g_1\in\re^{p-1}:\;\|g_1\|\leq 0.9k_1^{-1}k}
&&\left|\tilde{C}_L^{nk}\{h_{*n}(g_1)\}-{C}_L^k(g_1)\right|\vee\\ 
&&\left|\tilde{C}_U^{nk}\{h_{*n}(g_1)\}-{C}_U^k(g_1)\right|\vee
\left|\tilde{C}_R^{nk}\{h_{*n}(g_1)\}-{C}_R^k(g_1)\right|=0,
\end{eqnarray*}
where $C_L^{nk}$, $C_U^{nk}$, and $C_R^{nk}$ are from (\ref{eq:Cnk}). Due to space, we relegate the detailed proof of the sufficiency of $\mathcal P_n(k)$ and the associated convergence of level sets to Section \ref{asec:sufficiency_level_sets} in the Supplementary Material.
These results further imply that for every closed subset $L$ of the open (non-empty) set $\tilde{R}_0(k)$, $\lim_{n\rightarrow\infty}\sup_{g_1\in L}$ $\inf_{g_1'\in \tilde{R}_n(k)}\|g_1-g_1'\|=0$, and for every $\eta>0$, 
\[\lim_{n\rightarrow\infty}\;\inf_{g_1\in \left[\left(\tilde{R}_0^{k}\right)^{\eta}\right]^c}\;\inf_{g_1'\in \tilde{R}_n(k)}\|g_1-g_1'\|\geq \eta,\] 
where superscript $\eta$ on a set denotes the $\eta$-open enlargement of that set and a superscript $c$ denotes complement. Combining this with previous arguments, we now have that $\hat{h}_{3n}(k)\rightarrow \overline{\omega}_0 \hat{g}_1$ and $\hat{h}_{4n}(k)\rightarrow \hat{g}_2$, as $n\rightarrow\infty$, and the proof is complete.$\Box$

\section{Inference}
\label{sec:inference}
In this section, we develop a parametric bootstrap approach to inference. Define $\hat{U}_i=\hat{\omega}_n'X_i-\hat{\gamma}_n$ and $\hat{\epsilon}_i=Y_i-\ind\{\hat{U}_i\leq 0\}\hat{\beta}_n'Z_i-\ind\{\hat{U}_i>0\}\hat{\delta}_n'Z_i$, for $i=1,\ldots,n$, and let $\tilde{U}_n=n^{-1}\sum_{i=1}^n\hat{U}_i$, $\tilde{\epsilon}_n=n^{-1}\sum_{i=1}^n\hat{\epsilon}_i$, $\hat{\tau}_n^2=n^{-1}\sum_{i=1}^n(\hat{U}_i-\tilde{U}_n)^2$, and $\hat{\sigma}_n^2=n^{-1}\sum_{i=1}^n(\hat{\epsilon}_i-\tilde{\epsilon}_n)^2$. Define also $\tilde{\Sigma}_{1n}=\hat{\sigma}_n^2\big[n^{-1}\sum_{i=1}^n$ $\ind\{\hat{U}_i\leq 0\}Z_iZ_i'\big]^{-1}$, $\tilde{\Sigma}_{2n}=\hat{\sigma}_n^2\left[n^{-1}\sum_{i=1}^n\ind\{\hat{U}_i> 0\}Z_iZ_i'\right]^{-1}$, and $\tilde{M}_{n}$ which is an estimator of $\overline{\omega}_0$ using Gram-Schmidt orthogonalization (for example) wherein the columns of the $p\times (p-1)$ matrix $\tilde{M}_{n}$ are orthonormal to each other and orthogonal to $\hat{\omega}_n$. We next need to define several density estimators:
\begin{eqnarray*}
\hat{f}_{n0}=\int_{\re}\frac{1}{\hat{\eta}_{n1}}\phi\left(\frac{t}{\hat{\eta}_{n1}}\right)d\hat{F}_{n1}(t)&\;\;\mbox{and}\;\;&\hat{\xi}_{n}(u)=\int_{\re}\frac{1}{\hat{\eta}_{n2}}\phi\left(\frac{u-t}{\hat{\eta}_{n2}}\right)d\hat{F}_{n2}(t),
\end{eqnarray*}
where $\hat{\eta}_{n1}=2\hat{\tau}_n n^{-1/5}$, $\hat{\eta}_{n2}=2\hat{\sigma}_n n^{-1/5}$, $\phi$ is the standard normal density, $\hat{F}_{n1}(t)=n^{-1}\sum_{i=1}^n\ind\{\hat{U}_i\leq t\}$ and $\hat{F}_{n2}(t)=n^{-1}\sum_{i=1}^n\ind\{\hat{\epsilon}_i-\tilde{\epsilon}_n\leq t\}$. Also, let $\tilde r_n$ be a sequence such that $r_n n^{-1/2}\rightarrow\infty$ and $\tilde r_n/n\rightarrow 0$, and define $\hat{t}_n=\sup\left\{t>0:\;\sum_{i=1}^n\ind\{|\hat{U}_i|\leq t\}\leq \tilde r_n\right\}$ and $\tilde{G}_n=\{(X_i,Z_i):\; |\hat{U}_i|\leq\hat{t}_n\}$. 

We are ready to define our parametric bootstrap. As $W_1$, $W_2$, and the process $g\mapsto Q_{02}(g)$ are all independent, we can generate the components independently. Realizations of $W_j$ can be approximately generated by drawing from a mean zero Gaussian vector with covariance $\tilde{\Sigma}_{jn}$, for $j=1,2$. It is easy to verify that such a random variable, conditional on the observed sample, will converge weakly to $W_j$, $j=1,2$, provided we verify that $\hat{\Sigma}_{jn}$ is consistent for the covariance of $W_j$, for $j=1,2$. The main difficulty is generating random realizations of $g\mapsto Q_{02}(g)$, applying the two midpoint estimators, and verifying that they have the correct conditional limiting distribution. Accordingly, define the process
\begin{eqnarray*}
\tilde{Q}_{02}^{\ast}(g)&=&\int_0^{\infty}\left[B(t)\ind\{g_1'\tilde{X}_{\ast}(t)-g_2>t\}\tilde{E}^{-}_{\ast}(t)\right. \left. +(1-B(t))\ind\{-g_1'\tilde{X}_{\ast}(t)+g_2\geq t\}\tilde{E}^{+}_{\ast}(t)\right]d\tilde{N}(t),
\end{eqnarray*}
where $t\mapsto \tilde{N}(t)$ is a homogeneous Poisson process on $[0,\infty)$ with intensity $2\hat{f}_{n0}$, $B(t)$ is a white-noise type Bernoulli random variable with success probability $1/2$ (as previously defined), and where $\tilde{X}_{\ast}(t)$, $\tilde{E}^{-}_{\ast}(t)$ and $\tilde{E}^{+}_{\ast}(t)$ are also white-noise type stochastic processes similar to processes used in the proof of Theorem~\ref{argmax.compactness}. In this setting, the relevant random draws only occur at jump times in $\tilde{N}$. We define $\tilde{E}^{-}_{\ast}(t)$ and $\tilde{E}^{+}_{\ast}(t)$ as follows:
\begin{eqnarray*}
\tilde{E}^{-}_{\ast}(t)&=&\left[(\hat{\beta}_n-\hat{\delta}_n)'\tilde{Z}_{\ast}(t)\right]^2+2\tilde{\epsilon}_{\ast}(t)(\hat{\beta}_n-\hat{\delta}_n)'\tilde{Z}_{\ast}(t),\;\mbox{and}\\
\tilde{E}^{+}_{\ast}(t)&=&\left[(\hat{\beta}_n-\hat{\delta}_n)'\tilde{Z}_{\ast}(t)\right]^2-2\tilde{\epsilon}_{\ast}(t)(\hat{\beta}_n-\hat{\delta}_n)'\tilde{Z}_{\ast}(t),
\end{eqnarray*}
for all $g\in\re^{p-1}\times\re$, where the pair $(\tilde{X}_{\ast}(t),\tilde{Z}_{\ast}(t))=(\tilde{M}_{n}'X,Z)$ for $(X,Z)$ drawn independently and with replacement from $\tilde{G}_n$; and the draws for $\tilde{\epsilon}_{\ast}(t)$ are taken independently from the density $\hat{\xi}_n$. Note that drawing from $\hat{\xi}_n$ is equivalent to first randomly selecting uniformly an $i\in\{1,\ldots,n\}$ and then adding a normal random variable with standard deviation $\hat{\eta}_{n2}$ to $\hat{\epsilon}_i-\tilde{\epsilon}_n$. The reason we subtract off $\tilde{\epsilon}_n$ is to ensure that the generated random draws have mean zero, conditional on the data, for all $n$ almost surely. 

We can now define $\tilde{g}_{\ast}\in\arg\max_{g\in\re^{p-1}\times\re}\tilde{Q}_{02}^{\ast}(g)$. Note that we could also express $\tilde{Q}_{02}^{\ast}$ in its dual form as an infinite sum, in precisely the same way as $Q_{02}$ is the dual of $\tilde{Q}_{02}$, as described in the proof of Theorem~\ref{argmax.compactness}. We will denote this infinite sum dual form of $\tilde{Q}_{02}^{\ast}$ as $Q_{02}^{\ast}$. We do not provide the details but simply mention this to verify that both the mean-midpoint $\arg\min$ $\hat{g}_{\ast}$ and the mode-midpoint $\arg\min$ $\check{g}_{\ast}$, both based on $\tilde{g}_{\ast}$,  $\tilde{Q}_{02}^{\ast}$ and its infinite-sum dual $Q_{02}^{\ast}$, are well defined, using the same approach as was applied to $Q_{02}$ and $\tilde{Q}_{02}$ to define $\hat{g}$ and $\check{g}$. We now have the main result for this section which verifies that the forgoing parametric bootstrap algorithm provides an asymptotically valid approach for conducting inference on $\hat{\theta}_n$ and $\check{\theta}_n$. The proof is relegated to Section \ref{asec:inference} of the Supplementary Material.
\begin{theorem}\label{parametric.bootstrap}
Under assumptions C1, C2', C3, and C4--C6, we have that $\hat{\Sigma}_{jn}\rightarrow\Sigma_j$ in probability, as $n\rightarrow\infty$, for $j=1,2$; and that, conditional on the observed sample data, $(\tilde{M}_{n}\hat{g}_{1\ast},\hat{g}_{2\ast})$ converges weakly to $(\overline{\omega}_0\hat{g}_1,\hat{g}_2)$, with probability going to one as $n\rightarrow\infty$, where $(\hat{g}_{1\ast},\hat{g}_{2\ast})\equiv \hat{g}_{\ast}$. If, moreover, C3 is strengthened to C3', and either $\hat{\phi}_n$ or $\check{\phi}_n$ is used to define $\hat{U}_i$, then we also have that, conditional on the observed sample data, $(\tilde{M}_{\ast}\check{g}_{1\ast},\check{g}_{2\ast})$ converges weakly to $(\overline{\omega}_0\check{g}_1,\check{g}_2)$, with probability going to one as $n\rightarrow\infty$.
\end{theorem}

%

\section{Numerical estimation}
\label{sec:numeric}
\subsection{A uniform-search approach}
\label{sec:simest_uniform}
In this section, we provide numerical estimation approaches. With non-convexity of the objective function in (\ref{eq:mest}), identifying the global minimum requires an exhaustive search. However, an exhaustive search means exponentially many computations, implying impracticality even for moderate sample sizes, and hence, our algorithms contain relaxation. One of them involves a uniform search, and the other is a modification of mixed integer programming (MIP).

Below, we first give a high-level description of the uniform search algorithm. The idea is to evaluate the objective function (\ref{eq:mest}) using sufficiently many, randomly drawn $(\omega,\gamma)$ samples and to recursively narrow down the sampling frame to the neighborhood of the best-known candidate. This algorithm takes advantage of matrix algebraic tricks and relies on trigonometric transformations for drawing uniform random variables on a local region of a hypersphere. For the detailed algorithm, See Section \ref{asec:search} of the Supplementary Material.




The optimization of $M_n(\theta)$ in (\ref{eq:mest}) can be reduced to finding $(\omega, \gamma)$ that maximizes the total regression sum of squares ($\text{SSR}_n$) without explicitly estimating the regression parameters and then finding the optimal $(\beta,\eta)$ given the estimated $(\omega, \gamma)$. The total $\text{SSR}_n$ is simply the sum of the SSRs (sum of regression squares) of each region segregated by the change-plane characterized by $(\omega, \gamma)$; or $\text{SSR}_n(A^-_n(\omega,\gamma)) + \text{SSR}_n(A^+_n(\omega,\gamma))$, where $\text{SSR}_n(A)=\sum_{i\in A}^n Y_iZ_i' (\sum_{i\in A}^n Z_i Z_i')^-\sum_{i\in A}Z_iY_i$ for $A=A^-_n(\omega,\gamma),A^+_n(\omega,\gamma)$, and $A^-_n(\omega,\gamma) = \{i=1,2,...,n:\omega'X_i-\gamma \le 0\}$ and $A^+_n(\omega,\gamma) = \{i=1,2,...,n:\omega'X_i-\gamma > 0\}$. To simplify this using matrix algebra, let $D_n(\omega,\gamma)$ denote a diagonal matrix with elements $\{\ind{\{\omega'X_i - \gamma \le 0\}}:i=1,2,...,n\}$. Let $ \bs Z_{n+}(\omega, \gamma) = D_n(\omega,\gamma)  \bs Z_n$ and $M_{n+}(\omega, \gamma) = \bs Z_{n+}(\omega, \gamma)'\{\bs Z_{n+}(\omega, \gamma)  \bs Z_{n+}(\omega, \gamma)'\}^-$ $\bs Z_{n+}(\omega, \gamma)$, where $\bs Z_n$ is the vector with elements $\{Z_i\}_{i=1}^{n}$ and superscript $-$ denotes the Moore-Penrose matrix inverse. Similarly, define $ \bs Z_{n-}(\omega,\gamma)$ and $M_{n-}(\omega, \gamma)$ using $\{I-D_n(\omega,\gamma)\}$ in place of $D_n(\omega,\gamma)$.
Then, the sample $\text{SSR}_n(\omega, \gamma) =\bs  Y_n' \{M_{n+}D(\omega,\gamma) + M_{n-}D(\omega,\gamma)\}  \bs Y_n$, where $\bs Y_n$ is the vector with elements $\{Y_i\}_{i=1}^n$. We want to find $(\omega, \gamma)$ that maximizes $\text{SSR}_n(\omega, \gamma)$. Because, given an $\omega$ value, there can be $n-1$ unique possible $D_n(\omega,\gamma)$ values, and each value corresponds to $\gamma$ lying between a neighboring pair of the ordered $\{\omega'X_i\}_{i=1}^n$ values (inclusive of the left limit), we only consider the midpoints of the neighboring pairs for $\gamma$ values. Thus, we aim at finding $\omega$ that maximizes $\text{SSR}_n(\omega) \equiv \text{SSR}_n\{\omega, \tilde\gamma(\omega)\}$, where $\tilde\gamma(\omega) = \arg\max_{\gamma\in \frac 1 2 \{(\omega'X)_{(i)} + (\omega'X)_{(i+1)}\}_{i=1}^{n-1}}\text{SSR}_n(\omega, \gamma)$.

Optimization of the $\text{SSR}_n$ function can be done by searching over a sufficiently large set of uniform random values in $\mathcal {S}^{p-1}$. 
This can be done by projecting $n_0$ $p$-variate standard Gaussian random variables onto the sphere. Once we find the optimal $\omega$ value out of the $n_0$ candidates that maximizes the objective function, we narrow the search span to the neighborhood of the incumbent optimal $\omega$ value and find a better $\omega$ value by sampling uniform random values in the neighborhood. This procedure is repeated until no new best objective value has been achieved over the last $N_0$ iterations---e.g., $N_0 = 20$. The pair of ($\tilde\omega$, $\tilde\gamma$) is, hence, obtained, and $(\tilde\beta, \tilde\eta)$ is immediately obtained through least squares estimation for each resulting subgroup.

Although uniform random sampling on a subset of $\mathcal S^{p-1}$ can still be done through importance sampling after the projection, it could be highly inefficient for a narrow subset. To circumvent this problem, uniform random values can be drawn by first making $m_0$ uniform angles within a range $a_N[-\pi/2, \pi/2)$ for each dimension, where $a_N$ characterizing the range becomes narrower from 1 as the iteration goes large---e.g., $a_N = 0.8 a_{N-1}$---and by randomly drawing $M$ out of the $m_0^{p}$ values with probability proportional to their inverse density. These randomly drawn angles are then transformed to points in the $S^{p-1}$ space that has the most updated $\tilde\omega$ as the origin.
As $\omega$ and $-\omega$ give the same $\text{SSR}_n$, the actual support of the random variables made of $p$ angles, each element of which ranges $[-\pi/2, \pi/2)$, exhausts the whole $p$-dimensional hypersphere.
The density of a random variable $\omega$, of which Cartesian coordinate angles $(\varphi_1,...,\varphi_p)$ follow a uniform density over $[-\pi/,\pi/2)$, is given as $d(\omega) = \prod_{j=1}^{p-1} \cos^{p-j} (\varphi_j)$. 

Once $\tilde\theta$ is obtained, the mean- and mode-argmin estimates are obtained using the same uniform random draw technique because the importance sampling can be highly inefficient. When the sample size ($n$) is large and the change-plane covariates are continuous, the size of the level set becomes tiny. Let $\Omega_N = \{\omega_1, \omega_2, ..., \omega_{M_N}\}$ denote a set of $M_N$ uniform random candidates for $\omega$ in the level set such that $\langle \omega, \tilde\omega\rangle \le b_N$ for some sequence of decreasing constants $b_N\le 1$ indexed by the number of iterations $N$.
Then obtain $\{[C_L^n(\omega_j), C_R^n(\omega_j), R_n(\omega_j)]: j = 1,2,...,M_N\}$. This routine is repeated until $\sum_{j=1}^{M_N} \ind \{R_n(\omega_j) > 0\} \ge R_0$ is satisfied for some large integer $R_0 > 0$, the resolution constant. At the final iteration, the estimates are obtained as $\hat\omega = \frac{\sum_{j = 1}^{M_N} \omega_j \{ R_n(\omega_j) \vee 0 \}}{\|\sum_{j = 1}^{M_N} \omega_j\{ R_n(\omega_j) \vee 0 \}\|}$, $\check \omega = \arg\max_{\omega_j\in\Omega_N} R_n(\omega_j)$, $\hat\gamma = 0.5 C_L(\hat\omega) + 0.5 C_R(\hat\omega)$, and  $\check\gamma = 0.5 C_L(\check\omega) + 0.5 C_R(\check\omega)$.

\subsection{Mixed integer programming approaches}
\label{sec:simest_mip}
We describe a mixed integer programming (MIP) approach and its relaxation. The detailed algorithms can be found in Section \ref{asec:mip} of the Supplementary Material.
While the solutions for $\omega$ are infinitely many in general, 
the number of level sets with distinct values of $\bs V_n \equiv (V_1,...,V_n)$ is much smaller;  especially the maximal number of feasible solutions for $\bs V_n$ is $\sum_{i=0}^p{n\choose i}$, where the feasibility of $\bs V_n$ means that there exists $(\omega,\gamma)$ that satisfies $V_i = \ind\{\omega' X_i - \gamma \le 0\}$ for all $i=1,2,...,n$. An MIP, thus, can help improve the efficiency of the search procedure. 
We include $\bs V_n$ as one of the decision variables on top of $\theta\equiv(\omega,\gamma,\beta,\delta)$. One complication in the application of MIP is that the commercial MIP softwares allow up to a quadratic objective function. While our objective function involves either a matrix inverse---$\bs Y_n'D_n(\omega, \gamma) Z_n\{Z_n'D_n(\omega, \gamma) Z_n\}^-Z_n'D_n(\omega, \gamma)Y_n + Y_n'\{I_n-D_n(\omega, \gamma)\}Z_n [Z_n'\{I_n-D_n(\omega, \gamma)\}Z_n]^-Z_n'(I_n-D_n(\omega, \gamma)\}Y_n$---or a cubic function of the decision variables: $(Y_n-\beta'Z_n)'D_n\{(Y_n-\beta'Z_n) + (Y_n-\delta'Z_n)'\{I_n-D_n(\omega, \gamma)\}(Y_n-\beta'Z_n)$. To make the objective function amenable to the commercial software setup, we introduce dummy variables $q_i = \beta_i 'X_i\ind\{V_i=1\}$ and $s_i=\delta_i'X_i\ind\{V_i=-1\}, i=1,2,...,n,$ and denote $\bs q_n = (q_1,...,q_n)'$ and $\bs s_n=(s_1,...,s_n)'$. 
With these decision variables, it can be shown that the minimizer of the following objective function, after the normalization of $\omega$, is the same as that of (\ref{eq:mest}).
$$\text{To minimize } \bs q_n'\bs q_n + \bs s_n'\bs s_n - 2\bs q_n'Y_n - 2\bs s_n' Y_n,$$  
subject to
\begin{align*}
    &V_i (\omega'X_i - \gamma) \le 0, & i=1,2,...,n,\\
    &\omega'\omega > 0.1, -1 \le \omega_j \le 1, &j=1,2,...,d,\\
    &q_i = \beta_i'Z_i \ind\{V_i=1\}, & i=1,2,...,n, \text{ and}\\
    &s_i = \delta_i'Z_i\ind\{V_i=-1\}, & i=1,2,...,n,
\end{align*}
where the second constraint pushes the solution for $\omega$ away from the origin (all zeros), while by not putting a strong constraint such as $\|\omega\|=1$, the numerical search is done with more flexibility. The $(\omega,\gamma)$ part of the final solution is then rescaled so that $\|\omega\|=1$.

Although MIP provides an efficient search framework, it still requires an exponentially large amount of computation as the sample size becomes large. To relax the computation issue, we follow \cite{lee2021factor}, who introduced the block coordinate descent (BCD) algorithm for finding the optimal $\boldsymbol \theta$. The idea is to optimize $(\omega,\gamma)$ after fixing the current best $(\beta, \delta)$ value, and given the new best $(\omega,\gamma)$, the new best  $(\beta, \delta)$ value is obtained. This routine is repeated until no improvement in the objective function is achieved. We have an algebraic formulation for the MIP set-up that is distinct from \cite{lee2021factor}. See 
Section \ref{asec:mip} in the Supplementary Material for the algebraic formulation of the MIP objective and constraints for both the baseline (i.e., no relaxation) and BCD MIP algorithms.

\section{Simulation study}
\label{sec:sim}

We investigate the finite sample performance of the proposed estimation method in Monte Carlo simulations. Throughout this section, the BCD MIP algorithm (Section \ref{sec:simest_mip}, specifically Algorithm \ref{al:mip_cda} in the Supplementary Material) was used. The uniform search algorithm results in, by and large, similar outcomes, and they are provided in Section \ref{asec:sim} of the Supplementary Material. We consider three models representing different dimensions of the change-plane covariates $X$, and for each model, two sets of parameters are given, resulting in six scenarios. See Table \ref{table:scenario}.
\begin{table}[ht]
    \centering
    \begin{tabular}{cccccccc}
    \hline
          & & & &  \multicolumn{2}{c}{Scenario 1} & \multicolumn{2}{c}{Scenario 2}\\
         Model & $X$ & $Z$ & $(\omega, \gamma)$ & $\beta$ & $\delta$ &  $\beta$ & $\delta$\\
     \hline
         Model 1 & $U_1(-2,2)$ &  $\{1,\text{Bern}(0.5)\}$ & $(1,1)$ & $\boldsymbol 1_2$ & $-\boldsymbol 1_2$ &  $\boldsymbol {1.5}_2$ & $\boldsymbol {0.5}_2$ \\
         Model 2 & $\{U_1(-3,3), \text{Bern}(0.5)\}$ &  $\{1,\text{Bern}(0.5)\}$
         & $\frac 1 {\sqrt{2}} (1,-1,1)$ & $\boldsymbol 1_2$ & $-\boldsymbol 1_2$ &  $\boldsymbol {1.5}_2$ & $\boldsymbol {0.5}_2$ \\
         Model 3 & $U_3(-2,2)$ &  $\{1, U_2(-2,2)\}$ & $\frac 1 {\sqrt{3}} (1,-1,-1,1)$ & $\boldsymbol 1_3$ & $-\boldsymbol 1_3$ &  $\boldsymbol {1.5}_3$ & $\boldsymbol {0.5}_3$\\
     \hline
    \end{tabular}
    \caption{Models and parameter specification for the simulations. Model 1 is a change-point model, and Models 2 and 3 are change-plane models with discrete variables and only with continuous covariates, respectively. $U_p$ is the $p$-dimensional uniform random variables over the following basis support, $\text{Bern}$ is a Bernoulli random variable, $\boldsymbol 1_d$ is a vector of $d$ many $1$'s. For all the settings, $\epsilon$ follows the standard normal distribution.}
    \label{table:scenario}
\end{table}
In Scenario 1, the regression coefficients $(\beta,\delta)$ have opposite signs between the two subgroups, while in Scenario 2, the associations share the same direction but have different strengths.
The ratio of subgroup sizes is $2:1$ for Models 1 and 3, and $3:1$ for Model 2, respectively. Sample sizes of $n=$ 125, 250, 500, 1000, and 2000 are used, and each scenario is repeated $n_{\text{rep}}=300$ times.

The first part of the simulations demonstrates the consistency and the rate of convergence of the numerical values of the point estimates. The second part of the simulations is to examine the weak convergence of  $r_n(\hat \theta_n - \theta_0)\equiv\left[\sqrt{n}(\hat{\zeta}_n-\zeta_0),n(\hat{\phi}-\phi_0)\right]$ and $r_n(\check \theta_n - \theta_0)\equiv\left[\sqrt{n}(\check{\zeta}_n-\zeta_0),n(\check{\phi}-\phi_0)\right]$. A sample of size $n_{\text{rep}}$ is drawn from the limiting distribution, where the limiting random variables are numerically obtained using an approach analogous to what is described in Section \ref{sec:est}. The empirical distributions were first compared marginally between the point estimates ($r_n(\hat\theta - \theta_0)$) and a random sample from the limiting distribution. The joint weak convergence is studied using the Cram\'er-Wold device; the marginal comparison approach was repeated for two random linear combinations of the parameters. The third part of the simulations is to examine the empirical validity of the parametric bootstrap procedure. We investigate the coverage probability of the resulting confidence intervals based on the mean-argmins. In the final part of the simulations, we consider a few more scenarios to examine the impact of the imbalance of the change-plane on the estimation quality.

\begin{figure}
    \centering
    \includegraphics[width = \textwidth] {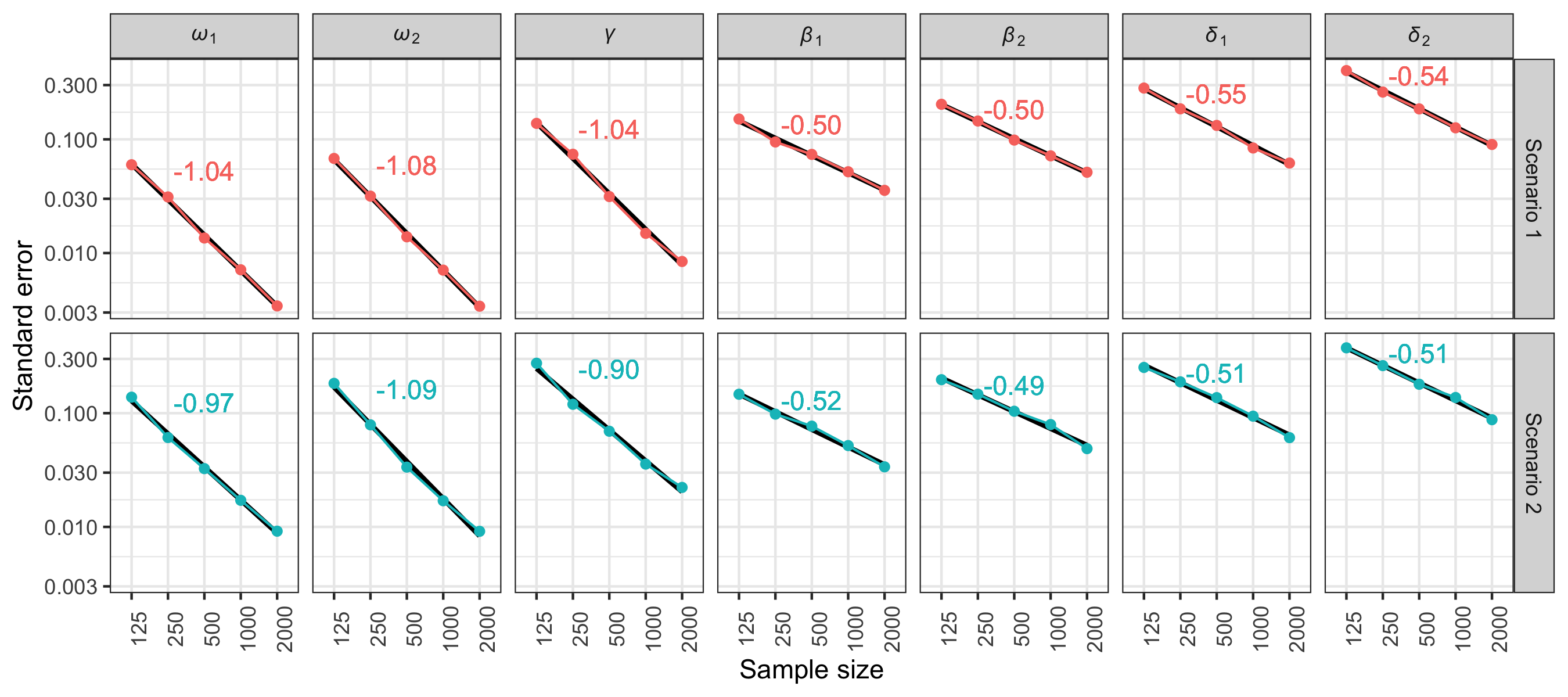}
    \caption{Rate of convergence results for mean-argmin estimators for Model 2. Both standard errors ($y$-axis) and sample sizes ($x$-axis) are presented on the log scale (with base 2). The lines and the annotated numbers are the least squares regression, and the corresponding slope estimates represent the exponents of the rate of convergence. M1, Model 1; M2, Model 2; M3, Model 3}
    \label{fig:rate1}
\end{figure}

Figure \ref{fig:rate1} illustrates that the estimation error of all estimates decreases polynomially with an increase in sample size for Model 2. Similar patterns are observed for the other models. See Section \ref{asec:sim} of the Supplementary Material. Also, the theoretical rate of convergence---$n^{-1}$ for $\hat\phi$ and $\check\phi$ and $n^{-1/2}$ for $\hat\zeta$ and $\check\zeta$---is well observed in the numerical results for both mean- and mode-argmin estimators in most settings.

Scenario 2 has a lower jump size ($\Xi(\beta_0,\delta_0)\approx1.5$) than Scenario 1 ($\Xi(\beta_0,\delta_0)\approx2.9$) for all models. As can be seen from Figure \ref{fig:rate1} (and Figure \ref{fig:rate2} in the Supplementary Material for the other models), the estimation quality for $\phi$ is worse for Scenrio 2 (e.g., SE is higher for Scenario 2). To examine the finite sample estimation quality for imbalanced group sizes, we further consider two more scenarios, which are modifications of Model 2 Scenario 1. Namely, a perfectly balanced scenario and a highly imbalanced scenario. Simulation results indicate that a group imbalance also lowers the estimation quality. See Section \ref{asec:balance} of the Supplementary Material.


In Figure \ref{fig:weak}, the weak convergence simulation results are presented for Model 2, and the Supplementary Materials contain the results for the other models. To see the joint convergence, per the Cram\'er-Wold device, the cumulative distribution function of linear combinations of the estimates with random coefficients were compared to the corresponding limiting distribution. 
The numerical results support the weak convergence theory.

\begin{figure}
    \centering
    \includegraphics[width = \textwidth]{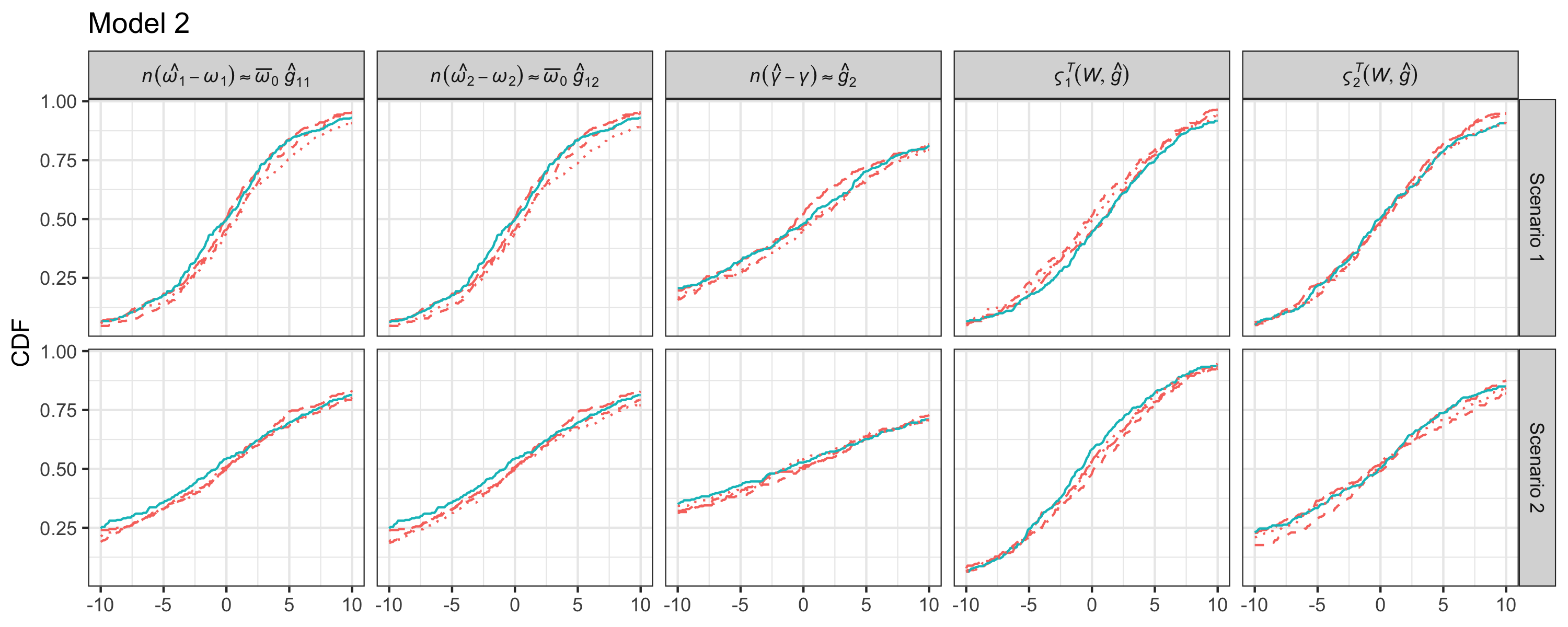}
    \includegraphics[width = \textwidth]{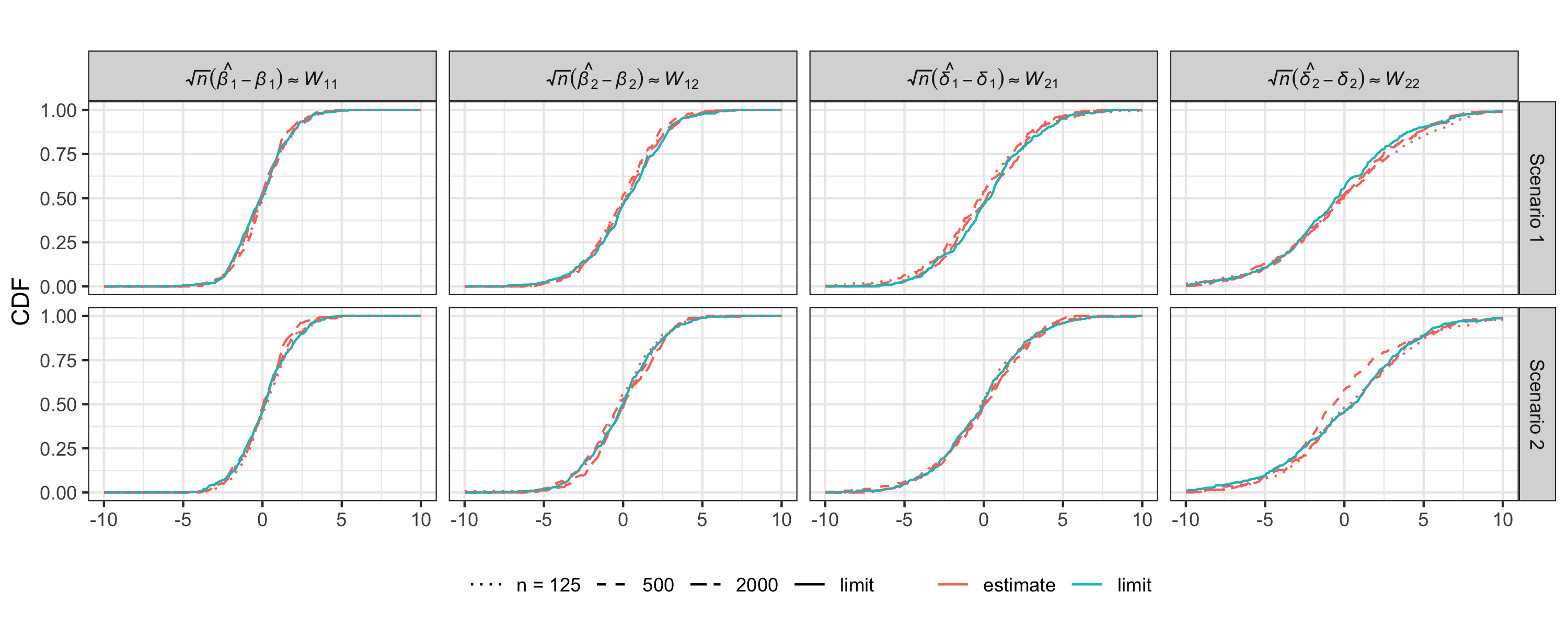}
    \caption{The mean-argmin estimated and limiting CDFs for Model 2. The random coefficients are  $\varsigma_1 = (-0.63, 0.40, 0.15, -0.66, 0.89, 0.89, -0.74)'$ and $\varsigma_2 = ( 0.67, -0.06, 0.10, 0.11, -0.52, 0.52, -0.64)'$ for Model 2.}    
    \label{fig:weak}
\end{figure}

We give the parametric bootstrap simulation results in Table \ref{table:BS}, where we used $\tilde r_n=n^{2/3}$ in the definition of $\hat t_n$ in Section \ref{sec:inference}. For most settings and parameters, the coverage probabilities are close to the nominal confidence level, 95\%. The estimated confidence intervals for $\phi$ are sometimes conservative, but there is a tendency that they become closer to the nominal level for larger sample sizes. This is likely a consequence of either finite sample sizes or the non-exhaustive search algorithm for higher-dimensional parameter spaces. We leave it as future work to develop a numerically more precise algorithm for improved coverage probability. The confidence intervals for the linear combinations ($\zeta_1$ and $\zeta_2$) have good coverage probabilities, which implies that the parametric bootstrapping framework works reasonably well for both marginal and joint distributions. More detailed results, such as confidence interval diagrams (Figure \ref{fig:ci_plot}) and coverage probabilities under alternative parameter values (Table \ref{tab:alternative_power}), are available in Section \ref{asec:sim} of the Supplementary Material.

\begin{table}[ht]
\centering
\begin{tabular}{rrrrrrrrrrrrrrrr}
  \hline
M & S &  $n$ & $\hat\omega_1$ & $\hat\omega_2$ & $\hat\omega_3$ & $\hat\gamma$ & $\hat\beta_1$ & $\hat\beta_2$ & $\hat\beta_3$ & $\hat\delta_1$ & $\hat\delta_2$ & $\hat\delta_3$ & $\hat\varsigma_1$ & $\hat\varsigma_2$ \\ 
  \hline
1 & 1 & 125 &  &  &  & 1.00 & 0.92 & 0.94 &  & 0.95 & 0.92 &  & 0.91 & 0.91 \\ 
  1 & 1 & 500 &  &  &  & 0.99 & 0.93 & 0.95 &  & 0.93 & 0.94 &  & 0.93 & 0.94 \\ 
  1 & 1 & 2000 &  &  &  & 0.98 & 0.96 & 0.94 &  & 0.94 & 0.94 &  & 0.95 & 0.96 \\ 
  1 & 2 & 125 &  &  &  & 0.98 & 0.92 & 0.94 &  & 0.93 & 0.96 &  & 0.95 & 0.93 \\ 
  1 & 2 & 500 &  &  &  & 0.99 & 0.94 & 0.95 &  & 0.94 & 0.94 &  & 0.94 & 0.94 \\ 
  1 & 2 & 2000 &  &  &  & 0.98 & 0.91 & 0.91 &  & 0.95 & 0.95 &  & 0.93 & 0.96 \\ 
  2 & 1 & 125 & 0.98 & 0.98 &  & 0.99 & 0.93 & 0.95 &  & 0.92 & 0.90 &  & 0.93 & 0.91 \\ 
  2 & 1 & 500 & 0.98 & 0.97 &  & 0.98 & 0.93 & 0.96 &  & 0.91 & 0.94 &  & 0.94 & 0.93 \\ 
  2 & 1 & 2000 & 0.97 & 0.98 &  & 0.97 & 0.94 & 0.94 &  & 0.95 & 0.95 &  & 0.93 & 0.95 \\ 
  2 & 2 & 125 & 0.95 & 0.95 &  & 0.94 & 0.93 & 0.96 &  & 0.94 & 0.92 &  & 0.93 & 0.95 \\ 
  2 & 2 & 500 & 0.98 & 0.99 &  & 0.98 & 0.92 & 0.95 &  & 0.93 & 0.94 &  & 0.92 & 0.94 \\ 
  2 & 2 & 2000 & 0.96 & 0.96 &  & 0.94 & 0.96 & 0.96 &  & 0.95 & 0.96 &  & 0.93 & 0.95 \\ 
  3 & 1 & 125 & 0.99 & 0.99 & 0.98 & 0.99 & 0.94 & 0.92 & 0.92 & 0.92 & 0.92 & 0.92 & 0.97 & 0.95 \\ 
  3 & 1 & 500 & 0.97 & 0.97 & 0.98 & 0.99 & 0.94 & 0.95 & 0.94 & 0.95 & 0.94 & 0.94 & 0.94 & 0.93 \\ 
  3 & 1 & 2000 & 0.98 & 0.96 & 0.96 & 0.97 & 0.95 & 0.92 & 0.92 & 0.95 & 0.96 & 0.93 & 0.90 & 0.93 \\ 
  3 & 2 & 125 & 0.98 & 0.98 & 0.98 & 0.98 & 0.92 & 0.89 & 0.93 & 0.90 & 0.93 & 0.92 & 0.97 & 0.95 \\ 
  3 & 2 & 500 & 0.99 & 0.99 & 0.98 & 0.99 & 0.95 & 0.95 & 0.94 & 0.94 & 0.93 & 0.97 & 0.96 & 0.93 \\ 
  3 & 2 & 2000 & 0.95 & 0.95 & 0.96 & 0.97 & 0.95 & 0.94 & 0.94 & 0.97 & 0.93 & 0.91 & 0.95 & 0.95 \\ 
   \hline
\end{tabular}
\caption{The coverage probability of the change-plane mean-argmin estimators based on 300 replicates. M, Model; S, Scenario; $\hat\varsigma_1$ and $\hat\varsigma_2$ are the linear combinations of the estimates with the coefficients previously described.}
\label{table:BS}
\end{table}


\section{Application to the ACTG175 AIDS study}
\label{sec:data_ACTG175}
We apply the change-plane regression model to the AIDS Clinical Trials Group 175 Study (ACTG175) data \citep{hammer1997controlled} to demonstrate the application of the proposed method to a precision medicine problem. The ACTG175 study is a double-blind, controlled, randomized trial that compared the effect of four daily regimes in adults infected with human immunodeficiency virus type 1 (HIV-1) whose CD4 cell counts were from 200 to 500 per cubic millimeter. A lower CD4 cell count is indicative of a weakened immune system.
These data have been used to study individualized treatment recommendations in a few recent studies, including 
 \cite{fan2017change}. In our analysis, we include 1,046 adults who have been randomized into either 600 mg of zidovudine (ZDV) $+$ 400 mg of didanosine(ddI) ($T=1$) or 600 mg of zidovudine $+$2.25 mg of zalcitabine(zal) ($T=0$) as done in
 \cite{fan2017change}. We aim to separate these patients based on their heterogeneous treatment effects and find the best treatment recommendation that maximizes the expected CD4 cell count at $20 \pm 5$ weeks.

We model CD4 cell counts at 20 weeks ($Y$) as a function of treatment ($T$), age ($A$), and homosexual activity ($H$) for two subgroups divided by a hyperplane characterized by the last two covariates: 
\begin{align}
\label{eq:ACTG175_ch2}
Y= & \ind\{\omega'X-\gamma\leq 0\} (\beta_0 +\beta_TT + \beta_AA + \beta_HH) +  \ind\{\omega'X-\gamma>0\} (\delta_0+\delta_TT + \delta_AA + \delta_HH) +\epsilon,
\end{align}
where $X=(A, H)$ and $\omega'X=\omega_A A +\omega_H H$. In Table \ref{tab:actg}, the mean-argmin estimates and the 95\% parametric bootstrap confidence intervals are presented.

In the context of precision medicine, the conditional average treatment effect (CATE) given $X=x$ is defined by $\Delta(x)=E\{Y^{*}(1)-Y^{*}(0)|X=x\}$, where  $Y^{*}(t)$ is the potential outcome for treatment $T=t$. The estimation of $\Delta(x)$ requires the stable unit treatment value assumption (SUTVA), no unmeasured confounders assumption (that is, $T \perp \!\!\! \perp \{Y^{*}(0), Y^{*}(1)\} |X$), and positivity, i.e., $\inf_{t,x} \Pr(T=t|X=x) > 0$  \citep{rubin2005causal, hernan2020}, all of which are believed to be satisfied for the ACTG 175 study. The CATE implied by Model (\ref{eq:ACTG175_ch2}) is given as 
$\Delta(x) =E(Y|T=1, X=x) -E(Y|T=0, X=x)
= \ind\{\omega'x-\gamma \leq 0\} \beta_T + \ind\{\omega'x-\gamma > 0\}\delta_T 
=\beta_T +\ind\{\omega'x-\gamma > 0\}(\delta_T -\beta_T),$
 and its estimate at $x=(A, H)$ is given as  
\begin{equation*}
\hat{\Delta}(A,H)=-6.50 + \ind\{0.077A -0.997H-1.889 > 0\} \times 69.94.
\end{equation*}

\begin{table}[]
    \centering
{\scriptsize
    \begin{tabular}{ccccccccccccc}
         & $\omega_A$ & $\omega_H$ & $\gamma$ & $\beta_0$ & $\beta_T$& $\beta_A$& $\beta_H$ & $\delta_0$ & $\delta_T$& $\delta_A$ & $\delta_H$&  $\delta_T-\beta_T$\\
         \hline
         $\hat\theta$ &
           0.077 &
         -0.997 &
          1.889 &
        409.81 &
         -6.50 &
         -2.02 &
         51.78 &
        369.44 &
         63.44 &
         -0.27 &
        -14.31 &
        69.94 \\
         $\hat\theta_L$ &
          -0.013 &
         -1.004 &
         -0.926 &
        324.97 &
        -32.45 &
         -5.63 &
         -5.24 &
        302.71 &
         38.74 &
         -1.92 &
        -42.03 &
        37.65
         \\
          $\hat\theta_U$ &
            0.148 &
          -0.992 &
           4.084 &
         493.56&
          17.51 &
           1.33 &
         107.48&
         435.13&
          85.29 &
           1.40 &
          14.65 &
          103.09
         \\
         \hline
    \end{tabular}
}
    \caption{The mean-argmin estimates and the 95\% confidence intervals of model (\ref{eq:ACTG175_ch2}) of the ACTG175 study.}
    \label{tab:actg}
\end{table}

$\Delta(A, H) \geq 0$ indicates a greater expected CD4 cell count (a better immune system) under $T=1$ than that under $T=0$ for the patients with $(A, H)$. For 571 (54.6\%) patients in the data who satisfy $0.077 A - 0.997 H -1.889 > 0$, $T=1$ (ZDV $+$ ddI) is recommended over $T=0$ (ZDV $+$ zal) to maximize the CD4 cell count at 20 weeks, with an estimated treatment effect $\hat{\Delta}=\hat{\delta}_T=63.44$. For the other subgroup of 475 patients, $T=0$ is recommended, as $T=1$ compared to $T=0$ has a negative estimated treatment effect, $\hat{\Delta}=\hat{\beta}_T=-6.50$. The difference of treatment effects in two identified subgroups, $\hat\delta_T - \hat\beta_T = 69.94$, is relatively large to the degree that the confidence interval, $(37.7, 103.1)$, does not contain zero ($p < 0.002$). We provide details of the derivation and interpretation of $\Delta(X)$ in Section H of the Supplementary Material.

\section{Discussion}
\label{sec:conclusion}
We studied the asymptotic behavior of change-plane estimators. Although there are quite a few theoretical results developed recently, to the best of our knowledge, ours is the first to precisely characterize the limiting distribution of the prototypical change-plane estimator without a surrogate loss or smoothing. The $n$-rate of convergence of the change-plane parameter estimates was proved and illustrated by simulation. Furthermore, with the limiting distributions being established, valid inferences on the change-plane and regression parameters were developed using a parametric bootstrap. 

There are a number of interesting directions for future work. Primarily among them are extensions of the present results to multiple change-plane models \citep{li2021multithreshold}, and shared regression slope models. Also, because it may make more sense to have a non-linear change-plane in a moderately high-dimensional $p$ setting, one might consider developing a change-manifold model, as one of the referees mentioned, or a kernel change-plane model.
From a computational point of view, it would be of value to develop efficient algorithms for enumerating all hyperplanes which split a given data cloud into two non-empty sets.  Sample size calculation and statistical tests for the existence of a change-plane (in the spirit of \cite{fan2017change}) is another important direction in which this work might be taken. Applying the change-plane regression model to correlated or heteroscedastic data---eg. time-series data as in \cite{enikeeva2019high} and \cite{chan2021optimal}---can create additional complications and also warrants future work.


 \section*{Acknowledgements}
Kang was supported in part by the Central Research Development Fund, University of Pittsburgh. Laber and Kosorok were funded in part by USA NCI grant CA142538. The authors gratefully acknowledge the support from the University of Pittsburgh Center for Research Computing through the resources provided. The authors thank anonymous reviewers who helped us develop the idea of estimation quality and helped explore alternative numerical procedures.

\section*{supplement}
Supplementary Materials available online include a list of notations, the proofs and supporting lemmas in Sections \ref{sec:est}--\ref{sec:inference}, the detailed algorithm introduced in Section \ref{sec:numeric}, the comprehensive simulation results in Section \ref{sec:sim}, and the detailed interpretation of the case study in Section \ref{sec:data_ACTG175}.

\bibliographystyle{imsart-nameyear}
\bibliography{mcp2.bib}       

\pagebreak

\begin{frontmatter}
\title{Supplementary material for ``Inference for change-plane regression" by Kang et al. (2024)}
\runtitle{Supplementary material for Kang et al. (2024)}

\end{frontmatter}


\appendix   

\section{List of notations} 
\label{asec:notation}
\subsection{List of notations}
\label{asec:list_notation}



\begin{longtable}{p{.15\textwidth}   p{.80\textwidth}} 
    \hline
    Section 2 \\
    \hline
         $n$ & Sample size  \\
         $X$ & $p-$ dimensional change-plane covariates \\
         $Z$ & $d-$ dimensional regression covariates \\
         $\epsilon$ & A continuous random variable, with mean zero and variance $0 <\sigma^2 <\infty$; independent of $(X, Z)$\\
         $(\beta_0, \delta_0)$ & $d-$ dimensional regression parameters \\
         $(\omega_0, \gamma_0)$ & $\omega\in S^{p-1}$ and $\gamma_0\in\re$ are the change-plane parameters \\
         $S^{p-1}$ &  $p-1$ dimensional unit sphere embedded in $\re^p$ (i.e., $\left\lbrace \omega \in \re^p\,:\,\|\omega\|=1\right\rbrace$), for $p\geq 2$. \\
         $S^0$ &  $S^0=\{-1,1\}$ \\
          $\zeta$ & $\zeta=(\beta,\delta)$, regression parameters \\
           $\phi$ & $\phi=(\omega,\gamma)$, change-plane parameters \\
          $\theta$ &  $\theta=(\zeta,\phi)$, regression and change-plane parameters \\
          $\ol{\omega}_0$ &  $p\times(p-1)$ matrix consisting of the orthonormal basis vectors
for the $p-1$-dimensional subspace in $\re^p$, 
            which is orthogonal to the linear span of $\omega_0$ \\
      $\re^{p-1}_{\ol{\omega}_0}$ & $p-1$-dimensional subspace orthogonal to $\omega_0$ in $\re^p$\\
      $U$  & $U\equiv\omega_0'X-\gamma_0$; continuous random variable with a continuous density $f$  in a neighborhood of zero  \\
      $f_0$ &  $f_0\equiv f(0)$,  $0<f_0<\infty$. \\
      $G$ &  Limiting distribution of $(Z,X)$ given $U=u$ for which $u\rightarrow 0$ \\
      $\overline{U}$ &  $\overline{U}=\overline{\omega}_0'X$ \\
       $\xi$ & Density of $\epsilon$\\          
       \hline
    Section 3 \\
    \hline
     $M_n(\theta)$& $M_n(\theta) \equiv \ep_n m_{\theta}(Y, Z, X)$, where 
        $m_{\theta}(Y, Z, X) = \ind\{\omega'X-\gamma\leq 0\}(Y-\beta'Z)^2+\ind\{\omega'X-\gamma>0\}(Y-\delta'Z)^2 $ \\
        $\ep_n$ & Standard empirical measure (e.g. $\ep_n X=n^{-1}\sum_{i=1}^{n} X_i$) \\
         $K$ &  The range of $\theta$ ; $K=K_1\times K_2$, where $K_1=\re^d\times\re^d$ and $K_2=S^{p-1}\times[l_0,u_0]$. $\gamma\in[l_0,u_0]$, where $l_0=-k_1-\rho$ 
         and $u_0=k_1+\rho$, for any $\rho>0$, $\|X\|\leq k_1<\infty$ almost surely; $\|Z\|\leq k_2<\infty$ almost surely \\
         $A_0(\phi)$ & $A_0(\phi)=\ind\{\omega'X-\gamma\leq 0\}$ \\
         $A_1(\phi)$ & $A_1(\phi)=\ind\{\omega'X-\gamma>0\}$ \\
         $D_{0n}(\phi)$ &  $D_{0n}(\phi)=\ep_n[ZZ'A_0(\phi)]/[\ep_n A_0(\phi)]$ \\
         $D_{1n}(\theta)$ & $D_{1n}(\theta)=\ep_n[ZZ'A_1(\phi)]/[\ep_n A_1(\phi)]$. \\
         $K_{2n}$ &  The set of $\phi\in K_2$ s.t. the smallest eigenvalues of $D_{0n}(\phi)$ and $D_{1n}(\phi)$ are bounded below by $c_3$, 
          for a fixed a $c_3\in(0,c_2)$, where $c_2$ comes from assumption C4 \\
         $K_{2n}'$ &  $K_{2n}' \supset K_{2n}$ be similarly defined but with the weaker requirement 
           that the minimum eigenvalues of $D_{0n}(\phi)$ and $D_{1n}(\phi)$ are both positive. \\
          $K_n$ & $K_n=K_1\times K_{2n}$ \\
            $K_n'$ & $K_n'=K_1\times K_{2n}'$ \\
         $K_n^{\ast}$&  $K_n^{\ast} = K_n$ if $K_n$ is non-empty, $K_n^{\ast}=K_n'$ if $K_n$ is empty but  $K_n'$ is non-empty,  and $K_n^{\ast} = K$ otherwise. \\
          $\tilde{\theta}_n$ & The $\arg\min$ over $K_n^{\ast}$ of $M_n(\theta)$ \\
        &   $\tilde{\theta}_n=(\tilde{\zeta}_n,\tilde{\phi}_n)$, $\tilde{\zeta}_n=(\tilde{\beta}_n,\tilde{\delta}_n)$ and  $\tilde{\phi}_n=(\tilde{\omega}_n,\tilde{\gamma}_n)$. \\
          $V$  & $V \equiv (V_1, \ldots, V_n), V_i=\ind\{\omega_1'X_i-\gamma_1>0\} -\ind\{\omega_1'X_i-\gamma_1\leq 0\}$, $1\leq i\leq n$ for a fixed point $(\omega_1, \gamma_1)$.  \\
          $\Phi_n$ & A level set. $\Phi_n=\Big\{(\omega,\gamma)\in S^{p-1}\times[l_0,u_0]:$ \\ & \hspace{3cm}$ \ind\{\omega'X_i-\gamma> 0\}-\ind\{\omega'X_i-\gamma\leq 0\}=V_i,\;1\leq i \leq n\Big\}$ \\
           $C_R^n(\omega)$ &  $C_R^n(\omega)=C_U^n(\omega)-C_L^n(\omega)$, where  $C_L^n(\omega)=\max_{1\leq i\leq n:V_i=-1}\omega'X_i$ and $C_U^n(\omega)=\min_{1\leq i\leq n:V_i=1}\omega'X_i$ \\
           $R_n$ & $R_n=\{\omega\in S^{p-1}:C_R^n(\omega)>0\}$\\
           $\hat{\omega}$ & Mean-midpoint of  $\Phi_n$. $\hat{\omega}=\frac{\int_{R_n}\omega C_R^n(\omega)d\nu(\omega)}{\left\|\int_{R_n}\omega C_R^n(\omega)d\nu(\omega)\right\|}$;
            $\nu$ is the uniform measure on $S^{p-1}$\\
            $\hat{\gamma}$ & Estimation of $\gamma$ using the mean-midpoint. $\hat{\gamma}= \left[C_L^n(\hat{\omega})+C_U^n(\hat{\omega})\right]/2$ \\
            $(\check{\omega},\check{\gamma}^{\ast})$  &  Mode-midpoint of  $\Phi_n$.  The argmax of the function $(\omega,\gamma)\mapsto\lambda$ over $\Phi_n$, subject to  $C_L^n(\omega)-\gamma\leq 0$ and  $C_U^n(\omega)-\gamma\geq\lambda$
            \\
             $\check{\gamma}$ & Estimation of $\gamma$ using the mode-midpoint $\check{\gamma}=\left[C_L^n(\check{\omega})+C_U^n(\check{\omega})\right]/2$ \\

             $\hat{\theta}_n$ & Sequence of estimators using the mean-midpoint. $\hat{\theta}_n=(\tilde{\beta}_n,\tilde{\delta}_n,\hat{\omega}_n,\hat{\gamma}_n)$, where $(\hat{\omega}_n,\hat{\gamma}_n)=(\hat{\omega},\hat{\gamma})$ \\
            $\check{\theta}_n$ &  Sequence of estimators using the mode-midpoint. $\check{\theta}_n=(\tilde{\beta}_n,\tilde{\delta}_n,\check{\omega}_n,\check{\gamma}_n)$, where $(\check{\omega}_n,\check{\gamma}_n)=(\check{\omega},\check{\gamma})$ \\
             \hline
Section 4\\
    \hline
        $M(\theta)$ &  $M(\theta) =Pm_{\theta}(X,Y,Z)$ \\
    $\Theta_{\eta}$ & $\Theta_{\eta}\equiv\{\theta\in K:\tilde{d}(\theta,\theta_0)\leq\eta\}$  for any $\eta>0$,  where $\tilde{d}^2(\theta, \theta_0)=\|\zeta-\zeta_0\|^2+\|\omega-\omega_0\|+|\gamma-\gamma_0|$. \\
    $\mathbb{G}_n$ & $\mathbb{G}_n= n^{1/2} (\ep_n -P)$, the empirical process associated with a sample of the random variables of size $n$\\
    ${\cal M}_{\eta}$ & ${\cal M}_{\eta}=\left\{m_{\theta}(Y,Z,X)-m_{\theta_0}(Y,Z,X)\;:\;\theta\in\Theta_{\eta}\right\}$ \\
    \hline
Section 5\\
    \hline
      $H_n$& $H_n=\re^d\times\re^d\times n(S^{p-1}-\omega_0)\times n([l_0,u_0]-\gamma_0)$, the index $h=(h_1,h_2,h_3,h_4)\in H_n$  \\
      $Q_n(h)$& $Q_n(h)=n\left[M_n(\beta_0+h_1/\sqrt{n},\delta_0+h_2/\sqrt{n},\omega_0+h_3/n,\gamma_0+h_4/n)-M_n(\theta_0)\right]$ \\
      $\tilde{h}_n$ & $\tilde{h}_n\in\arg\min_{h\in H_n} Q_n(h)$, and  by construction, \\ 
 & $\tilde{h}_n=\left(\sqrt{n}(\tilde{\beta}_n-\beta_0),\sqrt{n}(\tilde{\delta}_n-\delta_0),n(\tilde{\omega}_n-\omega_0),n(\tilde{\gamma}_n-\gamma_0)\right)$\\ 
 $\Phi_n'$ & The level set $\Phi_n' \equiv n(\Phi_n-\phi_0)$ \\ 
  $W_1$, $W_2$ & Independent mean zero Gaussian random variables with respective covariances $\Sigma_1=\sigma^2\left(E\left[ZZ'\ind\{U\leq 0\}\right]\right)^{-1}$ and 
   $\Sigma_2=\sigma^2\left(E\left[ZZ'\ind\{U> 0\}\right]\right)^{-1}$ \\
  $Q_{02}(g)$ & $Q_{02}(g)=\sum_{j=1}^{\infty}\ind\{-g_1'\tilde{X}_j^{-}+g_2<-\tilde{U}_j^{-}\leq 0\}\tilde{E}_j^{-}+\ind\{0<\tilde{U}_j^{+}\leq-g_1'\tilde{X}_j^{+}+g_2\}\tilde{E}_j^{+}$ \\
   $\tilde{E}_j^{-}$& $\tilde{E}_j^{-}=\left[(\beta_0-\delta_0)'Z_j^{-}\right]^2+2\epsilon_j^{-}(\beta_0-\delta_0)'Z_j^{-}$, $j\geq 1$ \\
    $(\epsilon_j^{-},\;j\geq 1)$ & i.i.d. realizations of the residual $\epsilon$ from the model~(1).\\ 
    $\tilde{X}_j^{-}$ & $\tilde{X}_j^{-}=\overline{\omega}_0'X_j^{-}$, $j\geq 1$ \\
    $(Z_j^{-},X_j^{-})$ & $(Z_j^{-},X_j^{-})\in\re^d\times\re^p$, $j\geq 1$, are i.i.d. draws from the distribution $G$ of Condition C3\\
     $\tilde{U}_j^{-}$ & $\tilde{U}_j^{-}=\sum_{k=1}^j M_k$, $j\geq 1$, where the $(M_k,k\geq 1)$ are i.i.d. exponential random variables with mean $f_0^{-1}$ \\
     $\tilde{E}_j^{+}$ & $\tilde{E}_j^{+}=\left[(\beta_0-\delta_0)'Z_j^{+}\right]^2-2\epsilon_j^{+}(\beta_0-\delta_0)'Z_j^{+}$, $j\geq 1$ \\
     $(\epsilon_j^{+},\;j\geq 1)$ & $(\epsilon_j^{+},\;j\geq 1)$ is an independent random replication of the sequence $(\epsilon_j^{-},\;j\geq 1)$\\
       $\tilde{X}_j^{+}$ & $\tilde{X}_j^{+}=\overline{\omega}_0'X_j^{+}$, $j\geq 1$ \\
      $(Z_j^{+},X_j^{+})$ & An independent random replication of $(Z_j^{-},X_j^{-})$ \\ 
      $\tilde{U}_j^{+}$ & An independent random realization of $\tilde{U}_j^{-}$\\
      $N(t)$ & A homogeneous Poisson process on $[0, \infty)$ with intensity 2$f_0$\\
      $B(t)$ & A white-noise type Bernoulli random variable with success probability $1/2$ \\
      $\epsilon(t)$ & A white-noise type process with distribution the same as the residual in the model (1) \\
      $g$ & $g \in \re^{p-1} \times \re$ \\
        $\tilde{\Phi}_0$ & $\tilde{\Phi}_0=\arg\min_{g\in\re^{p-1}\times\re}Q_{02}(g)$ \\
        $V_j^{-}(g), V_j^{+}(g)$ &   $V_j^{-}(g)=\ind\{g_1'\tilde{X}_j^{-}-g_2-\tilde{U}_j^{-}>0\}-\ind\{g_1'\tilde{X}_j^{-}-g_2-\tilde{U}_j^{-}\leq 0\}$ and \\
        &   $V_j^{+}(g)=\ind\{g_1'\tilde{X}_j^{+}-g_2+\tilde{U}_j^{+}>0\}-\ind\{g_1'\tilde{X}_j^{+}-g_2+\tilde{U}_j^{+}\leq 0\}$ \\
       $\tilde{C}_R(g_1)$ & $\tilde{C}_R(g_1)=\tilde{C}_U(g_1)-\tilde{C}_L(g_1)$, where \\
       & $\tilde{C}_L(g_1)=\left(\max_{1\leq j\leq\tilde{J}_0^{-}:V_j^{-}=-1}g_1'\tilde{X}_j^{-}-\tilde{U}_j^{-}\right)\vee\left(\max_{1\leq j\leq\tilde{J}_0^{+}:V_j^{+}=-1}g_1'\tilde{X}_j^{+}+\tilde{U}_j^{+}\right)$, and \\
          & $\tilde{C}_U(g_1)=\left(\min_{1\leq j\leq\tilde{J}_0^{-}:V_j^{-}=1}g_1'\tilde{X}_j^{-}-\tilde{U}_j^{-}\right)\wedge\left(\min_{1\leq j\leq\tilde{J}_0^{+}:V_j^{+}=1}g_1'\tilde{X}_j^{+}+\tilde{U}_j^{+}\right)$.   \\
     $\hat{g}=(\hat{g}_1.\hat{g}_2)$ &  $\hat{g}_1=\frac{\int_{\tilde{R}_0}g_1\tilde{C}_R(g_1)d\mu(g_1)}{\int_{\tilde{R}_0}\tilde{C}_R(g_1)d\mu(g_1)}$ $\hat{g}_2=\left[\tilde{C}_L(\hat{g}_1)+\tilde{C}_U(\hat{g}_1)\right]/2$\\ 
      $\check{g}=(\check{g}_1,\check{g}_2)$ &  $\check{g}_1=\arg\max_{g_1\in\tilde{R}_0}\tilde{C}_R(g_1)$  $\check{g}_2=\left[\tilde{C}_L(\check{g}_1)+\tilde{C}_U(\check{g}_1) \right]/2$ \\    
 \hline
Sections 5.1--5.2\\
    \hline
       $h_{\ast n}$ & A function $\re^{p-1}\mapsto n(S^{p-1}-\omega_0)$ defined by
$h_{\ast n}(g_1)=\overline{\omega}_0 g_1-n\left\{1-\sqrt{\left(1-\frac{\|g_1\|^2}{n^2}\right)_{+}}\right\}\omega_0$ \\
$r_n$ & For $k \in (0, \infty)$,  $r_n(k)=n\left[1-\sqrt{\left(1-\frac{(0.9 k_1^{-1}k)^2}{n^2}\right)_{+}}\right]$ \\
 $H_{20}^{\ast}$ &  For $k \in (0, \infty)$, $H_{20}^{\ast}(k)\equiv\{(g_1,g_2)\in\re^{p-1}\times\re:\:\|g_1\|\leq 0.9 k_1^{-1}k,\;|g_2|\leq k\}$ \\ 
  $H_{2n}$ &  For $k \in (0, \infty)$, $H_{2n}(k)\equiv\{(h_3,h_4)\in n(S^{p-1}-\omega_0)\times n([l,u]-\gamma_0):\;\|h_3\|^2\leq (0.9 k_1^{-1}k)^2+r_n^2(k),\;|h_4|\leq k\}$ \\
  $H_0^{\ast}$&   For $k \in (0, \infty)$, $H_0^{\ast}(k)=\{(h_1,h_2,g_1,g_2)\in\re^d\times\re^d\times H_{20}^{\ast}(k):\; \|h_1\|\vee\|h_2\|\leq k\}$ \\
   $H_n$ &  For $k \in (0, \infty)$, $H_n(k)=\{(h_1,h_2,h_3,h_4)\in\re^d\times\re^d\times H_{2n}(k):\|h_1\|\vee\|h_2\|\leq k\}$ \\
$\tilde{h}_n$ &  For $k \in (0, \infty)$, $\tilde{h}_n(k)=(\tilde{h}_{1n}(k),\tilde{h}_{2n}(k),\tilde{h}_{3n}(k),\tilde{h}_{4n}(k))\in\arg\min_{h\in H_n(k)} Q_n(h)$  \\
    $\Phi_n'$ &  For $k \in (0, \infty)$, $\Phi_n'(k)$ is the level set containing $(\tilde{h}_{3n}(k),\tilde{h}_{4n}(k))$. \\
$m_{2n}^{-}, m_{2n}^{+}$ &  $m_{2n}^{-}(k)=\sum_{i=1}^n \ind\{-4k \leq  nU_i<-2k\}$, $m_{2n}^{+}(k)=\sum_{i=1}^n\ind\{2k<nU_i\leq 4k\}$ \\
  $F_n$ &  For $k \in (0, \infty)$,  $F_n(k)=\ind\{m_{2n}^{-}(k)\wedge m_{2n}^{+}(k)\geq 1\}$ \\
  $V(k)$ & $V_i(k)\equiv\ind\{\tilde{h}_{3n}'X_i-\tilde{h}_{4n}+nU_i>0\}-\ind\{\tilde{h}_{3n}'X_i-\tilde{h}_{4n}+nU_i\leq 0\}$ \\
  $C_R^{nk}$ & $C_R^{nk}(h_3)=C_U^{nk}(h_3)-C_L^{nk}(h_3)$, where 
    $C_L^{nk}(h_3)=\max_{1\leq i\leq n:\;V_i(k)=-1}h_3'X_i+nU_i$, ~ $C_U^{nk}(h_3)=\min_{1\leq i\leq n:\;V_i(k)=1}h_3'X_i+nU_i$  \\
   $R_n$ & For  $k \in (0, \infty)$, $R_n(k)=\{h_3\in n(S^{p-1}-\omega_0):\; C_R^{nk}(h_3)>0,\; \|h_3\|^2\leq (0.9 k_1^{-1}k)^2+r_n^2(k)\}$ \\
 $\nu_n$ & The uniform measure on $n(S^{p-1}-\omega_0)$\\
  $\hat{h}_n$ &  Mean-midpoint estimator $\hat{h}_n(k)\equiv(\hat{h}_{1n}(k),\hat{h}_{2n}(k),\hat{h}_{3n}(k),\hat{h}_{4n}(k))$, where 
   $\hat{h}_{3n}(k)=\frac{\int_{R_n(k)} h C_R^{nk}(h)d\nu_n(h)}{\int_{R_n(k)} C_R^{nk}(h)d\nu_n(h)}$, and  $\hat{h}_{4n}(k)=\left[C_L^{nk}(\hat{h}_{3n}(k))+C_U^{nk}(\hat{h}_{3n}(k))\right]/2$. \\
   $\tilde{h}_{j0}$& $\tilde{h}_{j0}(k)=\arg\min_{h\in\re^d:\;\|h\|\leq k}h'\Sigma_j h-2 h'\Sigma_j W_j$, $j=1,2$. \\
 $\tilde{g}$ &$\tilde{g}(k)=(\tilde{g}_{1}(k),\tilde{g}_{2}(k))$, where $\tilde{g}(k)\in\arg\min_{g\in H_{20}^{\ast}(k)} Q_{02}(g)$\\
 $\tilde{h}_0$ & $\tilde{h}_0(k)=(\tilde{h}_{10}(k),\tilde{h}_{20}(k),\overline{\omega}_0\tilde{g}_{1}(k),\tilde{g}_{2}(k))$ \\

  $\tilde{\Phi}_0$  &  $\tilde{\Phi}_0$ is the level set in $\re^{p-1}\times\re$ which contains $\tilde{g}_0(k)$, 
   i.e., $\tilde{\Phi}_0(k)$ is the set of all $g\in\re^{p-1}\times\re$ such that $V_{j}^{-}(g)=\tilde{V}_{j}^{-}(k)$ and $V_{j}^{+}(g)=\tilde{V}_{j}^{+}(k)$, for all $j\geq 1$. \\

 $m_{20}^{-}, m_{20}^{+}$ &  $m_{20}^{-}(k)$ =$\sum_{j=1}^{\infty}\ind\{2k<\tilde{U}_j^{-}\leq 4k\}$ and $m_{20}^{+}(k)=\sum_{j=1}^{\infty}\ind\{2k<\tilde{U}_j^{+}\leq 4k\}$\\
 $F_0$ & $F_0(k)=\ind\{m_{20}^{-}(k)\wedge m_{20}^{+}(k)\geq 1\}$\\
  $\tilde{J}_0^{-}$, $\tilde{J}_0^{+}$ & $\tilde{J}_0^{-}(k)=\min\{j\geq 1:\; \tilde{U}_j^{-}\geq 4k\}$ and $\tilde{J}_0^{+}(k)=\min\{j\geq 1:\; \tilde{U}_j^{+}\geq 4k\}$ \\
   $\tilde{C}_R^k$ & $\tilde{C}_R^k(g_1)=\tilde{C}_U^k(g_1)-\tilde{C}_L^k(g_1)$, where \\
    & $\tilde{C}_L^k(g_1)=\left(\max_{1\leq j\leq\tilde{J}_0^{-}(k):\tilde{V}_{j}^{-}(k)=-1}g_1'\tilde{X}_j^{-}-\tilde{U}_j^{-}\right)\vee$ \\
    & \hspace{0.55in} $\left(\max_{1\leq j\leq\tilde{J}_0^{+}(k):\tilde{V}_{j}^{+}(k)=-1}g_1'\tilde{X}_j^{+}+\tilde{U}_j^{+}\right)$ and \\ & $\tilde{C}_U^k(g_1)=\left(\min_{1\leq j\leq\tilde{J}_0^{-}(k):\tilde{V}_{j}^{-}(k)=1}g_1'\tilde{X}_j^{-}-\tilde{U}_j^{-}\right)\wedge$ \\ 
    & \hspace{0.55in} $\left(\min_{1\leq j\leq\tilde{J}_0^{+}(k):\tilde{V}_j^{+}(k)=1}g_1'\tilde{X}_j^{+}+\tilde{U}_j^{+}\right)$ \\
    $\tilde{R}_0$ & $\tilde{R}_0(k)=\{g_1\in\re^{p-1}:\;\tilde{C}_R^k(g_1)>0,\;\|g_1\|\leq 0.9 k_1^{-1}k\}$ \\
    $\hat{h}_0$ & The mean-midpoint limiting estimator. $\hat{h}_0(k)\equiv (\tilde{h}_{10}(k),\tilde{h}_{20}(k), \overline{\omega}_0\hat{g}_1(k),\hat{g}_{2}(k))$, where 
   $\hat{g}_1(k)=\frac{\int_{\tilde{R}_0(k)} g \tilde{C}_R^{k}(g)d\mu(g)}{\int_{\tilde{R}_0(k)} \tilde{C}_R^{k}(g)d\mu(g)}$, and  $\hat{g}_2(k)=\left[\tilde{C}_L^{k}(\hat{g}_{1}(k))+\tilde{C}_U^{k}(\hat{g}_{1}(k))\right]/2$. \\
   $\check{h}_{3n}$ & $\check{h}_{3n}(k)$ is the unique argmax of $C_R^{nk}(h_3)$ over $R_n(k)$ \\
    $\check{h}_4$ &$\check{h}_4(k)=\left[C_L^{nk}(\check{h}_{3n}(k))+C_L^{nk}(\check{h}_{3n}(k))\right]/2$ \\
    $\check{h}_n$ & $\check{h}_n(k)=(\tilde{h}_{1n}(k),\tilde{h}_{2n}(k),\check{h}_{3n}(k),\check{h}_{4n}(k))$ \\
     $\check{h}_0$ &$\check{h}_0(k)=(\tilde{h}_{10}(k),\tilde{h}_{20}(k),\overline{\omega}_0\check{g}_1(k),\check{g}_2(k))$, where  
        $\check{g}_1(k)$ is the unique argmax of $\tilde{C}_R^k(g_1)$ over $g_1\in\tilde{R}_0(k)$ and $\check{g}_2(k)=\left[\tilde{C}_L^k(\check{g}_1(k))+ \tilde{C}_U^k(\check{g}_1(k))\right]/2$.\\
  $F$ &  $F\subset\re^{2d+p+1}$ is a closed set.  \\
 $D_n$ & A diagonal matrix made of a vector $(\sqrt n \boldsymbol 1_{2d}^\top, n\boldsymbol 1_{p+1}^\top)^\top$ \\
  $\rho$ &  $\rho\equiv\sqrt{\Pr(\omega_0'X - \gamma_0)\{1-\Pr(\omega_0'X - \gamma_0)\}}$\\ 
   $\Xi(\beta_0,\delta_0)$ & $\Xi(\beta_0,\delta_0)\equiv \lim_{u\searrow0}E\Big[|\beta_0 - \delta_0|'Z\Big| U=u\Big]$\\
   \hline
Section 5.3\\
    \hline
   $J_n^{-}, J_n^{+}$ & $J_n^{-}(k)=\{1\leq i\leq n:\;-4k\leq nU_i\leq 0\}$, $J_n^{+}(k)=\{1\leq i\leq n:\;0< nU_i\leq 4k\}$ \\
     $m_n^{-}, m_n^{+}$& $m_n^{-}(k)=m_{1n}^{-}(k)+m_{2n}^{-}(k)$ and $m_n^{+}(k)=m_{1n}^{+}(k)+m_{2n}^{+}(k)$,  where \\
     & $m_{1n}^{-}(k)=\sum_{i=1}^n\ind\{-2k\leq nU_i\leq 0\}$, $m_{1n}^{+}(k)=\sum_{i=1}^n\ind\{0<nU_i\leq 2k\}$, \\
    &  $m_{2n}^{-}(k)=\sum_{i=1}^n \ind\{-4k \leq  nU_i<-2k\}$, and $m_{2n}^{+}(k)=\sum_{i=1}^n\ind\{2k<nU_i\leq 4k\}$ \\
   $\hat{\Sigma}_{1n}, \hat{\Sigma}_{2n}$&   $\hat{\Sigma}_{1n}=\{n^{-1}\sum_{i=1}^n\ind\{U_i\leq 0\} Z_iZ_i'\}^{-1}$, $\hat{\Sigma}_{2n}=\{n^{-1}\sum_{i=1}^n \ind\{U_i>0\} Z_iZ_i'\}^{-1}$ \\ 
   $W_{1n}, W_{2n}$ & $W_{1n}=n^{-1/2}\sum_{i=1}^n\ind\{U_i\leq 0\}\epsilon_iZ_i$, and $W_{2n}=n^{-1/2}\sum_{i=1}^n\ind\{U_i>0\}\epsilon_iZ_i$ \\
   ${\cal P}_n$ &
    ${\cal P}_n(k)=\left\{m_n^{-}(k),(X_i,Z_i,\epsilon_i,-nU_i):\;i\in J_n^{-}(k), m_{1n}^{-}(k), m_{2n}^{-}(k), \hat{\Sigma}_{1n}, W_{1n};\;\right.$\\
   & \hspace{0.5in} $\left.m_n^{+}(k),(X_i,Z_i,\epsilon_i,nU_i):\;i\in J_n^{+}(k), m_{1n}^{+}(k), m_{2n}^{+}(k), \hat{\Sigma}_{2n}, W_{2n}\right\}$\\
 $m_0^{-}, m_0^{+}$ & $m_0^{-}(k)=m_{10}^{-}(k)+m_{20}^{-}(k)$, and $m_0^{+}(k)=m_{10}^{+}(k)+m_{20}^{+}(k)$, where \\
 & $m_{10}^{-}(k)=\sum_{j=1}^{\infty}\ind\{0\leq \tilde{U}_j^{-}\leq 2k\}$, $m_{10}^{+}(k)=\sum_{j=1}^{\infty}\ind\{0<\tilde{U}_j^{+}\leq 2k\}$ \\
 & $m_{20}^{-}(k)$ =$\sum_{j=1}^{\infty}\ind\{2k<\tilde{U}_j^{-}\leq 4k\}$ and $m_{20}^{+}(k)=\sum_{j=1}^{\infty}\ind\{2k<\tilde{U}_j^{+}\leq 4k\}$ \\
  ${\cal P}_0$ &
   ${\cal P}_0(k)=\left\{m_0^{-}(k),(X_j^{-},Z_j^{-},\epsilon_j^{-},\tilde{U}_j^{-}):\;1\leq j\leq m_0^{-}(k), m_{10}^{-}(k), m_{20}^{-}(k), \Sigma_1, W_1;\;\right.$\\
 &  \hspace{0.5in} $\left.m_0^{+}(k),(X_j^{+},Z_j^{+},\epsilon_j^{+},\tilde{U}_j^{+}):\;1\leq j\leq m_0^{+}(k), m_{10}^{+}(k), m_{20}^{+}(k), \Sigma_2, W_2\right\}$\\
 $\mathbb{D}_{\infty}^q$ & The space of infinite sequences $x_1,x_2,\ldots$, with $x_j\in\re^q$ for all $j\geq 1$ \\
  $\mathbb{Z}^{0+}$ & The set of non-negative integers \\
 $\mathbb{E}^q$ & $\mathbb{E}^q$ is consisting of the set of elements $(x_0,\{x_j:\;j\geq 1\})\in\mathbb{Z}^{0+}\times\mathbb{D}_{\infty}^q$ such that $x_j=0$ for all $j>x_0$ \\
$d_{\ast}$ &   For $x,y\in\mathbb{E}^q$, define the metric $d_{\ast}(x,y)=|x_0-y_0|+\max_{1\leq j\leq x_0\wedge y_0}\|x_j-y_j\|$, where $d_{\ast}(x,y)=0$ when $x_0=y_0=0$. \\
$W_{1n}^{\ast}, W^{\ast}_{2n}$ &  $W_{1n}^{\ast}=n^{-1/2}\sum_{i=1}^n\ind\{nU_i<-4k\}\epsilon_i Z_i$ and $W^{\ast}_{2n}=n^{-1/2}\sum_{i=1}^n\ind\{nU_i>4k\}\epsilon_i Z_i$ \\
${\cal P}_n^{\ast}$ &${\cal P}_n^{\ast}(k)=\left\{m_n^{-}(k),(X_i,Z_i,\epsilon_i,-nU_i):\;i\in J_n^{-}(k), m_{1n}^{-}(k), m_{2n}^{-}(k), \hat{\Sigma}_{1n}, W_{1n}^{\ast};\;\right.$\\
&  \hspace{0.5in}$\left.m_n^{+}(k),(X_i,Z_i,\epsilon_i,nU_i):\;i\in J_n^{+}(k), m_{1n}^{+}(k), m_{2n}^{+}(k), \hat{\Sigma}_{2n},W_{2n}^{\ast}\right\}$ \\
$N$ & $N=n-m_n$ \\ 
$\re_{\omega_0}$ & The linear span of $\omega_0$ \\
   \hline      
Section 6\\
    \hline
      $\hat{U}_i$ &  $\hat{U}_i=\hat{\omega}_n'X_i-\hat{\gamma}_n$ \\
       $\hat{\epsilon}_i$ &  $\hat{\epsilon}_i=Y_i-\ind\{\hat{U}_i\leq 0\}\hat{\beta}_n'Z_i-\ind\{\hat{U}_i>0\}\hat{\delta}_n'Z_i$ \\
       $\tilde{U}_n$ & $\tilde{U}_n=n^{-1}\sum_{i=1}^n\hat{U}_i$ \\
       $\tilde{\epsilon}_n$ & $\tilde{\epsilon}_n=n^{-1}\sum_{i=1}^n\hat{\epsilon}_i$ \\
       $\hat{\tau}_n^2$ & $\hat{\tau}_n^2=n^{-1}\sum_{i=1}^n(\hat{U}_i-\tilde{U}_n)^2$ \\
        $\hat{\sigma}_n^2$ &  $\hat{\sigma}_n^2=n^{-1}\sum_{i=1}^n(\hat{\epsilon}_i-\tilde{\epsilon}_n)^2$. \\
       $\tilde{\Sigma}_{1n}$ & $\tilde{\Sigma}_{1n}=\hat{\sigma}_n^2\big[n^{-1}\sum_{i=1}^n$ $\ind\{\hat{U}_i\leq 0\}Z_iZ_i'\big]^{-1}$ \\
        $\tilde{\Sigma}_{2n}$ & $\tilde{\Sigma}_{2n}=\hat{\sigma}_n^2\left[n^{-1}\sum_{i=1}^n\ind\{\hat{U}_i> 0\}Z_iZ_i'\right]^{-1}$ \\ 
        $\tilde{M}_{n}$  &  A $p\times (p-1)$ matrix which is an estimator of $\overline{\omega}_0$ using Gram-Schmidt orthogonalization \\
        $\hat{\eta}_{n1}, \hat{\eta}_{n2}$ & $\hat{\eta}_{n1}=2\hat{\tau}_n n^{-1/5}$, $\hat{\eta}_{n2}=2\hat{\sigma}_n n^{-1/5}$ \\
       $\hat{F}_{n1}, \hat{F}_{n2}$ &   $\hat{F}_{n1}(t)=n^{-1}\sum_{i=1}^n\ind\{\hat{U}_i\leq t\}$ and $\hat{F}_{n2}(t)=n^{-1}\sum_{i=1}^n\ind\{\hat{\epsilon}_i-\tilde{\epsilon}_n\leq t\}$ \\
       $\phi$ & The standard normal density \\
        $\tilde r_n$ & A sequence such that $\tilde r_n n^{-1/2}\rightarrow\infty$ and $\tilde r_n/n\rightarrow 0$ \\
        $\hat{t}_n$ & $\hat{t}_n=\sup\left\{t>0:\;\sum_{i=1}^n\ind\{|\hat{U}_i|\leq t\}\leq r_n\right\}$ \\
        $\tilde{G}_n$ & $\tilde{G}_n=\{(X_i,Z_i):\; |\hat{U}_i|\leq\hat{t}_n\}$ \\
         $\tilde{N}$ & A homogeneous Poisson process on $[0,\infty)$ with intensity $2\hat{f}_{n0}$, where $\hat{f}_{n0}=\int_{\re}\frac{1}{\hat{\eta}_{n1}}\phi\left(\frac{t}{\hat{\eta}_{n1}}\right)d\hat{F}_{n1}(t)$ \\
         $\tilde{E}^{-},  \tilde{E}^{+}_{\ast}$ &    $ \tilde{E}^{-}_{\ast}(t)=\left[(\hat{\beta}_n-\hat{\delta}_n)'\tilde{Z}_{\ast}(t)\right]^2+2\tilde{\epsilon}_{\ast}(t)(\hat{\beta}_n-\hat{\delta}_n)'\tilde{Z}_{\ast}(t)$ and  \\
         & $\tilde{E}^{+}_{\ast}(t)=\left[(\hat{\beta}_n-\hat{\delta}_n)'\tilde{Z}_{\ast}(t)\right]^2-2\tilde{\epsilon}_{\ast}(t)(\hat{\beta}_n-\hat{\delta}_n)'\tilde{Z}_{\ast}(t)$\\
     $(\tilde{X}_{\ast}(t),\tilde{Z}_{\ast}(t))$ &  $(\tilde{X}_{\ast}(t),\tilde{Z}_{\ast}(t))=(\tilde{M}_{n}'X,Z)$ for $(X,Z)$ drawn independently and with replacement from $\tilde{G}_n$ \\
     $\tilde{\epsilon}_{\ast}(t)$ & Random samples taken independently from the density $\hat{\xi}_n$, where $\hat{\xi}_{n}(u)=\int_{\re}\frac{1}{\hat{\eta}_{n2}}\phi\left(\frac{u-t}{\hat{\eta}_{n2}}\right)d\hat{F}_{n2}(t)$ \\
     $\tilde{g}_{\ast}$ &  $\tilde{g}_{\ast}\in\arg\max_{g\in\re^{p-1}\times\re}\tilde{Q}_{02}^{\ast}(g)$ \\ 
   \hline   
\end{longtable}

\subsection{Some emprical process notations for Section \ref{sec:rateconv}}
\label{asec:empirical_notation}
We give a few definitions---such as outer expectation $E^*[\cdot]$, normalized empirical process operator $\mathbb G_n$, and the uniform norm---needed for one of the conditions in Corollary 14.5 used in Section \ref{sec:rateconv}. For an arbitrary, possibly non-measurable, map $X:\Omega\mapsto\bar{\mathbb R}$, where $\Omega$ is an event space and $\bar\cdot$ denotes the closure of a space, 
$\mathbb E^*[X]\equiv \inf_{U\ge \Upsilon}\mathbb E[U]$ is the outer expectation, 
where $\Upsilon$ is a collection of measurable maps, 
$\mathbb E\|X\|_A$ denotes the uniform norm of a random map $X$ over a set $A$, 
and $\mathbb G_n(\cdot) = \sqrt{n}(\mathbb P_n - P)(\cdot)$.

\FloatBarrier

\section{Proofs of Section \ref{sec:est}} 
\label{asec:est}
We show that the mean-midpoint exists, is well-defined, and is contained in $\Phi_n$.
\begin{lemma}\label{meanmidexistence}
Assume $n$ is large enough so that $K_{2n}'$ is non-empty and $(\omega_1,\gamma_1)\in K_{2n}'\cap\Phi_n$. Then the following are true:
\begin{enumerate}
\item The set $R_n$
\begin{enumerate}
\item is open and contains an open ball on $S^{p-1}$, 
\item is the intersection of a finite number of open half-spheres $R_{nj}$, $j=1,\ldots,m$, and thus if $\omega\in R_n$ then $-\omega\not\in R_n$ (a property we  call ``strict positivity''), and
\item is geodesically connected (specifically, the geodesic segment on $S^{p-1}$ joining any two points in $R_n$ is contained in $R_n$);
\end{enumerate}
\item $\Phi_n=\left\{(\omega,\gamma)\in K_2: \omega\in R_n,\;\mbox{and}\;C_L^n(\omega)\leq\gamma<C_U^n(\omega)\right\}$; and
\item $(\hat{\omega},\hat{\gamma})\in\Phi_n$ with $\hat{\gamma}\in[l_0,u_0]$.
\end{enumerate}
\end{lemma}

\begin{proof}[Proof of Lemma \ref{meanmidexistence}]
Let $J_R(n)$ be the set of all pairs of indices $1\leq j,k\leq n$ such that for each 
$(j,k)\in J_R(n)$, $V_j=-1$ and $V_k=1$. It can be seen that $J_R(n)$ is finite and also non-empty, when $n$ is
large enough so that $K_{2n}'$ is non-empty and $(\omega_1,\gamma_1)\in K_{2n}'$. By definition of $\Phi_n$, 
we know that $X_j\neq X_k$. Let $\omega_1=(X_k-X_j)/ \|X_k-X_j\|$ and use the plane in $\re^p$ orthogonal to 
$\omega_1$ to split $S^{p-1}$ into two equal parts, which exclude the intersection of the plane and the sphere.
It is easy to see that each part is an open half-sphere, is geodesically connected, and strictly positive. 
Take the half-sphere which includes $\omega_1$, and label it $R_{n1}$. Now do this for all pairs in $J_R(n)$, 
and, when there are duplicate open half-spheres, choose only one of them. This leads to the set of all unique 
open half-spheres, $R_{nl},l=1,\ldots,m$, which are generated in this manner by pairs in $J_R(n)$. Let 
$\tilde{R}_n=\cap_{l=1}^m R_{nl}$. $\tilde{R}_n$ is by construction the intersection of a finite number of 
open half-spheres. Moreover, since properties of openness, geodesic connectedness, strict positivity are 
preserved when taking intersections, $\tilde{R}_n$ also has all of these properties. 

The first part of the lemmas is complete if we show that $\tilde{R}_n=R_n$. Let $\omega\in \tilde{R}_n$. We have by construction that for every $(j,k)\in J_R(n)$, $\omega'(X_k-X_j)>0$, and thus $C_R^n(\omega)>0$, which implies that $\omega\in R_n$, and thus $\tilde{R}_n\subset R_n$. Let $\omega\in R_n$. Then $C_R^n(\omega)>0$, which implies that 
$\max_{1\leq i\leq n:V_i=-1} \omega'X_i<\min_{1\leq i\leq n:V_i=1} \omega'X_i$. Thus $\omega'(X_k-X_j)>0$ for all $(j,k)\in J_R(n)$, and hence $\omega\in R_n$. Thus $\tilde{R}_n=R_n$. The fact that $R_n$ contains an open ball follows automatically because non-empty open sets must contain open balls. Thus Part 1 of the Lemma is proved. Part 2 follows from the definitions involved. 

For Part 3, consider the function $f_{R}^n(\omega)=C_R^n(\omega)/\int_{\omega\in R_n} C_R^n(\omega)d\nu(\omega)$, and note that this is a well-defined bounded density with respect to $\nu(\omega)$, and that it has support only on $R_n$. As the boundary of $R_n$ is included in the complement of $R_n$ and because $C_R^n(\omega)$ is strictly positive for all $\omega\in R_n$, the projection of $\int_{R_n}f_R^n(\omega)d\nu(\omega)$ onto $S^{p-1}$, which is precisely $\hat{\omega}$, must be $\in R_n$. The fact that $\hat{\gamma}\in[l_0,u_0]$ follows from the fact that $-k_1\leq C_L^n(\hat{\omega})\leq C_U^n(\hat{\omega})\leq k_1$ by Condition~C1.
\end{proof}

The following lemma states that $(\check{\omega},\check{\gamma})$ is unique and contained in $\Phi_n$.
\begin{lemma}\label{modemidexistence}
Assume $n$ is large enough so that $K_{2n}'$ is non-empty and $(\omega_1,\gamma_1)\in K_{2n}'\cap\Phi_n$. Then the following are true:
\begin{enumerate}
\item $(\check{\omega},\check{\gamma})$, where 
$\check{\gamma}=\left[C_L^n(\check{\omega})+C_U^n(\check{\omega})\right]/2,$
is the unique maximizer over $\Phi_n$ of (\ref{eq1.modemin}), and
\item 
$(\check{\omega},\check{\gamma})\in\Phi_n$.
\end{enumerate}
\end{lemma}

\begin{proof}[Proof of Lemma \ref{modemidexistence}]
    The first part of the lemma follows from Lemma~\ref{modemidexistence.unique} below by substituting $m=n$, $q=p$, $R=S^{p-1}$, and both $x_j=X_j$ and $v_j=\ind\{\omega_1'X_j-\gamma_1> 0\}-\ind\{\omega_1'X_j-\gamma\leq 0\}$, for $1\leq j\leq n$. The fact that $\check{\gamma}\in[l_0,u_0]$ follows from the fact that $-k_1\leq C_L^n(\check{\omega})\leq k_1$. The second part of the lemma follows from the straightforward observation that 
    $C_L^n(\check{\omega})-\check{\gamma}\leq 0$ and $C_U^n(\check{\omega})-\check{\gamma}>0$, and thus the parameter estimator is in $\Phi_n$.
\end{proof}

\begin{lemma}\label{modemidexistence.unique}
Let $(x_j,v_j)\in\re^q\times\{-1,1\}$, $1\leq j\leq m<\infty$, and let $R\subset S^{q-1}$ be geodesically connected and closed. Assume there exists $(k,l)\in\{1,\ldots,m\}$ such that $v_k=-1$ and $v_l=1$, and that there exists a $(\omega_1,u_1)\in R\times\re$ such that
\[v_j=\ind\{\omega_1'x_j-u_1> 0\}-\ind\{\omega_1'x_j-u_1\leq 0\},\]
for all $1\leq j\leq m$. Then there exists a unique $\arg\max$ over $R\times\re$, $(\check{\omega},\check{u})$, of
\begin{eqnarray}
&&(\omega,u)\mapsto\lambda,\;\mbox{subject to}\nonumber\\
&&\max_{1\leq j\leq m:v_j=-1} \omega'x_j-u\leq -\lambda\; \mbox{and}\;\min_{1\leq j\leq m:v_j=1}\omega'x_j-u\geq\lambda.\label{eq2.modemid}
\end{eqnarray}
Moreover, $\check{u}=\{\max_{1\leq j\leq m: v_j=-1}\check{\omega}'x_j + \min_{1\leq j\leq m: v_j=1}\check{\omega}'x_j\}/2$, and the inequalities in the constraints become equalities when evaluated at this $\arg\max$. 
\end{lemma}

\begin{proof} [Proof of Lemma \ref{modemidexistence.unique}]
Define the cone $C=\{x\in\re^q:x=rp,\;\mbox{where}\;p\in R,\; r\in[0,\infty)\}$, and note that $R\subset C$. We first verify that $C$ is convex. Let $u,v\in C$. If $u=v=0$, then it is obvious that $\alpha u+(1-\alpha)v=0\in C$ for all $0<\alpha<1$. If one is zero but one is not, then it is also easy to verify that $\alpha u+(1-\alpha)v\in C$ for all $0<\alpha<1$. Now suppose both are non-zero. Let $u_1=u/\|u\|$ and $v_1=v/\|v\|$. If $u_1=v_1$, then all points of the form $ru_1$ are in $C$, and thus $\alpha u+(1-\alpha)v\in C$ for all $0<\alpha<1$. If $u_1=-v_1$, then the inclusion also holds. Now assume that $u_1\neq v_1$ and $u_1\neq -v_1$. It is easy to verify that the segment in $\re^q$ joining $u_1$ and $v_1$ excludes the zero point and has a projection on $S^{q-1}$ which is precisely the shortest geodesic line joining $u_1$ and $v_1$, which geodesic we denote $l(u,v)$. It is now clear that, for any $0<\alpha<1$, $\alpha u+(1-\alpha)v$ projects onto $S^{q-1}$ at some point on $l(u,v)$. Hence, $C$ is convex. 

Now note that $(\omega_1,u_1)\in C\times\re$. Define the functions $C_L^{\ast}(\omega)=\max_{1\leq j\leq m:v_j=-1}\omega'X_j$, $C_U^{\ast}=\min_{1\leq j\leq m:v_j=1}\omega'X_j$, and $C_R^{\ast}(\omega)=C_U^{\ast}(\omega)-C_L^{\ast}(\omega)$. Note that the minimum of a finite set of linear functions is convex and the maximum is concave, hence $C_R^{\ast}(\omega)$ is also convex. Note that, by assumption, $C_L^{\ast}(\omega_1)-u_1\leq 0$ and $C_U^{\ast}(\omega_1)-u_1>0$. Hence, $C_R^{\ast}(\omega_1)>0$. Let $\eta=C_R^{\ast}(\omega_1)$, define $\omega_2=\omega_1/\eta$, and note that now $C_R^{\ast}(\omega_2)\geq 1$. Consider minimizing over $C$ the function $\omega\mapsto \|\omega\|^2$ subject to the constraint that $C_R^{\ast}(\omega)\geq 1$. Because $\omega\mapsto\|\omega\|^2$ is strictly convex, the constraint is convex, and we have a feasible point $\omega_2$: then there must exist a unique solution $\check{\omega}_1\in C$. Note that $C_R^{\ast}(0)=0$, and this contradicts the constraint. Thus, $\check{\omega}_1\neq 0$. Now suppose $C_R^{\ast}(\check{\omega}_1)>1$. Then there exists  $0<\alpha<1$ such that $C_R^{\ast}(\alpha\check{\omega}_1)=1$. Because $\alpha<1$, the previous sentence implies that $\alpha\check{\omega}_1$ is both feasible and smaller in norm than $\check{\omega}_1$, but this would contradict the uniqueness of the minimum. Thus, $C_R^{\ast}(\check{\omega}_1)=1$. It follows that the minimization we employed is equivalent to minimizing $\|\omega\|^2$ subject to $C_R^{\ast}(\omega/\|\omega\|)\geq 1/\|\omega\|$, where the inequality becomes an equality for the minimizer. This is now equivalent to maximizing $\lambda$ over $\omega\in R$ subject to $C_R^{\ast}(\omega)\geq$ $2\lambda$, where the equality is again achieved for the maximizer. What we have now shown is that there exists a unique $\arg\max$ over $R$, $\check{\omega}_2=\check{\omega}_1/\|\check{\omega}_1\|$, of $\omega\mapsto\lambda$ subject to $C_R(\omega)\geq 2\lambda$, where the equality is achieved.

Now we will verify that $(\check{\omega},\check{u})=\left[ \check{\omega}_2,\{C_L^*(\check{\omega}_2) + C_U^*(\check{\omega}_2)\}/2\right]$ is the unique maximizer of (\ref{eq2.modemid}) with equalities in the constraints. Let $\check{\lambda}$ be the corresponding value of $\lambda$ attained by this maximizer, or $\check{\lambda} = C_R^*(\check \omega)/2$. Suppose there exists some $(\check{\omega}_3,\check{u}_3)\in R\times\re$ which is not equal to $(\check{\omega},\check{u})$ but which satisfies~(\ref{eq2.modemid}) for some  $\check{\lambda}_3\geq\check{\lambda}$. We know that $C_R^{\ast}(\check{\omega})\geq 2 \check{\lambda}$, and thus $\check{\omega}_3=\check{\omega}$ by the previously established uniqueness of the maximizer of $C_R^{\ast}$ over $R$. Hence, $\check{\lambda}_3=\check{\lambda}$. The form of~(\ref{eq2.modemid}) now forces $\check{u}_3=\check{u}$, and thus the desired uniqueness of the maximizers, the form of $\check{u}$, and the equalities in the constraints are established.
\end{proof}
\FloatBarrier


\section{Proofs of Section \ref{sec:rateconv}} 
\label{asec:rateconv}
\subsection{Proof of Theorem \ref{thm.1}}
Before giving the proof of Theorem~\ref{thm.1}, define
\begin{eqnarray*}
A_{00}(\phi)&=&\ind\{\omega_0'X-\gamma_0\leq 0,\omega'X-\gamma\leq 0\},\\
A_{01}(\phi)&=&\ind\{\omega_0'X-\gamma_0\leq 0,\omega'X-\gamma>0\},\\
A_{10}(\phi)&=&\ind\{\omega_0'X-\gamma_0>0,\omega'X-\gamma\leq 0\},\;\mbox{and}\\
A_{11}(\phi)&=&\ind\{\omega_0'X-\gamma_0>0,\omega'X-\gamma>0\},
\end{eqnarray*}
and note that $A_0(\phi)=A_{00}(\phi)+A_{10}(\phi)$ and $A_1(\phi)=A_{01}(\phi)+A_{11}(\phi)$, where $A_j(\phi)$, $j=0,1$, was defined in Section~\ref{sec:est}. 

\begin{proof}[Proof of Theorem~\ref{thm.1}]
Let $\tilde{\theta}_n=(\tilde{\beta}_n,\tilde{\delta}_n,\tilde{\omega}_n,\tilde{\gamma}_n)$ be any sequence in $K_n^{\ast}$ satisfying the conditions of Theorem~\ref{thm.1}. Note that $K_n^{\ast}$ can be compactified by applying the metric
\[d_1(x,y)=\left\|\frac{x}{1+\|x\|}-\frac{y}{1+\|y\|}\right\|,\]
on $\re^d$, to both copies of $\re^d$ in $K$, applying the usual Euclidean metric to $S^{p-1}$ and $[l_0,u_0]$, and then utilizing the resulting product metric. Note that $d_1$ is equivalent to the Euclidean norm on restrictions to compact subsets of $\re^d$.

Define
\begin{eqnarray*}
H_{00}^n&=&A_{00}(\tilde{\phi}_n)\left\{[\epsilon-(\tilde{\beta}_n-\beta_0)'Z]^2\right\}/(1+\|\tilde{\beta}_n\|^2),\\
H_{01}^n&=&A_{01}(\tilde{\phi}_n)\left\{[\epsilon-(\tilde{\delta}_n-\beta_0)'Z]^2\right\}/(1+\|\tilde{\delta}_n\|^2),\\
H_{10}^n&=&A_{10}(\tilde{\phi}_n)\left\{[\epsilon-(\tilde{\beta}_n-\delta_0)'Z]^2\right\}/(1+\|\tilde{\beta}_n\|^2),\;\mbox{and}\\
H_{11}^n&=&A_{11}(\tilde{\phi}_n)\left\{[\epsilon-(\tilde{\delta}_n-\delta_0)'Z]^2\right\}/(1+\|\tilde{\delta}_n\|^2),\\
\end{eqnarray*}
and let $H_0^n=H_{00}^n+H_{10}^n$ and $H_1^n=H_{01}^n+H_{11}^n$.
Note that standard empirical process methods show that the classes $\{A_{ij}(\phi):(\omega,\gamma)\in S^{p-1}\times[l,u]\}$, $0\leq i,j\leq 1$, are Donsker and bounded, and hence Glivenko-Cantelli. Standard arguments also reveal that
$\{[\epsilon-(\beta-\beta_0)'Z]^2/(1+\|\beta\|^2):\beta\in\re^d\}$ is Glivenko-Cantelli and bounded, and thus
$H_{00}^n$ is contained in a Glivenko-Cantelli class. Similar arguments verify that $H_{ij}^n$, for all $0\leq i,j\leq 1$,
are contained in Glivenko-Cantelli classes. Hence $|H_{ij}^n-h_{ij}^n|\rightarrow 0$ and $|H_j^n-h_j^n|\rightarrow 0$, almost surely, where
$h_{ij}^n=EH_{ij}^n$ and $h_j^n=EH_j^n$, for all $0\leq i,j\leq 1$. Note that, as discussed at the beginning of Section \ref{sec:est}, we have $K_n^{\ast}=K_n$ for all $n$ large enough almost surely. Let $\Omega$ be the subset of the probability space where these convergence results hold, and note that the inner probability measure of $\Omega$ is 1 by the definition of (outer) almost sure convergence. Fix an event $\kappa\in\Omega$.

Let $a_{ij}^n=EA_{ij}(\tilde{\phi}_n)$ and $a_j^n=EA_j(\tilde{\phi}_n)$, $0\leq i,j\leq 1$, and let $\tilde{n}$ be a subsequence wherein
$a_{ij}^{\tilde{n}}\rightarrow\tilde{a}_{ij}$, $d_1(\tilde{\beta}_{\tilde{n}},\tilde{\beta})\rightarrow 0$, $d_1(\tilde{\delta}_{\tilde{n}},\tilde{\delta})\rightarrow 0$, $\|\tilde\omega_{\tilde{n}}-\tilde\omega\|\rightarrow 0$, and $|\tilde\gamma_{\tilde{n}}-\tilde\gamma|\rightarrow 0$, for some $0\leq \tilde{a}_{ij}\leq 1$, $\tilde{\beta},\tilde{\delta}\in\overline{\re^d}$, $(\tilde{\omega},\tilde{\gamma})\in S^{p-1}\times[l,u]$, and $0\leq i,j\leq 1$, where $\overline{\re^d}$ denotes inclusion of infinite values. Let $\tilde a_j=\tilde a_{0j}+\tilde a_{1j}$, $j=0,1$, and $\nu_0(X)=\beta_0\ind\{\omega_0'X-\gamma_0\leq 0\}+\delta_0\ind\{\omega_0'X-\gamma_0>0\}$. Note by construction that $\tilde a_0+\tilde a_1=1$. We will now show by contradiction that both $\tilde a_0>0$ and $\tilde a_1>0$. First assume that $\tilde a_1=0$. Accordingly, because $\kappa\in\Omega$ and we are tracing along the subsequence $\tilde{n}$,
\begin{eqnarray*}
\lefteqn{\liminf_{n\rightarrow\infty}M_{\tilde{n}}(\theta_0)\;\geq\;\liminf_{n\rightarrow\infty}M_{\tilde{n}}(\tilde{\theta}_{\tilde{n}})}&&\\ 
&\geq&\liminf_{n\rightarrow\infty}\left[E\left\{\left[c_3\|\tilde{\beta}_{\tilde{n}}-\nu_0(X)\|^2+\epsilon^2\right]\left(A_{00}(\tilde{\phi}_{\tilde{n}})+A_{10}(\tilde{\phi}_{\tilde{n}})\right)\right\}+o(1+\|\tilde{\beta}_{\tilde{n}}\|^2)\right]\\
&\geq&\liminf_{n\rightarrow\infty}\left\{c_3\left[\|\tilde{\beta}_{\tilde{n}}-\beta_0\|^2 \tilde a_{00}+\|\tilde{\beta}_{\tilde{n}}-\delta_0\|^2 \tilde a_{10}\right]+\sigma_0^2+o(1+\|\tilde{\beta}_{\tilde{n}}\|^2)\right\},
\end{eqnarray*}
which implies that
\[\liminf_{n\rightarrow\infty}\left\{c_3\left[\|\tilde{\beta}_{\tilde{n}}-\beta_0\|^2 \tilde a_{00}+\|\tilde{\beta}_{\tilde{n}}-\delta_0\|^2 \tilde a_{10}\right]+o(1+\|\tilde{\beta}_{\tilde{n}}\|^2)\right\}\;\leq\; 0,\]
because $M_{\tilde{n}}(\theta_0)\rightarrow\sigma_0^2$. But this leads to a contradiction as
C2 implies that both $\tilde a_{00}$ and $\tilde a_{10}$ are $>0$ and C5 ensures that $\beta_0\neq\delta_0$. Thus $\tilde a_1>0$. We can also symmetrically argue that $\tilde a_0>0$. Hence we can conclude that $\tilde a_0>0$, $\tilde a_1>0$; and, by reapplying C2, that both $\tilde a_{00}+\tilde a_{01}>0$ and $\tilde a_{10}+\tilde a_{11}>0$.

By recycling previous arguments, we now have, generally, that
\begin{eqnarray*}
\lefteqn{\limsup_{n\rightarrow\infty}M_{\tilde{n}}(\theta_0)\;\geq\;\limsup_{n\rightarrow\infty}M_{\tilde{n}}(\tilde{\theta}_{\tilde{n}})}&&\\ 
&\geq&\limsup_{n\rightarrow\infty}\left\{c_3\left[\|\tilde{\beta}_{\tilde{n}}-\beta_0\|^2 \tilde a_{00}+\|\tilde{\beta}_{\tilde{n}}-\delta_0\|^2 \tilde a_{10}\right.\right.\\
&&+\left.\left. \|\tilde{\delta}_{\tilde{n}}-\beta_0\|^2 \tilde a_{01}+\|\tilde{\delta}_{\tilde{n}}-\delta_0\|^2 \tilde a_{11}
\right]+\sigma_0^2+o(1+\|\tilde{\beta}_{\tilde{n}}\|^2+\|\tilde{\delta}_{\tilde{n}}\|^2)\right\},
\end{eqnarray*}
which implies that
\begin{eqnarray}
\limsup_{n\rightarrow\infty}\left\{\tilde a_{00}\|\tilde{\beta}_{\tilde{n}}-\beta_0\|^2+\tilde a_{10}\|\tilde{\beta}_{\tilde{n}}-\delta_0\|^2
+\tilde a_{01}\|\tilde{\delta}_{\tilde{n}}-\beta_0\|^2\right.&&\nonumber\\
\left.+\tilde a_{11}\|\tilde{\delta}_{\tilde{n}}-\delta_0\|^2+o(1+\|\tilde{\beta}_{\tilde{n}}\|^2+\|\tilde{\delta}_{\tilde{n}}\|^2)\right\}&\leq& 0.\label{eqs4e2}
\end{eqnarray}
The first thing to note is, by the positivity conclusions of the previous paragraph, both $\|\tilde{\beta}_{\tilde{n}}\|=O(1)$ and $\|\tilde{\delta}_{\tilde{n}}\|=O(1)$. 

Now assume that $\tilde a_{00}>0$ and $\tilde a_{11}>0$, and note that this already satisfies the positivity conclusions. Next, note that this assumption forces both $\tilde a_{01}=0$ and $\tilde a_{10}=0$ since, otherwise, we have a contradiction by assumption C5, since a sequence cannot converge simultaneously to two distinct values. Hence $\tilde{\beta}=\beta_0$ and $\tilde{\delta}=\delta_0$. Moreover, we also now have, from the fact that $\tilde a_{01}=0$, that
\begin{eqnarray*}
0&\geq & \limsup_{n\rightarrow\infty}E\left[A_{01}(\tilde{\phi}_{\tilde{n}})\right]\\
&\geq&\liminf_{n\rightarrow\infty}P\left[\omega_0'X-\gamma_0< 0,\tilde{\omega}_{\tilde{n}}'X-\tilde{\gamma}_{\tilde{n}}>0\right]\\
&\geq&P\left[\omega_0'X-\gamma_0<0,\tilde{\omega}'X-\tilde{\gamma}>0\right],
\end{eqnarray*}
where, for the last inequality, we used the Portmanteau theorem applied to the (trivial) weak convergence of $\tilde{\omega}_{\tilde{n}}'X-\tilde{\gamma}_{\tilde{n}}$ to $\tilde{\omega}'X-\tilde{\gamma}$. We can use similar arguments to show that $\tilde a_{10}=0$ implies that $P[\omega_0'X-\gamma_0>0,\tilde{\omega}'X-\tilde{\gamma}<0]=0$. Thus, by assumption C3, $\tilde{\omega}=\omega_0$ and $\tilde{\gamma}=\gamma_0$.

We next assume that $\tilde a_{01}>0$ and $\tilde a_{10}>0$. Note that this also satisfies the positivity conclusions. Arguing along the same lines as in the previous paragraph, we have that $\tilde a_{00}=0$ and $\tilde a_{11}=0$ and, moreover, that $\tilde{\beta}=\delta_0$ and $\tilde{\delta}=\beta_0$. Now redefine $(\beta_0,\delta_0)$ as $(\delta_0,\beta_0)$ and redefine $(\omega_0,\gamma_0)$ as $(-\omega_0,-\gamma_0)$, and note, with this redefinition and after reapplying the assumptions accordingly, we once again conclude that $\tilde{\theta}=\theta_0$. 

\begin{remark}\label{remark4.1}
Our model is unable to discriminate between $(\beta_0,\delta_0,\omega_0,\gamma_0)$ and $(\delta_0,\beta_0,$ $-\omega_0,-\gamma_0)$. Fortunately, this does not matter, provided the assumptions are applied to at least one of the two versions, since the two models generate identical data. Since these two possibilities are indistinguishable, we will simply refer to both versions as $\theta_0=(\beta_0,\delta_0,\omega_0,\gamma_0)$.
\end{remark}

By careful inspection, we can now deduce that the only way to satisfy the positivity conclusions and not obtain a contradiction is for one of the following two possibilities to hold: (i) $\tilde a_{00}>0$, $\tilde a_{11}>0$, and $\tilde a_{01}=\tilde a_{10}=0$; or (ii), $\tilde a_{01}>0$, $\tilde a_{10}>0$, and $\tilde a_{00}=\tilde a_{11}=0$. Allowing for the two possible versions of $\theta_0$, we have now established that $\tilde{\theta}_{\tilde{n}}\rightarrow\theta_0$. Because the convergent subsequence was arbitrary, we have that $\tilde{\theta}_n\rightarrow\theta_0$. As $\kappa\in\Omega$ was arbitrary, and because the inner probability of $\Omega$ is 1, we now have that $\tilde{\theta}_n\rightarrow\theta_0$, outer almost surely.
\end{proof}

\subsection{Proof of Propositions \ref{prop.rate1} and \ref{prop.rate2} }
\label{asec:con}

We now prove the two propositions in reverse order: 

\begin{proof}[Proof of Proposition~\ref{prop.rate2}]
    Note first that
\begin{eqnarray*}
m_{\theta}(Y,Z,X)-m_{\theta_0}(Y,Z,X)&=&A_{00}(\phi)\left[-2\epsilon(\beta-\beta_0)'Z+(\beta-\beta_0)'ZZ'(\beta-\beta_0)\right]\\
&&+A_{01}(\phi)\left[-2\epsilon(\delta-\beta_0)'Z+(\delta-\beta_0)'ZZ'(\delta-\beta_0)\right]\\
&&+A_{10}(\phi)\left[-2\epsilon(\beta-\delta_0)'Z+(\beta-\delta_0)'ZZ'(\beta-\delta_0)\right]\\
&&+A_{11}(\phi)\left[-2\epsilon(\delta-\delta_0)'Z+(\delta-\delta_0)'ZZ'(\delta-\delta_0)\right]\\
&\equiv&B_1(\theta)+B_2(\theta)+B_3(\theta)+B_4(\theta).
\end{eqnarray*}
Standard arguments can be used to verify that each ${\cal B}_j\equiv\{B_j(\theta):\tilde{d}(\theta,\theta_0)\leq\eta\}$ is a VC class, $1\leq j\leq 4$. It is also not hard to verify that both ${\cal B}_1$ and ${\cal B}_4$ have envelop
$F_1=\eta(b_1|\epsilon|+b_2)$ for some finite constants $0<b_1,b_2<\infty$. It is also not hard to argue
that $A_{01}(\phi)\leq A_{01}^{\ast}\equiv\ind\{-b_3(\|\omega-\omega_0\|+|\gamma-\gamma_0|)<\omega_0'X-\gamma_0\leq 0\}$ for some $0<b_3<\infty$ and all $\theta$ such that $\tilde{d}(\theta,\theta_0)$ is small enough. Hence there exist constants $0<b_4,b_5<\infty$
such that ${\cal B}_2$ has envelop $F_2=A_{01}^{\ast}(b_4|\epsilon|+b_5)$, where we note that
$E F_2^2\leq EA_{01}^{\ast}(2b_4^2\sigma_0^2+2b_5^2)$ and $EA_{01}^{\ast}\leq f_0b_3(\|\omega-\omega_0\|+|\gamma-\gamma_0|)$. 
We can similarly argue that ${\cal B}_3$ has envelope $F_3=A_{10}^{\ast}(b_4|\epsilon|+b_5)$ for $A_{10}^{\ast}\equiv\ind\{0<\omega_0'X-\gamma_0\leq b_3(\|\omega-\omega_0\|+|\gamma-\gamma_0|)\}$, after increasing slightly $b_3$, $b_4$ and $b_5$ if needed. Hence $F=2F_1+F_2+F_3$ is an envelop for ${\cal M}_{\eta}$ with $\sqrt{EF^2}\leq b_6\eta $, for some $0<b_6<\infty$. Now Theorem 11.1 of \cite{kosorok2008introduction} yields the desired result with $c_5=b_6$ since ${\cal M}_{\eta}$ is a VC class.
\end{proof}

\begin{proof}[Proof of Proposition~\ref{prop.rate1}]
    
This proof requires the following three lemmas, which we prove at the end of this section:
\begin{lemma}\label{lthm:cr0}
Under assumptions C2 and C3, the conditional joint distribution of $(\nu^{-1}U,(X,Z))$ given $|U|\leq\nu$ converges weakly, as $\nu\downarrow 0$, to $(U_0,W_0)$, where $U_0\perp W_0$, $U_0$ is uniform on $[-1,1]$, and $W_0\sim G$.
\end{lemma}

\begin{lemma}\label{lthm:cr1}
Let $(W,B) \in \mathbb{R}^q\times\re$ be a jointly distributed pair of random variables with $P(B>0)=1$ and $E(W)=\mu$ finite,
and suppose there exists a constant $\eta >0$ such that
$E|u'(W-\mu)| \geq \eta$ for all $u \in S^{q-1}$.
Then, there exists a constant $c > 0$ such that
\begin{equation}
E\left[|a'W-b|(B\wedge 1)\right] \geq c(\|a\|+|b|)
\end{equation}
for all $a \in \mathbb{R}^q$ and $b \in \mathbb{R}$.
\end{lemma}

\begin{lemma}\label{lthm:cr2}
For any $\omega_1,\omega_2\in S^{p-1}$, $|(\omega_2-\omega_1)'\omega_1|=\|\omega_2-\omega_1\|^2/2$.
\end{lemma}

We will utilize these lemmas later in the proof. We first note that
that Conditions {C1, C2, C3, C4 and C5} imply that for all $0<\eta$ small enough, we have for any $\theta\in\Theta_{\eta}$ that
\begin{eqnarray}
M(\theta)-M(\theta_0)&\geq& (c_2/2) a_{00}(\phi)\|\beta-\beta_0\|^2+(c_2/2) a_{11}(\phi)\|\delta-\delta_0\|^2\nonumber\\
&&+(c_2/2)a_{01}(\phi)\|\delta-\beta_0\|^2+(c_2/2)a_{10}(\phi)\|\beta-\delta_0\|^2,\label{prop.rate1.eq1}
\end{eqnarray}
where $a_{ij}(\phi)=EA_{ij}(\phi)$, for all $0\leq i,j\leq 1$. Moreover, we can also verify that there exists an additional constant $0<c_3\leq c_2/2$ such that $(c_2/2)a_{00}(\phi)\geq c_3$, $(c_2/2)a_{11}(\phi)\geq c_3$, $(c_2/2)\|\delta-\beta_0\|^2\geq c_3$, and $(c_2/2)\|\beta-\delta_0\|^2\geq c_3$, for all $\theta\in\Theta_{\eta}$, decreasing $\eta$ slightly as needed to $\eta_1>0$. Thus, for some $0<c_3,\eta_1<1$,
\begin{eqnarray}
M(\theta)-M(\theta_0)&\geq& c_3\left(\|\beta-\beta_0\|^2+\|\delta-\delta_0\|^2+a_{01}(\phi)+a_{10}(\phi)\right),\label{rate.eq2}
\end{eqnarray}
for all $\theta\in\Theta_{\eta_1}$. The proposition is proved if we can establish that there exists
$0<c_4,\eta_2<\infty$ such that
\begin{eqnarray}
a_{01}(\phi)+a_{10}(\phi)&\geq& c_4\left(\|\omega-\omega_0\|+|\gamma-\gamma_0|\right)\label{rate.eq3}
\end{eqnarray}
for all $\theta\in\Theta_{\eta_2}$.

To this end, define $\nu=(k_1+1)\eta^2$, $(X_0,Z_0)\equiv W_0$ where $W_0$ is as defined in Lemma~\ref{lthm:cr0}, and let $(a)_{+}$ and $(a)_{-}$ be the positive and negative parts of $a\in\re$ respectively. We now have, for all $\theta\in\Theta_{\eta}$, for $\eta>0$ sufficiently small, that
\begin{eqnarray*}
a_{01}(\phi)&=&P[\omega_0'X-\gamma_0\leq 0<\omega'X-\gamma]\\
&=&P[0\leq -\omega_0'X+\gamma_0<(\omega-\omega_0)'X-\gamma+\gamma_0]\\
&=&P[-(\omega-\omega_0)'X+\gamma-\gamma_0<\omega_0'X-\gamma_0\leq 0]\\
&\geq&P[-(\omega-\omega_0)'X+\gamma-\gamma_0<\omega_0'X-\gamma_0<0]\\
&=&P\left[\left.-\nu^{-1}\{(\omega-\omega_0)'X-\gamma+\gamma_0\}<\nu^{-1}U<0\right| -\nu\leq U\leq\nu\right]\cdot P[|U|\leq\nu]\\
&=&\left(P\left[-\nu^{-1}\{(\omega-\omega_0)'X_0-\gamma+\gamma_0\}<U_0<0\right]+o(1)\right)\cdot 2f_0\nu(1+o(1))\\
&=&P\left[-\nu^{-1}\{(\omega-\omega_0)'X_0-\gamma+\gamma_0\}<U_0<0\right]2f_0\nu+o(\nu)\\
&=&\frac{1}{2}E\left(-\nu^{-1}\{(\omega-\omega_0)'X_0-\gamma+\gamma_0\}\right)_{-}\cdot 2f_0\nu+o(\nu)\\
&=&f_0E\left(-(\omega-\omega_0)'X_0+\gamma-\gamma_0\right)_{-}+o(\eta^2),
\end{eqnarray*}
where the first three equalities and the first inequality follow from definitions and straightforward calculations. The first equality after the inequality follows from the use of $U=\omega_0'X-\gamma_0$ as defined in Condition C2 combined with the fact that 
\[\left|(\omega-\omega_0)'X-\gamma+\gamma_0\right|\leq \nu\] 
as a consequence of Condition C1. The next equality follows from Lemma~\ref{lthm:cr0}, the portmanteau theorem (see, e.g., Theorem~7.6 of \citep{kosorok2008introduction}), and Condition C2. The second-to-last equality follows from the independence of $U_0$ and $X_0$ combined with the fact that $U_0$ is uniform on $[-1,1]$ (this also explains the leading 1/2 term). The last equality follows from the definition of $\nu$.

We can use similar arguments to verify that
\[a_{10}(\phi)\;\geq\; f_0E\left(-(\omega-\omega_0)'X_0+\gamma-\gamma_0\right)_{+}+o(\eta^2).\]
Thus
\[a_{01}(\phi)+a_{10}(\phi)\;\geq\;f_0E|(\omega-\omega_0)'X_0-\gamma+\gamma_0|+o(\eta^2).\]
Since $\omega_0\omega_0'+\overline{\omega}_0\overline{\omega}_0'$ equals the identity by definition of $\overline{\omega}_0$, we obtain that
\begin{eqnarray}
(\omega-\omega_0)'X_0&=&(\omega-\omega_0)'\overline{\omega}_0\overline{\omega}_0'X_0+(\omega-\omega_0)'\omega_0\gamma_0\\ \nonumber
&=&(\omega-\omega_0)'\overline{\omega}_0\overline{\omega}_0'X_0+O(\eta^4), \label{new.e1}
\end{eqnarray}
where the second term on the right of the first equality follows from the fact that $\omega_0'X_0-\gamma_0=0$ and the second equality follows from Lemma~\ref{lthm:cr2}. Letting $\tilde{\omega}=\overline{\omega}_0'(\omega-\omega_0)$ and $\tilde{X}_0=\overline{\omega}_0'X_0$, and using the fact that $\tilde{X}_0$ is full rank, we obtain by application of Lemma~\ref{lthm:cr1} that there exists a $c_4>0$ such that
\begin{eqnarray*}
E|(\omega-\omega_0)'X_0-\gamma+\gamma_0|&\geq&E|\tilde{\omega}'\tilde{X}_0-\gamma+\gamma_0|+o(\eta^2)\\
&\geq&(c_4/f_0)(\|\tilde{\omega}\|+|\gamma-\gamma_0|)+o(\eta^2)\\
&\geq&(c_4/f_0)(\|\omega-\omega_0\|+|\gamma-\gamma_0|)+o(\eta^2),
\end{eqnarray*}
where the second inequality follows from the fact that 
\[\|\tilde{\omega}\|^2=|(\omega-\omega_0)'(\overline{\omega}_0\overline{\omega}_0'+\omega_0\omega_0')(\omega-\omega_0)|-O(\eta^4).\]
Thus~(\ref{rate.eq3}) follows. Note that in using Lemma~\ref{lthm:cr1}, the required condition of the lemma can be shown to be satisfied for $W = \tilde X_0$ by contradiction. Suppose there does not exist $\eta > 0$ such that $E|u'(W-\mu)| >\eta$ for all $u\in S^{q-1}$. Then there exists a $u\in S^{q-1}$ such that $E|u'(W-\mu)| = 0$. However, since $W-\mu$ is full rank, this is not possible for any $u\in S^{q-1}$. Thus, the condition of Lemma 5.5 is satisfied.

Combining this with~(\ref{rate.eq2}), we now conclude that for some $0<c_5,\eta_5<1$,
\[M(\theta)-M(\theta_0)\geq c_5\tilde{d}^2(\theta,\theta_0),\]
for all $\{\theta:\;\tilde{d}(\theta,\theta_0)\leq\eta_5\}$, which is the desired conclusion.$\Box$

\end{proof}

\subsection{Proof of Lemmas~\ref{lthm:cr0}--\ref{lthm:cr2}}
\begin{proof}[Proof of Lemma~\ref{lthm:cr0}]
Let $W=(X,Z)$ and $W(u)$ be a random variable with distribution equal to the conditional distribution of $W$ given $U=u$. Let $U_{\nu}$ be a random variable with distribution equal to the conditional distribution of $U$ give $|U|\leq\nu$. Then $(\nu^{-1}U,W)$ given $|U|\leq\nu$ has the same distribution as $(\nu^{-1}U_{\nu},W(U_{\nu}))$. Note that by Condition~C3, $W(u)\leadsto W(0)=W_0\sim G$, as $u\rightarrow 0$. It is also easy to verify by Condition~C2 that $\nu^{-1}U_{\nu}\leadsto U_0$ as $\nu\downarrow 0$. For a complete metric space $(\mathbb{B},d)$, let $BL_1(\mathbb{B})$ be the space of all Lipshitz continuous functions $f:\mathbb{B}\mapsto\re$ such that $|f|\leq 1$ and $\sup_{x,y\in\mathbb{B}:d(x,y)\leq\rho}|f(x)-f(y)|\leq \rho$, for all $\rho>0$. 

We now have that
\begin{eqnarray*}
\lefteqn{\sup_{f\in BL_1(\re^{d+p+1})}\left|Ef\left(\nu^{-1}U_{\nu},W(U_{\nu})\right)-Ef(U_0,W_0)\right|}&&\\
&\leq&\sup_{f\in BL_1(\re^{d+p+1})}\left|Ef\left(\nu^{-1}U_{\nu},W(U_{\nu})\right)-Ef(\nu^{-1}U_{\nu},W_0)\right|\\
&&+\sup_{f\in BL_1(\re^{d+p+1})}\left|Ef(\nu^{-1}U_{\nu},W_0)-Ef(U_0,W_0)\right|\\
&=&A_{\nu}+B_{\nu}.
\end{eqnarray*}
If we can show that $A_{\nu}\rightarrow 0$ and $B_{\nu}\rightarrow 0$, as $\nu\downarrow 0$, the desired conclusion for the lemma will follow from the Portmanteau theorem.

Note that for any $f\in BL_1(\re^{d+p+1})$,
\begin{eqnarray*}
\lefteqn{\left|Ef\left(\nu^{-1}U_{\nu},W(U_{\nu})\right)-Ef(\nu^{-1}U_{\nu},W_0)\right|}&&\\
&\leq&E\left|E\left[\left. f\left(\nu^{-1}U_{\nu},W(U_{\nu})\right)-f(\nu^{-1}U_{\nu},W_0)\right|U_{\nu}\right]\right|\\
&\leq&E\left\{\sup_{g\in BL_1(\re^{d+p})}\left|E\left[\left. g(W(U_{\nu}))-g(W(0))\right|U_{\nu}\right]\right|\right\}\\
&\rightarrow& 0,
\end{eqnarray*}
as $\nu\downarrow 0$, by the fact that $W(U_{\nu})\leadsto W(0)$ combined with the Portmanteau theorem. Since the right-hand-side of the last inequality does not depend on $f$, we obtain that $A_{\nu}\rightarrow 0$.

Note also that for any $f\in BL_1(\re^{d+p+1})$, 
\begin{eqnarray*}
\lefteqn{\left|Ef(\nu^{-1}U_{\nu},W_0)-Ef(U_0,W_0)\right|}&&\\
&\leq&E\left|E\left[\left. f(\nu^{-1}U_{\nu},W_0)-f(U_0,W_0)\right|W_0\right]\right|\\
&\leq&E\left\{\sup_{g\in BL_1(\re)}\left|E\left[\left. g(\nu^{-1}U_{\nu})-g(U_0)\right| W_0\right]\right|\right\}\\
&\rightarrow& 0,
\end{eqnarray*}
as $\nu\downarrow 0$, using similar arguments as before combined with the fact that $W_0$ is independent of both $\nu^{-1}U_{\nu}$ and $U_0$. Since, once again, the right-hand-side of the last inequality does not depend on $f$, we obtain that $B_{\nu}\rightarrow 0$. Hence the desired conclusion of the lemma now follows.
\end{proof}

\begin{proof}[Proof of Lemma~\ref{lthm:cr1}]
Suppose $\displaystyle \frac{|a'\mu+b|}{\|a\|} \geq
\frac{\eta}{2}$. Then,
\begin{eqnarray*}
E|a'W+b|&=& \|a\| E \left| \frac{a'(W-\mu)}{\|a\|}
+\frac{a'(\mu+b)}{\|a\|} \right| \\
& \geq & \|a\| \left| E \left (
\frac{a'(W-\mu)}{\|a\|}+\frac{a'(\mu+b)}{\|a\|}\right)
\right|~~ \text{ (by Jensen's Inequality)}\\
&=& \|a\| \left|\frac{a'(\mu+b)}{\|a\|} \right| \\
& \geq & \frac{\eta}{2} \|a\|.
\end{eqnarray*}

\noindent Now suppose $\displaystyle \frac{|a'\mu+b|}{\|a\|} <
\frac{\eta}{2}$. Then,
\begin{eqnarray*}
E|a'W+b|&=& \|a\| E \left| \frac{a'(W-\mu)}{\|a\|}
+\frac{a'(\mu+b)}{\|a\|} \right| \\
&\geq & \|a\| \left( E \left| \frac{a'(W-\mu)}{\|a\|}\right| -
\left| \frac{a'(\mu+b)}{\|a\|} \right|\right) \\
& \geq & \|a\| \left( \eta-\frac{\eta}{2}\right) \\
& \geq & \frac{\eta}{2} \|a\|.
\end{eqnarray*}

\noindent Therefore, for all $a \in \mathbb{R}^q$ and $b \in \mathbb{R}$,
\begin{equation}
\label{eqn:cr2}
E|a'W+b| \geq \frac{\eta}{2} \|a\|.
\end{equation}

\noindent Now suppose $\mu=0$. Then
\begin{equation*}
E|a'W+b| \geq |E(a'W+b)| = ~| b | ~~\text{ (by Jensen's
Inequality)}.
\end{equation*}
\noindent Hence, combining with (\ref{eqn:cr2}), we have
\begin{equation}
\label{eqn:cr3} E|a'W + b| \geq (\frac{\eta}{2} \wedge 1)
\max{(\|a\|, |b|)}.
\end{equation}

\noindent Next, suppose $\mu \neq 0$. First, assume that $\|b\| > 2
\|\mu\|\|a\|$. Then,
\begin{eqnarray*}
E|a'W+b| &\geq & |b|-|a'\mu| \\
& \geq & |b| - \|a\|\|\mu\| ~~( |a'\mu| \leq \|a\|\|\mu\|) \\
& \geq & \frac{|b|}{2}.
\end{eqnarray*}
\noindent Therefore,
\begin{equation}
\label{eqn:cr3b} E|a'W+b| \geq \frac{\max{(|b|,
2\|\mu\|\|a\|)}}{2}.
\end{equation}

\noindent Now assume that $\|b\| \leq 2 \|\mu\|\|a\|$. Then, from
(\ref{eqn:cr2}), we have
\begin{eqnarray*}
E |a'W+b| &\geq  & \frac{\eta}{2} \|a\|  \\
& = & \frac{\eta}{2} \max{( \|a\|, \frac{|b|}{2 \|\mu\|})} \\
& = & \frac{\eta}{2} \frac{1}{2 \|\mu\|} \max { (|b|, 2 \|\mu\|\|a\|)}.
\end{eqnarray*}

\noindent Combining with (\ref{eqn:cr3b}), we have
\begin{eqnarray*}
E|a'W+b| & \geq & \left( \frac{1}{2} \wedge \frac{\eta}{4
\|\mu\|} \right)  \max {(|b|, 2\|\mu\|\|a\|)} \\
& \geq & \left( \frac{1}{2} \wedge \frac{\eta}{4 \|\mu\|} \right )
(1 \wedge 2 \|\mu\|) \max{ (\|a\|, |b|)}.
\end{eqnarray*}
\noindent Since $\max{(\|a\|, |b|)} \geq \frac{1}{2} (\|a\|+|b|)$,
we now have, combining everything together, that
\begin{equation}
E|a'W+b| \geq c_{\ast} ( \|a\|+|b|),\label{e1.lower}
\end{equation}
\noindent where
 \[ c_{\ast} =
 \begin{cases}
 ~\frac{1}{2} ( \eta/2 \wedge 1) , &\text{if~} \mu =0,\\
 ~\frac{1}{2} \left(  \frac{1}{2} \wedge \frac{\eta}{4\|\mu\|}
\right)(1 \wedge 2\|\mu\|),  &\text{if~} \mu \neq 0.
 \end{cases}
 \]
 
 Now suppose, for some $(a,b)\in\re^q\times\re:\;\|a\|^2+b^2=1$, that $E\left(B|a'W-b|\right)=0$. Then \newline $E\left(\left.|a'W-b|\,\right|B>0\right)=0$ since $P[B\leq 0]=0$. But this implies $E|a'W-b|=0$ which is a contradiction of~(\ref{e1.lower}). Thus there exists a $c_{\ast\ast}>0$ for which $E\left(B|a'W-b|\right)\geq c_{\ast\ast}$ for all $(a,b)\in\re^q\times\re:\;\|a\|^2+b^2=1$, and hence the desired result follows with $c=c_{\ast}\wedge (c_{\ast\ast}/2)$.
\end{proof}

\begin{proof}[Proof of Lemma~\ref{lthm:cr2}]
Note that the points $0$, $\omega_1$ and $\omega_2$ form an isosceles triangle with a base corresponding to the line segment from $\omega_1$ to $\omega_2$. This triangle divides at the midpoint of this base into two identical right triangles which are mirror images of each other, with the line from 0 to $\omega_1$ forming the hypotenuse of one of them and the line from 0 to $\omega_2$ forming the hypotenuse of the other. This is illustrated in Figure~\ref{fig:partangle1}. The angle $\theta$ between the base and the hypotenuse of either right triangle is the same as the angle between $\omega_1$ and $\omega_2-\omega_1$ and satisfies $\cos\theta=\|\omega_2-\omega_1\|/2$ since both $\omega_1$ and $\omega_2$ are in $S^{p-1}$. Hence
\[|(\omega_2-\omega_1)'\omega_1|=\|\omega_2-\omega_1\|\cdot \|\omega_1\|\cdot \cos\theta=\frac{\|\omega_2-\omega_1\|^2}{2}.\]

\begin{figure}[htp]
     \centering
     \includegraphics[angle=0,width=3in]{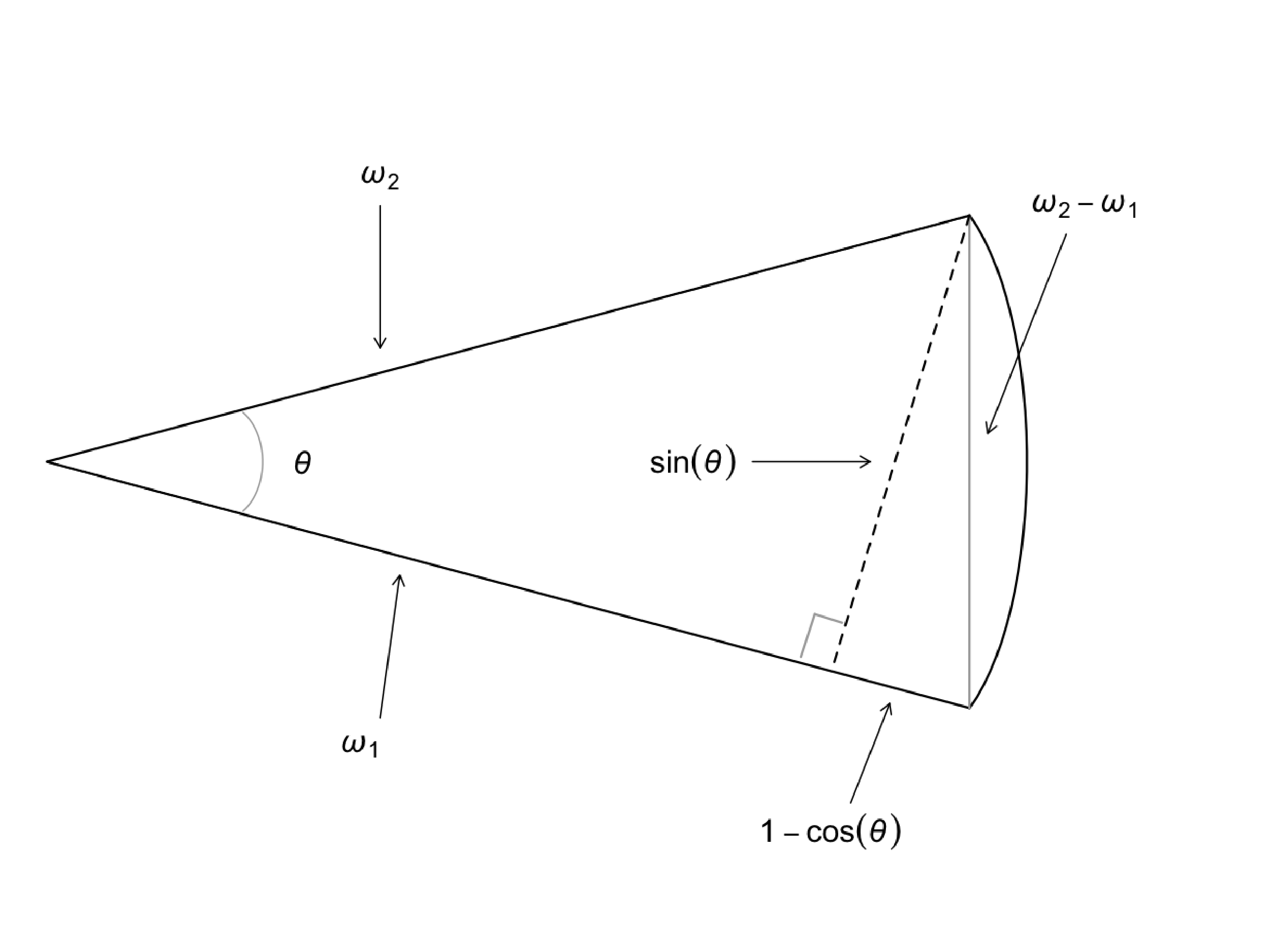}
\caption{Partial of the unit circle for the rate of convergence. $\theta$ is the angle between $\omega_1$  and $\omega_2$ which should be small.}
     \label{fig:partangle1}
\end{figure}
\end{proof}
\FloatBarrier

\section{Proofs of Section \ref{sec:WeakConvergence}} 
\label{asec:WeakConvergence}
\subsection{Proof of the preliminary results in Section \ref{sec:WeakConvergence}}
\label{asec:WeakConvergence_prelim}

\subsubsection{Lemma \ref{limiting.properties}}
\label{asec:limiting.properties}

In the following lemma, we provide several useful properties, such as convexity and boundedness, of the limiting functions and sets.
\begin{lemma}\label{limiting.properties}
The following are true:
\begin{enumerate}
\item $\tilde{\Phi}_0=\left\{(g_1,g_2)\in\re^{p-1}\times\re: \tilde{C}_R(g_1)>0\;\mbox{and}\;\tilde{C}_L(g_1)\leq g_2<\tilde{C}_U(g_1)\right\}$,
\item $\tilde{C}_U$, $-\tilde{C}_L$, and $\tilde{C}_R$ are concave functions almost surely, and
\item $\tilde{R}_0$ is convex and bounded almost surely.
\end{enumerate}
\end{lemma}

\begin{proof}[Proof of Lemma\ref{limiting.properties}]
By the definitions, we know that $(g_1,g_2)\in\tilde{\Phi}_0$ if and only if $g_1'\tilde{X}_j^{-}-g_2-\tilde{U}_j^{-}\leq 0$ for all $1\leq j\leq \tilde{J}_0^{-}:V_j^{-}=-1$; $g_1'\tilde{X}_j^{+}-g_2+\tilde{U}_j^{+}\leq 0$ for all $1\leq j\leq \tilde{J}_0^{+}:V_j^{+}=-1$; $g_1'\tilde{X}_j^{-}-g_2-\tilde{U}_j^{-}> 0$ for all $1\leq j\leq \tilde{J}_0^{-}:V_j^{-}=1$; and $g_1'\tilde{X}_j^{+}-g_2+\tilde{U}_j^{+}> 0$ for all $1\leq j\leq \tilde{J}_0^{+}:V_j^{+}=1$. Thus $\tilde{\Phi}_0$ is precisely the set of $(g_1,g_2)\in\re^{p-1}\times\re$ for which $\tilde{C}_L(g_1)-g_2\leq 0$ and $\tilde{C}_U(g_1)-g_2>0$. This now establishes Part 1. Recall that linear functions are both concave and convex. Since $\tilde{C}_U$ is the minimum of a finite number of concave functions, it is also concave; since $\tilde{C}_L$ is the maximum of a finite number of convex functions, it is also convex, and thus its negative is concave. Since $\tilde{C}_R$ is the sum of two concave functions, it is also concave. This establishes Part 2. The boundedness of $\tilde{R}_0$ follows from Part 1 of the lemma and also Part 1 of Theorem~\ref{argmax.compactness}. Suppose $g_{11}$ and $g_{12}$ are both in $\tilde{R}_0$, then $\tilde{C}_R(g_{11})>0$ and $\tilde{C}_R(g_{12})>0$. By the concavity of $\tilde{C}_R$, we have for any $\alpha\in(0,1)$, that $\tilde{C}_R(\alpha g_{11}+(1-\alpha) g_{12})\geq \alpha\tilde{C}_R(g_{11})+(1-\alpha)\tilde{C}_R(g_{12})>0$, and Part 3 is established.$\Box$
\end{proof}

\subsubsection{Proof of Lemma \ref{lemma.onetoone}}
\label{asec:lemma.onetoone}
\begin{proof}[Proof of Lemma\ref{lemma.onetoone}]
Continuity of both $h_{\ast n}$ and $h_3\mapsto \overline{\omega}_0'h_3$ is obvious. Suppose $h_{\ast n}(g_{11})=h_{\ast n}(g_{12})$ for $g_{11},g_{12}\in A_n$. Then, using the orthogonality of $\overline{\omega}_0$ and $\omega_0$, we have
\begin{eqnarray*}
0&=&\|h_{\ast n}(g_{11})-h_{\ast n}(g_{12})\|^2\\
&=&\|g_{11}-g_{12}\|^2+n^2\left(\sqrt{1-\frac{\|g_{11}\|^2}{n^2}}-\sqrt{1-\frac{\|g_{12}\|^2}{n^2}}\right)^2,
\end{eqnarray*}
which implies that $g_{11}=g_{12}$. Now, for $g_1\in A_n$, we can use the monotone increasingness property of $u\mapsto 1-\sqrt{1-u}$ (over positive $u$) to verify that
\begin{eqnarray*}
\|h_{\ast n}(g_1)\|^2&=&\|g_1\|^2+n^2\left(1-\sqrt{1-\frac{\|g_1\|^2}{n^2}}\right)^2\\
&\leq&n^2+n^2,
\end{eqnarray*}
and thus $h_{\ast n}(g_1)\in B_n$ for all $g_1\in A_n$. Now let $h_3\in B_n$, and note that $\overline{\omega}_0'h_3$ is the projection of $h_3$ onto the orthocomplement of $\omega_0$ (i.e., the linear span of $\overline{\omega}_0$ which equals $\re_{\overline{\omega}_0}^{p-1}$). Using the fact that $\omega_0+h_3/n\in S^{p-1}$, along with arguments similar to those used in Lemma~\ref{lthm:cr2}, we can now establish that the length of this projection is equal to 
\[n\sqrt{\frac{\|h_3\|^2}{n^2}-\frac{\|h_3\|^2}{2n^2}}=\frac{\|h_3\|}{\sqrt{2}}\leq n,\]
and thus $\overline{\omega}_0'h_3\in A_n$. Note also that for any $g_1\in A_n$, $\overline{\omega}_0'h_{\ast n}(g_1)=g_1$, and thus $h_{\ast n}^{-1}$ is the inverse of $h_{\ast n}$. Moreover, for any $h_3\in B_n$, 
\begin{eqnarray*}
h_{\ast n}(\overline{\omega}_0'h_3)&=&\overline{\omega}_0\overline{\omega}_0'h_3-n\left(1-\sqrt{1-\frac{(\overline{\omega}'h_3)^2}{n^2}}\right)\omega_0,\\
&\equiv&C_1 h_3 +c_{2n}\omega_0,
\end{eqnarray*}
where $C_1 h_3$ is the projection of $h_3$ onto the orthocomplement of $\omega_0$ and $c_{2n}$, and after recycling previous arguments, can be shown to be the length and sign of the projection of $h_3$ onto $\omega_0$. Thus $h_{\ast n}(\overline{\omega}_0'h_3)=h_3$, and $h_{\ast n}$ is also the inverse of $h_{\ast n}^{-1}$.
Combining the previous results, we now have that both $h_{\ast n}$ and $h_{\ast n}^{-1}$ as given are one to one and onto.
\end{proof}

\subsubsection{Theorem \ref{convergence-modemid}}
\label{asec:convergence-modemid}
We develop the compact convergence of the mode-midpoint estimator in the theorem below needed to prove Theorem~\ref{convergence-modemid}.  We begin with the following lemma:
\begin{lemma}\label{compact.modemid.def}
Assume Conditions {C1, C2, C3', C4 and C5} hold. Then the following hold:
\begin{enumerate}
\item When $F_n(k)=1$, $h_3\mapsto C_R^{nk}(h_3)$ has a unique maximum over $R_n(k)$ which we will denote $\check{h}_{3n}(k)$.
\item $(\check{h}_{3n}(k),\check{h}_{4n}(k))$, where $\check{h}_4(k)=\left[C_L^{nk}(\check{h}_{3n}(k))+C_L^{nk}(\check{h}_{3n}(k))\right]/2$, is contained in $\Phi_n'(k)$, for all $n$ large enough.
\item When $F_0(k)=1$, $g_1\mapsto\tilde{C}_R^k(g_1)$ over $g_1\in\tilde{R}_0(k)$ is concave and has a unique maximum which we will denote $\check{g}_1(k)$.
\item $(\check{g}_1(k),\check{g}_2(k))$, where $\check{g}_2(k)=\left[\tilde{C}_L^k(\check{g}_1(k))+ \tilde{C}_U^k(\check{g}_1(k))\right]/2$, is contained in $\tilde{\Phi}_0(k)$. 
\end{enumerate}
\end{lemma}

\begin{proof} [Proof of Lemma \ref{compact.modemid.def}]
When $F_n(k)=1$, $C_R^{nk}$ is well-defined, and the set $R_n^{\ast}(k)=\{h_3:C_R^{nk}(h_3)\geq C_R^{nk}(\tilde{h}_3(k))$ is nonempty, closed and geodesically connected, using arguments similar to those used in Lemma~\ref{meanmidexistence}. Hence the set $R_n^{\ast\ast}(k)=\{h_3\in R_n^{\ast}(k):\;\|h_3\|^2\leq (0.9k_1^{-1}k)^2+r_n^2(k)\}$ is also geodesically connected and closed and must contain the set of all $\arg\max$ of $C_R^{nk}$ over $R_{n}(k)$. Hence the set $\omega_0+R_n^{\ast\ast}(k)/n\subset S^{p-1}$ is also geodesically closed and connected. When $F_n(k)=1$, all of the conditions of Lemma~\ref{modemidexistence.unique} are satisfied, and thus $\check{h}_{3n}(k)$ is well-defined and the unique maximizer. Since $\check{h}_{3n}(k)$ is the $\arg\max$ of a function over a closed subset of $R_n(k)$, $\|\check{h}_{3n}(k)\|\leq k_1^{-1}k$ for all $n$ large enough, and thus 
$C_L^{nk}(\check{h}_{3n}(k))\leq k$ and $C_U^{nk}(\check{h}_{3n}(k))\leq k$,
and hence also $|\check{h}_{4n}(k)|\leq k$, for all $n$ large enough. By recycling arguments used in the proof of Lemma~\ref{meanmidexistence}, we obtain the desired containment in $\Phi_n'(k)$. By recycling previous arguments, we can verify that when $F_0(k)=1$, the minimums and maximums in the definitions of the maps $\tilde{C}_L^k$, $\tilde{C}_U^k$ and $\tilde{C}_R^k$ are over non-null sets and are thus well defined, and we can also readily verify the concavity of $\tilde{C}_R^k$. For $\tilde{C}_R^k$ to not have a unique maximum would require $h_3\mapsto\tilde{C}_L^k(h_3)$ and $h_3\mapsto\tilde{C}_U^k(h_3)$ to be parallel for at least two distinct points in $\tilde{R}_0(k)$, but this is impossible since the $\tilde{X}_j^{-}$ and $\tilde{X}_j^{+}$ values involved are continuously distributed and thus the probability of this happening is zero. Hence Part~(3) follows. Part~(4) follows by recycling previous arguments.
\end{proof}

\subsubsection{Theorem \ref{compact.convergence.mode}}
\label{asec:compact.convergence.mode}
Denote $\check{h}_n(k)=(\tilde{h}_{1n}(k),\tilde{h}_{2n}(k),\check{h}_{3n}(k),\check{h}_{4n}(k))$ and $\check{h}_0(k)=(\tilde{h}_{10}(k),\tilde{h}_{20}(k),\overline{\omega}_0\check{g}_1(k),\check{g}_2(k))$. 

Now we are ready for the following convergence theorem for restrictions over compact sets:
\begin{theorem}\label{compact.convergence.mode}
Assume conditions {C1, C2, C3', C4 and C5}. Then, for every $k\in (0,\infty)$,
\[\left(\begin{array}{c}F_n(k)\check{h}_n(k)\\ F_n(k)\end{array}\right)\leadsto \left(\begin{array}{c}F_0(k)\check{h}_0(k)\\ F_0(k)\end{array}\right),\]
where $\check{h}_n(k)=(\tilde{h}_{1n}(k),\tilde{h}_{2n}(k),\check{h}_{3n}(k),\check{h}_{4n}(k))$ and $\check{h}_0(k)=(\tilde{h}_{10}(k),\tilde{h}_{20}(k),\overline{\omega}_0\check{g}_1(k),\check{g}_2(k))$.
\end{theorem}
This gives us enough results to prove Theorem~\ref{convergence-modemid}.

{\bf Proof of Theorem~\ref{compact.convergence.mode}}. We can utilize essentially all of the results of the proof of Theorem~\ref{compact.convergence} above as Condition C3' is stronger than Condition C3. We will assume that we are in the new probability space define in that theorem and that $F_0(k)=1$ and $n$ is large enough (as needed). Lemma~\ref{compact.modemid.def} tells us that $h_3\mapsto C_R^{nk}(h_3)$ has a unique maximum over $h_3\in R_n(k)$ and hence also $g_1\mapsto \tilde{C}_R^{nk}(g_1)$ has a unique maximum over $g_1\in\tilde{R}_n(k)$ for all $n$ large enough. This same lemma also reveals that $g_1\mapsto\tilde{C}_R^k(g_1)$ has a unique maximum over $g_1\in \tilde{R}_0(k)$. Since we have previously established uniform convergence of $\tilde{C}_R^{nk}$ to $\tilde{C}_R^k$, we now have---after combining with and recycling previous arguments---that $\check{h}_{3n}(k)\rightarrow \overline{\omega}'\check{g}_1(k)$, as $n\rightarrow\infty$; and, consequently, also $\check{h}_{4n}(k)\rightarrow \check{g}_2(k)$. Thus $\check{h}_n(k)\rightarrow \check{h}_0(k)$, as $n\rightarrow\infty$, and the desired weak convergence follows.$\Box$



\subsection{Proof of Theorem \ref{convergence-modemid}}
\label{asec:proof.convergence-modemid}
\begin{proof}[Proof of Theorem \ref{convergence-modemid}]
Let $F\subset\re^{2d+p+1}$ be closed. As before, we need to show that  for any $\eta>0$, there exists a $0<k<\infty$ such that $\lim\inf_{n\rightarrow\infty}P(F_n(k)=1)\geq 1-\eta$, $P(F_0(k)=1)\geq 1-\eta$, $\lim\inf_{n\rightarrow\infty}P(F_n(k)\check{h}_n(k)=\red{D_n}(\check{\theta}_n-\theta_0))\geq 1-\eta$,
and $P(F_0(k)\check{h}_0(k)=\check{h}_0)\geq 1-\eta$. The first two of these follow since these were already shown under Condition C3. 

Now let $\Phi_n'$ be the level set containing $\red{n}(\check{\phi}_n-\phi_0)$. As previously, the entirety of $\Phi_n'$ is simultaneously asymptotically bounded in probability as a result of Theorem~\ref{theorem:rateofconvergence1}. Hence, as before, for any $\eta>0$, there exists a $k<\infty$ such that both $\lim\inf{n\rightarrow\infty}P(F_n(k)=1)\geq 1-\eta/2$ and $\lim\inf_{n\rightarrow\infty}P(\red{D_n}(\check{\theta}_n-\theta_0)\in H_n(k))\geq 1-\eta/2$, and thus
$P(F_n(k)\check{h}_n(k)=\red{D_n}(\check{\theta}_n-\theta_0))\geq 1-\eta$. We can apply similar arguments for $\check{h}_0$ via Theorem~\ref{argmax.compactness} to obtain that for any $\eta>0$ there exists a $k<\infty$ such that $P(F_0(k)\check{h}_0(k)=\check{h}_0)\geq 1-\eta$, and our required probability bounds for all $\eta>0$ are satisfied.

Hence, fixing an $\eta>0$  and finding a corresponding $k<\infty$ which simultaneously satisfy these criteria, we have as before that
\begin{eqnarray*}
\lefteqn{\lim\sup_{n\rightarrow\infty}P(\red{D_n}(\check{\theta}_n-\theta_0)\in F)}&&\\
&\leq&\lim\sup_{n\rightarrow\infty}P(\red{D_n}(\check{\theta}_n-\theta_0)\in F,\;F_n(k)=1,\;F_n(k)\check{h}_n(k)=\red{D_n}(\check{\theta}_n-\theta_0))\\
&&+\lim\sup_{n\rightarrow\infty}P(F_n(k)=0)+\lim\sup_{n\rightarrow\infty}P(F_n(k)\check{h}_n(k)\neq \red{D_n}(\check{\theta}_n-\theta_0))\\
&\leq&\lim\sup_{n\rightarrow\infty}P(F_n(k)\check{h}_n(k)\in K)+2\eta\\
&\leq&P(F_0(k)\check{h}_0(k)\in K)+2\eta\\
&\leq&P(\check{h}_0\in K,\;F_0(k)=1,\;F_0(k)\check{h}_0(k)=\check{h}_0)\\
&&+P(F_0(k)=0)+P(F_0(k)\check{h}_0(k)\neq \check{h}_0)+2\eta\\
&\leq&P(\check{h}_0\in K)+4\eta,
\end{eqnarray*}
which, as before, implies the desired week convergence because of the arbitrariness of $\eta>0$ and $K$, followed by the Portmanteau Theorem.
\end{proof}

\subsection{Proof of Theorem~\ref{argmax.compactness} and its related lemma}
\label{asec:WeakConvergence_compact}
{\bf Proof of Theorem~\ref{argmax.compactness}}.

Note that
\begin{eqnarray*}
Q_{02}(g)&=&\int_0^{\infty}\left[B(t)\ind\{g_1'\tilde{X}(t)-g_2>t\}\tilde{E}^{-}(t)\right.\\
&&\left. +(1-B(t))\ind\{-g_1'\tilde{X}(t)+g_2\geq t\}\tilde{E}^{+}(t)\right]dN(t)\\
&\equiv&\tilde{Q}_{02}(g),
\end{eqnarray*}
where $t\mapsto N(t)$ is a homogeneous Poisson process on $[0,\infty)$ with intensity $2f_0$, $B(t)$ is a white-noise type Bernoulli random variable with success probability $1/2$, and where $\tilde{X}(t)$, $\tilde{E}^{-}(t)$ and $\tilde{E}^{+}(t)$ are also white-noise type stochastic processes which we will define shortly. By ``white-noise type'' we mean that there is a new independent draw for every distinct value of $t$. In the given process, these random draws only need to occur at jump times in $N$. Define
$\tilde{E}^{-}(t)$ and $\tilde{E}^{+}(t)$ as follows:
\begin{eqnarray*}
\tilde{E}^{-}(t)&=&\left[(\beta_0-\delta_0)'\tilde{Z}(t)\right]^2+2\epsilon(t)(\beta_0-\delta_0)'\tilde{Z}(t),\;\mbox{and}\\
\tilde{E}^{+}(t)&=&\left[(\beta_0-\delta_0)'\tilde{Z}(t)\right]^2-2\epsilon(t)(\beta_0-\delta_0)'\tilde{Z}(t),
\end{eqnarray*}
where $\epsilon(t)$ is a white-noise type process with distribution the same as the residual in model~(\ref{eq:dmodel}) and $(\tilde{Z}(t),\tilde{X}(t))$ is joint white-noise type process with joint distribution equal to $G$. What we are doing is combining the Poisson process associated with the $\tilde{U}_j^{-}$ values with the Poisson process associated with the $\tilde{U}_j^{+}$ values. The combined process will have intensity equal to the sum of the two constituent intensities which, in this case, is $2f_0=f_0+f_0$. Because of the independence of these two processes and the equality of their intensities, the identity of the two constituent processes can be recaptured through the Bernoulli white-noise type process $B(t)$.

Part 1 of the theorem will follow if we can establish that
\begin{eqnarray}
\lim_{k\rightarrow\infty}P\left\{\inf_{g\in\re^{p-1}\times\re:\;\|g_1\|\vee|g_2|>k}\tilde{Q}_{02}(g)\leq 0\right\}&=&0,\label{eq:lcompact1} 
\end{eqnarray}
since we already know that $\tilde{Q}_{02}(0)=0$, and thus an element $\tilde{g}$ of the $\arg\min$ of $g\mapsto\tilde{Q}_{02}(g)$ can have $\|g_1\|\vee|g_2|>k$ when the event in the given probability statement happens. Once we establish Part 1, Part 3 will follow since for all $\tilde{U}_j^{-}>T_0(k_1+1)$, $\ind\{\tilde{g}_1'\tilde{X}_j^{-}-\tilde{g}_2-\tilde{U}_j^{-}\leq 0\}=1$, and also $\ind\{\tilde{g}_1'\tilde{X}_j^{+}-\tilde{g}_2+\tilde{U}_j^{+}> 0\}=1$ for all $\tilde{U}_j^{+}>T_0(k_1+1)$. Part 2 will now follow from the positive expectation and continuity of the $\tilde{E}_j^{-}$ and $\tilde{E}_j^{+}$ random variables, using arguments similar to those used in the proof of Theorem~\ref{compact.convergence}.  Thus the proof will be complete once we establish~(\ref{eq:lcompact1})).

To this end, for each $0<k<\infty$ and each integer $j\geq 1$, define $A_j(k)=\{g\in\re^{p-1}\times\re:\; kj^2<\|g_1\|\vee |g_2|\leq k(j+1)^2\}$, and note that
\[P\left\{\inf_{g\in\re^{p-1}\times\re:\;\|g_1\|\vee|g_2|>k}\tilde{Q}_{02}(g)\leq 0\right\}\leq \sum_{j=1}^{\infty}
P\left\{\inf_{g\in A_j(k)}\tilde{Q}_{02}(g)\leq 0\right\}.\]
However, since $\|g_1\|\vee|g_2|\leq k(j+1)^2$ for all $g\in A_j$, the restriction of $g\mapsto \tilde{Q}_{02}(g)$ to $g\in A_j$ involves integrating with respect to $dN(t)$ only up to $t=t_j(k)\equiv (k_1+1)k(j+1)^2$, since the involved indicator functions will all be zero for $t>t_j(k)$. Let $N_j(k)=\int_0^{t_j(k)}dN(t)$, and note that this is, by definition, a Poisson random variable with intensity $2f_0 t_j(k)$. Now, conditional on $N_j(k)=m$, and applying the restriction $g\in A_j(k)$, we have that
\begin{eqnarray*}
\tilde{Q}_{02}(g)&=&\sum_{i=1}^m \left[B_i\ind\{g_1'\tilde{X}_i^{\ast}-g_2>U_i^{\ast}\}\tilde{E}_i^{\ast-}\right.\\
&&\left.+(1-B_i)\ind\{-g_1'\tilde{X}_i^{\ast}+g_2\geq U_i^{\ast}\}\tilde{E}_i^{\ast+}\right],
\end{eqnarray*}
where $(B_i,\tilde{E}_i^{\ast-},\tilde{E}_i^{\ast+})$ are i.i.d. realizations of $(B(t),\tilde{E}^{-}(t),\tilde{E}^{+}(t))$ and $(\tilde{Z}_i,\tilde{X}_i^{\ast},\tilde{\epsilon}_i)$ are i.i.d. realizations of $(\tilde{Z}(t),\tilde{X}(t),\epsilon(t))$. Also, the $U_i^{\ast}$ values are independent of the other random variables and have a uniform distribution on the interval $[0,t_j(k)]$.

Define, for $i\geq 1$, the random quantity $\tilde{Y}_i=(B_i,\tilde{Z}_i,\tilde{X}_i^{\ast},\tilde{\epsilon}_i,U_i^{\ast})$, where the generic version of these random variables are denoted by omitting the $i$ subscript, and also define the class of functions 
\begin{eqnarray*}
{\cal F}_j(k)&=&\left\{B\ind\{g_1'\tilde{X}^{\ast}-g_2>U^{\ast}\}e^{-}(\tilde{Z},\tilde{\epsilon})\right.\\
&&\left. +(1-B)\ind\{-g_1'\tilde{X}^{\ast}+g_2\geq U^{\ast}\}e^{+}(\tilde{Z},\tilde{\epsilon}):\; g\in A_j(k)\right\},\;\mbox{where}\\
e^{-}(\tilde{Z},\tilde{\epsilon})&=&\left[(\beta_0-\delta_0)'\tilde{Z}\right]^2+2\tilde{\epsilon}(\beta_0-\delta_0)'\tilde{Z}\;\mbox{and}\\
e^{+}(\tilde{Z},\tilde{\epsilon})&=&\left[(\beta_0-\delta_0)'\tilde{Z}\right]^2-2\tilde{\epsilon}(\beta_0-\delta_0)'\tilde{Z}.
\end{eqnarray*}
Let $\mathbb{G}_m=m^{1/2}(\ep_m-P)$ be the empirical process associated with a sample of these random variables of size $m$, where $P$ depends on the context, and where a sum of zero elements is zero and thus $\mathbb{G}_0=0$ almost surely. Specifically, when we apply the empirical process to the class ${\cal F}_j(k)$, the distribution is the appropriate one defined by the above random variables with the $U_i^{\ast}$ being uniform on $[0,t_j(k)]$. We now have the following lemma of which proof can be found at the end of this subsection:
\begin{lemma}\label{entropy.control}
There exist constants $0<c_{\ast 1},c_{\ast 2}<\infty$ such that for every $j,m\geq 1$ and $0<k<\infty$, $\|\;\|\mathbb{G}_m\|_{{\cal F}_j(k)}^{\ast}\|_{P,2}\leq c_{\ast 1}$ and $\inf_{f\in {\cal F}_j(k)}Pf\geq c_{\ast 2}$, where superscript ${\ast}$ here denotes outer expectation, subscript $P,2$ indicates the $L_2(P)$ norm, and where $P$ depends on $j,k$ as specified above.
\end{lemma}

Continuing with the proof of Theorem~\ref{argmax.compactness}, fix $0<k<\infty$, and note that for each $j\geq 1$,
\begin{eqnarray*}
\lefteqn{P\left\{\inf_{g\in A_j(k)}\tilde{Q}_{02}(g)\leq 0\right\}}&&\\
&\leq&P\left\{\inf_{g\in A_j(k)}\tilde{Q}_{02}(g)\leq 0,\;N_j(k)\geq f_0t_j(k)\right\}+P\{N_j(k)<f_0 t_j(k)\}\\
&\leq&P\left\{\inf_{f\in{\cal F}_j(k)}\left(N_j(k)^{1/2}\mathbb{G}_{N_j(k)}f+N_j(k)Pf\right)\leq 0,\;N_j(k)\geq f_0t_j(k)\right\}\\
&&+P\{-(N_j(k)-2f_0t_j(k))>f_0t_j(k)\}\\
&\leq&E\left[P\left\{\left.N_j(k)^{1/2}\|\mathbb{G}_{N_j(k)}^{\ast}\|_{{\cal F}_j(k)}\geq N_j(k)\inf_{f\in {\cal F}_j(k)} Pf\right|\;N_j(k)\geq f_0 t_j(k) \right\}\right]+\frac{2f_0t_j(k)}{f_0^2t_j^2(k)}\\
&\leq&E\left\{\left. \frac{N_j(k) c_{\ast 1}}{N_j^2(k)\left[\inf_{f\in {\cal F}_j}Pf\right]^2}\right|\;N_j(k)\geq f_0 t_j(k)\right\}+\frac{2}{kf_0(k_1+1)(j+1)^2}\\
&\leq&\frac{c_{\ast 1}}{kc_{\ast 2}f_0(k_1+1)(j+1)^2}+\frac{2}{kf_0(k_1+1)(j+1)^2}\\
&=&\frac{c_{\ast 1}+2c_{\ast 2}}{kc_{\ast 2}f_0 (k_1+1)(j+1)^2},
\end{eqnarray*}
where the Markov inequality and  Lemma~\ref{entropy.control}  were utilized in the third and the fourth inequalities, respectively.
This now gives as that 
\[\lim_{k\rightarrow\infty}\sum_{j=1}^{\infty}
P\left\{\inf_{g\in A_j(k)}\tilde{Q}_{02}(g)\leq 0\right\}=0,\]
since $\sum_{j=1}^{\infty}(j+1)^{-2}<\infty$, and thus~(\ref{eq:lcompact1}) follows, completing the proof.$\Box$

\begin{proof} [Proof of Lemma \ref{entropy.control}]
Define ${\cal F}_j^1(k)=\left\{B\ind\{g_1'\tilde{X}^{\ast}-g_2>U^{\ast}\}e^{-}(\tilde{Z},\tilde{\epsilon}):\; g\in A_j(k)\right\}$ and \newline ${\cal F}_j^2(k)=\left\{(1-B)\ind\{-g_1'\tilde{X}^{\ast}+g_2\geq U^{\ast}\}e^{+}(\tilde{Z},\tilde{\epsilon}):\; g\in A_j(k)\right\}$, and note that \newline $\|\;\|\mathbb{G}_m\|_{{\cal F}_j(k)}^{\ast}\|_{P,2}\leq \|\;\|\mathbb{G}_m\|_{{\cal F}_j^1(k)}^{\ast}\|_{P,2}+\|\;\|\mathbb{G}_m\|_{{\cal F}_j^2(k)}^{\ast}\|_{P,2}$. Note that $\{g_1'\tilde{X}^{\ast}-g_2-U^{\ast}:\; g\in A_j(k)\}$ is a vector space of dimension $p$, and thus, by Lemma~9.6 of \cite{kosorok2008introduction}, is a VC-subgraph class of functions with VC-index $\leq p+2$. Hence, by Lemma~9.9 Part~(iii) of \cite{kosorok2008introduction}, the class $\{\ind\{g_1'\tilde{X}^{\ast}-g_2-U^{\ast}>0\}:\; g\in A_j(k)\}$ is also VC-subgraph with VC-index $\leq p+2$. Now applying Part~(vi) of the same lemma, we obtain that ${\cal F}_j^1(k)$ is thus VC-subgraph with VC-index $\leq 2p+3$. We can also readily verify that \[F_{j1}(k)\equiv B\left[\left\{(\beta_0-\delta_0)'\tilde{Z}\right\}^2+2|(\beta_0-\delta_0)'\tilde{Z}|\cdot|\tilde{\epsilon}|\right]\leq B(b_1+b_2|\tilde{\epsilon}|),\]
for universal constants $0<b_1,b_2<\infty$, is a measurable envelope for ${\cal F}_j^1(k)$, by Condition~C4. By universal, we mean that the constants depend only on $p$, $k_2$, $\beta_0$, and $\delta_0$, and not on $j$ or $k$. We can also similarly verify that ${\cal F}_j^2(k)$ is VC-subgraph with VC-index $\leq 2p+3$ and has measurable envelope $F_{j2}(k)\equiv (1-B)(b_1+b_2|\tilde{\epsilon}|)$. 

Now Theorem~9.3 of \cite{kosorok2008introduction} yields that for $l=1,2$,
\[\sup_Q\int_0^1\sqrt{1+\log N\left(s\|F_{jl}(k)\|_{Q,2},\;{\cal F}_j^l(k),\; L_2(Q)\right)}\;ds\leq b_3,\]
where the supremum is taken over all finitely discrete probability measures $Q$ for which $\|F_{jl}(k)\|_{Q,2}>0$, and where $0<b_3<\infty$ is also universal in that it only depends on $p$ and not on $j$ or $k$. We can now apply Theorem~11.1 of \cite{kosorok2008introduction} to obtain that
\begin{eqnarray*}
\|\;\|\mathbb{G}_m\|_{{\cal F}_j^l(k)}^{\ast}\|_{P,2}&\leq& b_3b_4\|F_{jl}(k)\|_{P,2}\\
&\leq& b_3b_4\sqrt{(b_1^2+2b_1b_2\sigma+b_2^2\sigma^2)/2}\\
&\equiv&b_5,
\end{eqnarray*}
for $l=1,2$, where $b_4$ is another universal constant depending only on $p$. Thus the first assertion of the lemma follows with $c_{\ast 1}=2b_5$.

To prove the second assertion, fix $j\geq 1$ and $0<k<\infty$, and let $g\in A_j(k)$. Now we have
\begin{eqnarray*}
\lefteqn{E\left[B\ind\{g_1'\tilde{X}^{\ast}-g_2>U^{\ast}\}e^{-}(\tilde{Z},\tilde{\epsilon})+(1-B)\ind\{-g_1'\tilde{X}^{\ast}+g_2\geq U^{\ast}\}e^{+}(\tilde{Z},\tilde{\epsilon}) \right]}&&\\
&=&E\left[\frac{\left\{g_1'\tilde{X}^{\ast}-g_2\right\}_{+}\left[(\beta_0-\delta_0)'\tilde{Z}\right]^2}{2k(k_1+1)(j+1)^2}+\frac{\left\{-g_1'\tilde{X}^{\ast}+g_2\right\}_{+}\left[(\beta_0-\delta_0)'\tilde{Z}\right]^2}{2k(k_1+1)(j+1)^2}\right]\\
&=&E\left[\frac{\left|g_1'\tilde{X}^{\ast}-g_2\right|\left[(\beta_0-\delta_0)'\tilde{Z}\right]^2}{2k(k_1+1)(j+1)^2}\right]\\
&\geq&\frac{b_6(\|g_1\|+|g_2|)}{2k(k_1+1)(j+1)^2}\\
&\geq&\frac{b_6kj^2}{2k(k_1+1)(j+1)^2}\\
&\geq&\frac{b_6}{8(k_1+1)},
\end{eqnarray*}
where the first equality follows from the fact the $U^{\ast}$ is uniform over $[0,t_j(k)]$; and the first inequality, with a new universal constant $0<b_6<\infty$ not depending on $j$ or $k$, follows from Lemma~\ref{lthm:cr1} combined with the assumption from Condition~C5 that $P\{(\beta_0-\delta_0)'Z=0\}=0$. The second to last inequality follows from the fact that $\|g_1\|\vee |g_2|>kj^2$ for all $g\in A_j(k)$. The last inequality follows from the fact that $j/(j+1)\geq 1/2$ for all $j\geq 1$. Since $g$, $j$ and $k$ were all arbitrary, we have now established the second assertion with $c_{\ast 2}=b_6[8(k_1+1)]^{-1}$, and the proof is complete.
\end{proof}

\subsection{Proof of Lemma \ref{compact.convergence.concatenated}}
\label{asec:proof.compact.convergence.concatenated}

{\bf Proof.} The fact that $\hat{\Sigma}_{jn}\rightarrow \Sigma_j$ almost surely, as $n\rightarrow\infty$, for $j=1,2$, follows from the strong law of large numbers. Define $W_{1n}^{\ast}=n^{-1/2}\sum_{i=1}^n\ind\{nU_i<-4k\}\epsilon_i Z_i$ and $W^{\ast}_{2n}=n^{-1/2}\sum_{i=1}^n\ind\{nU_i>4k\}\epsilon_i Z_i$. Define also 
\begin{eqnarray*}
{\cal P}_n^{\ast}(k)&=&\left\{m_n^{-}(k),(X_i,Z_i,\epsilon_i,-nU_i):\;i\in J_n^{-}(k), m_{1n}^{-}(k), m_{2n}^{-}(k), \hat{\Sigma}_{1n}, W_{1n}^{\ast};\;\right.\\
&& \left.m_n^{+}(k),(X_i,Z_i,\epsilon_i,nU_i):\;i\in J_n^{+}(k), m_{1n}^{+}(k), m_{2n}^{+}(k), \hat{\Sigma}_{2n},W_{2n}^{\ast}\right\}.
\end{eqnarray*}
If we can show that ${\cal P}_n^{\ast}(k)\leadsto{\cal P}_0$ uniformly, as $n\rightarrow\infty$, we will be done because
$E\|W_{1n}-W_{1n}^{\ast}\|^2\leq P\{0\leq -nU_i\leq 4k\}\sigma^2k_2^2\rightarrow 0$ and $E\|W_{2n}-W_{2n}^{\ast}\|^2\leq P\{0< nU_i\leq 4k\}\sigma^2k_2^2\rightarrow 0$, and thus ${\cal P}_n^{\ast}-{\cal P}_n\rightarrow 0$ uniformly in probability.

If we condition on $m_n^{-}(k)$ and $m_n^{+}(k)$, then the following four random variables are mutually independent: $W_{1n}^{\ast}$, $\{(X_i,Z_i,\epsilon_i,-nU_i):\; i\in J_n^{-}(k), m_{1n}^{-}(k), m_{2n}^{-}(k)\}$ (as a group), $W_{2n}^{\ast}$, and $\{(X_i,Z_i,\epsilon_i,nU_i):\; i\in J_n^{+}(k), m_{1n}^{+}(k), m_{2n}^{+}(k)\}$ (as a group). We will now establish weak convergence of $(m_n^{-}(k),m_{n}^{+}(k))$, and then show that the four random variables, conditional on $(m_n^{-}(k),m_{n}^{+}(k))$, converge weakly, and this will yield the desired joint weak convergence. It is easy to verify that $(m_n^{-}(k),m_{n}^{+}(k))$ come from a trinomial distribution with sample size $n$ and probabilities $p_{1n}=P[-4k\leq nU_i\leq 0]$,
$p_{2n}=P[0<nU_i\leq 4k]$, and $p_{3n}=1-p_{1n}-p_{2n}$. Now standard arguments combined with Condition~C2 yield that $(m_n^{-}(k),m_n^{+}(k))\leadsto (m_0^{-}(k),m_0^{+}(k))$. Standard arguments also yield that $\lim\sup_{n\rightarrow\infty}P[m_n^{-}(k)>n^{1/3}]=0$ and $\lim\sup_{n\rightarrow\infty}P[m_n^{+}(k)>n^{1/3}]=0$. 

Now consider the conditional distribution of $W_{1n}^{\ast}$ given $m_n^{-}(k)=m_n$ where $0\leq m_n\leq n$ is a sequence. Based on previous arguments, we can assume without loss of generality that $0\leq m_n\leq n^{1/3}$. Letting $N=n-m_n$, we have that
\[W_{1n}^{\ast}(k)=^{\text{d}}(N+m_n)^{-1/2}\sum_{j=1}^N \tilde{\epsilon}_j\tilde{Z}_j,\]
where $=^\text{d}$ denotes equality in distribution, $\tilde{\epsilon}_j$, for $1\leq j\leq N$, are i.i.d. realizations from the distribution of $\epsilon_i$, and $\tilde{Z}_j$ are i.i.d. draws from the distribution of $Z_i$ conditional on $nU_i<-4k$. It is relatively easy to verify that the conditions of the Lindeberg-Feller Theorem apply to the (clearly mean zero) $\tilde{\epsilon}_j\tilde{Z}_j$ terms for sample size $N$. Thus $N^{-1/2}\sum_{j=1}^N \tilde{\epsilon}_j\tilde{Z}_j\leadsto W_1$. Since $N/(N+m_n)\rightarrow 1$ due to the bounds on $m_n$, we now have that, conditional on $m_n^{-}(k)$, $W_{n1}^{\ast}\leadsto W_1$, where $W_1$ is independent of $m_0^{-}(k)$. We can argue similarly that conditional on $m_n^{+}(k)$, $W_{n2}^{\ast}\leadsto W_2$, where $W_2$ is independent of $m_0^{+}(k)$. 

Now we consider weak convergence of $\{(X_i,Z_i,\epsilon_i,-nU_i):\;i\in J_n^{-}(k)\}$ conditional on the number of observations in $J_n^{-}(k)$ being equal to $m_n^{-}(k)$. The terms $(X_i,Z_i,\epsilon_i,-nU_i)$, for $i\in J_n^{-}(k)$ are i.i.d. draws from their joint distribution given $0\leq -nU_i\leq 4k$. Because $\epsilon_i$ is independent from the other terms, its distribution is unaffected by the conditioning. As $m_n^{-}(k)=O_P(1)$, it is easy to apply a minor modification of Lemma~\ref{lthm:cr0} of Supplementary Material to obtain that the characteristic function of $\{(X_i,Z_i,\epsilon_i,-nU_i):\;i\in J_n^{-}(k)\}$ given $m_n^{-}(k)=m$ converges to 
\[\prod_{j=1}^m\left[E\left(e^{i\{t_{1j}'X_1^{-}+t_{2j}'Z_1^{-}\}}\right)\;E\left(e^{it_{3j}\epsilon_1}\right)\;E\left(e^{it_{4j}\tilde{U}_1^{\ast}}\right) \right],\]
where $i=\sqrt{-1}$ in the above expression, $(t_{1j},t_{2j},t_{3j},t_{4j})\in\re^{p+d+1+1}$, $j\geq 1$, $(X_1^{-}Z_1^{+})\sim G$, $\tilde{U}_1^{\ast}$ is a uniform$[0,4k]$ random variable, and the product when $m=0$ is 1. The minor modification of the lemma involves changing the range of $U$ from $[-\nu,\nu]$ to $[-\nu,0]$, with $\nu=4k/n$, and then replacing the interval $[-1,1]$ with $[-1,0]$. Thus $\{(X_i,Z_i,\epsilon_i,-nU_i):\;i\in J_n^{-}(k)\}$ conditional on $m_n^{-}(k)=m$ converges weakly to $(X_j^{-},Z_j^{-},\epsilon_j^{-},\tilde{U}_j^{-}):\;1\leq j\leq m$. We can similarly show that $\{(X_i,Z_i,\epsilon_i,nU_i):\;i\in J_n^{+}(k)\}$ conditional on $m_n^{+}(k)=m$ converges weakly to $(X_j^{+},Z_j^{+},\epsilon_j^{+},\tilde{U}_j^{+}):\;1\leq j\leq m$. Combining these results together, we obtain that ${\cal P}_n^{\ast}(k)\leadsto {\cal P}_0(k)$, and the proof is complete.$\Box$

\subsection{Proof of sufficiency of $\mathcal P_n(k)$ and the convergence of level sets}
\label{asec:sufficiency_level_sets}

Because we also reordered the indices in $J_n^{-}(k)$ and $J_n^{+}(k)$ by $|nU_i|$, and as $\tilde{U}_j^{-}$ and $\tilde{U}_{j}^{+}$ are continuous, we have that the indices in $J_n^{-}(k)$ correspond exactly to $1\leq j\leq m_0^{-}(k)$ and those in $J_n^{+}(k)$ correspond exactly to those in $1\leq j\leq m_0^{+}(k)$ for all $n$ large enough. We will assume $n$ is large enough going forward. This allows us to reconstitute ${\cal P}_n(k)$ as follows: replace $(m_n^{-}(k),m_{1n}^{-}(k),m_{2n}^{-}(k),m_n^{+}(k),$ $m_{1n}^{+}(k),m_{2n}^{+}(k))$ with $(m_0^{-}(k),m_{10}^{-}(k),m_{20}^{-}(k),m_0^{+}(k),m_{10}^{+}(k),m_{20}^{+}(k))$; for each $1\leq j\leq m_0^{-}(k)$, set $(X_{j-},Z_{j-},\epsilon_{j-},U_{j-})=(X_i,Z_i,\epsilon_i, |nU_i|)$ for the matching observation $i\in J_n^{-}(k)$; and for each $1\leq j\leq m_0^{+}(k)$, set $(X_{j+},Z_{j+},\epsilon_{j+},U_{j+})=(X_i,Z_i,\epsilon_i, nU_i)$ for the matching observation $i\in J_n^{+}(k)$. Thus
\begin{eqnarray*}
{\cal P}_n(k)&=&\left\{m_0^{-}(k),(X_{j-},Z_{j-},\epsilon_{j-},U_{j-}):\;1\leq j\leq m_0^{-}(k), m_{10}^{-}(k), m_{20}^{-}(k), \hat{\Sigma}_{1n}, W_{1n};\;\right.\\
&& \left.m_0^{+}(k),(X_{j+},Z_{j+},\epsilon_{j+},U_{j+}):\;1\leq j\leq m_0^{+}(k), m_{10}^{+}(k), m_{20}^{+}(k), \hat{\Sigma}_{2n}, W_{2n}\right\},
\end{eqnarray*}
where 
\[\lim\sup_{n\rightarrow\infty}\;\max_{1\leq j\leq m_0^{-}(k)}\left[\|X_{j-}-X_j^{-}\|\vee\|Z_{j-}-Z_j^{-}\|\vee|\epsilon_{j-}-\epsilon_j^{-}|\vee |U_{j-}-\tilde{U}_j^{-}|\right]=0\]
and
\[\lim\sup_{n\rightarrow\infty}\;\max_{1\leq j\leq m_0^{+}(k)}\left[\|X_{j+}-X_j^{+}\|\vee\|Z_{j+}-Z_j^{+}\|\vee|\epsilon_{j+}-\epsilon_j^{+}|\vee |U_{j+}-\tilde{U}_j^{+}|\right]=0,\]
where we take the maximum over a null set to be zero in this situation. Also, for each $1\leq j\leq m_0^{-}(k)$, set $E_{jn-}(h_1,h_2)=E_{in}^{-}(h_1,h_2)$ for the matching observation $i\in J_n^{-}(k)$; and for each $1\leq j\leq m_0^{+}(k)$, set $E_{jn+}(h_1,h_2)=E_{in}^{+}(h_1,h_2)$ for the matching observation $i\in J_n^{+}(k)$. We can also verify that 
\[\lim\sup_{n\rightarrow\infty}\max_{1\leq j\leq m_0^{-}(k)}\;\sup_{(h_1,h_2)\in\re^d\times\re^d:\|h_1\|\vee\|h_2\|\leq k}\left|E_{jn-}(h_1,h_2)-\tilde{E}_j^{-}\right|=0,\]
and 
\[\lim\sup_{n\rightarrow\infty}\max_{1\leq j\leq m_0^{+}(k)}\;\sup_{(h_1,h_2)\in\re^d\times\re^d:\|h_1\|\vee\|h_2\|\leq k}\left|E_{jn+}(h_1,h_2)-\tilde{E}_j^{+}\right|=0.\]

Note that with $F_n(k)=1$,
\[\arg\min_{h\in H_n(k)}\left[Q_n(h)=Q_{1n}(h_1,h_2)+Q_{2n}(h)\right]=\arg\min_{h\in H_n(k)}\left[Q_{1n}(h_1,h_2)+Q_{2n}^{\ast}(h)\right],\]
where 
\begin{eqnarray*}
Q_{2n}^{\ast}(h)&=&Q_{2n}^{\ast -}(h)+Q_{2n}^{\ast +}(h)\\
&\equiv&\sum_{i\in J_n^{-}(k)}\ind\{-h_3'X_i+h_4<nU_i\leq 0\}E_{in}^{-}(h_1,h_2) +\sum_{i\in J_n^{+}(k)}\ind\{0<nU_i\leq -h_3'X_i+h_4\}E_{in}^{+}(h_1,h_2),
\end{eqnarray*}
because the components in the sum defining $Q_{2n}(h)$ corresponding to indices $\{1\leq i\leq n:\;i\not\in J_n^{-}(k)\cup J_n^{+}(k)\}$ are constant over $h\in H_n(k)$. Hence the $\arg\min$ of the process $Q_{n}(h)$ only depends on the data in ${\cal P}_n(k)$, and thus $\tilde{h}_n(k)$ in (\ref{eq:tilde.hn.k}) is also a function of ${\cal P}_n(k)$. Since $(\hat{h}_{1n}(k),\hat{h}_{2n}(k))=(\tilde{h}_{1n}(k),\tilde{h}_{2n}(k))$, we have that 
$(\tilde{h}_3(k),\tilde{h}_4(k))\in\arg\min_{(h_3,h_4)\in H_{2n}(k)}\tilde{Q}_{2n}^{\ast}(h_3,h_4)$, where
\begin{eqnarray*}
\tilde{Q}_{2n}^{\ast}(h_3,h_4)&=&\sum_{i\in J_n^{-}(k)}\ind\{-h_3'X_i+h_4<nU_i\leq 0\}E_{in}^{-}(\hat{h}_{1n}(k),\hat{h}_{2n}(k))\\
&&+\sum_{i\in J_n^{+}(k)}\ind\{0<nU_i\leq -h_3'X_i+h_4\}E_{in}^{+}(\hat{h}_{1n}(k),\hat{h}_{2n}(k))\\
&=&\sum_{j=1}^{m_0^{-}(k)}\ind\{-h_3'X_{j-}+h_4<-U_{j-}\leq 0\}E_{jn-}(\hat{h}_{1n}(k),\hat{h}_{2n}(k))\\
&&+\sum_{j=1}^{m_0^{+}(k)}\ind\{0<U_{j+}\leq -h_3'X_{j+}+h_4\}E_{jn+}(\hat{h}_{1n}(k),\hat{h}_{2n}(k))\\
&\equiv&\tilde{Q}_{2n}^{\ast -}(h_3,h_4)+\tilde{Q}_{2n}^{\ast +}(h_3,h_4).
\end{eqnarray*}
Thus $(\tilde{h}_{3n}(k),\tilde{h}_{4n}(k))\in\arg\min_{(h_3,h_4)\in H_{2n}(k)}\tilde{Q}_{2n}^{\ast}(h_3,h_4)$ depends only on the data in ${\cal P}_n(k)$. Similarly, the level set $\Phi_n'(k)$ depends only on the data in ${\cal P}_n(k)$; and we can also verify that $V_i(k)$ as defined above, for $i\in J_n^{-}(k)\cup J_n^{+}(k)$, also depends on the data in ${\cal P}_n(k)$. Linking up indices as before, we can define $V_{j-}(k)=V_i(k)$ for $1\leq j\leq m_0^{-}(k)$ via linking with the appropriate $i\in J_n^{-}(k)$, and also $V_{j+}(k)=V_i(k)$ for $1\leq j\leq m_0^{+}(k)$ via the appropriate linking with $i\in J_n^{+}(k)$. Thus $\Phi_n'(k)=\Big\{(h_3,h_4)\in H_{2n}(k):\ind\{h_3'X_{j-}-h_4-U_{j-}>0\}-\ind\{h_3'X_{j-}-h_4-U_{j-}\leq 0\}=V_{j-}(k)$, $1\leq j\leq m_0^{-}(k),\;\mbox{and}\;\ind\{h_3'X_{j+}-h_4+U_{j+}>0\}-\ind\{h_3'X_{j+}-h_4+U_{j+}\leq 0\}=V_{j+}(k)$, $\;1\leq j\leq m_0^{+}(k)\Big\}$. Now we can re-express several previously defined quantities in terms of data in ${\cal P}_n(k)$, specifically, 
\begin{eqnarray*}
\lefteqn{h_3\mapsto C_L^{nk}(h_3)\;=}&&\\
&&\left(\max_{1\leq j\leq m_0^{-}(k):\;V_{j-}(k)=-1} h_3'X_{j-}-U_{j-}\right)\vee\left(\max_{1\leq j\leq m_0^{+}(k):\;V_{j+}(k)=-1} h_3'X_{j+}+U_{j+}\right),\\
\lefteqn{h_3\mapsto C_U^{nk}(h_3)\;=}&&\\
&&\left(\min_{1\leq j\leq m_0^{-}(k):\;V_{j-}(k)=1} h_3'X_{j-}-U_{j-}\right)\wedge\left(\min_{1\leq j\leq m_0^{+}(k):\;V_{j+}(k)=1} h_3'X_{j+}+U_{j+}\right),
\end{eqnarray*}
and $C_R^{nk}$, $R_n(k)$ and $\hat{h}_{3n}(k)$ are defined as before relative to $C_L^{nk}$ and $C_U^{nk}$ and thus depend only on data in ${\cal P}_n(k)$.
By Lemma~\ref{lemma.onetoone} and surrounding arguments, we have that for all $n$ large enough, $(\overline{\omega}_0'h_3,h_4)\in H_{20}^{\ast}(k)$ for all $(h_3,h_4)\in H_{2n}(k)$ and $(h_{\ast n}(g_1),g_2)\in H_{2n}(k)$ for all $(g_1,g_2)\in H_{20}^{\ast}(k)$. Thus the process $(h_3,h_4)\mapsto \tilde{Q}_{2n}^{\ast}(h_3,h_4)$ ranging over $H_{2n}(k)$ is equivalent to the process $(g_1,g_2)\mapsto \tilde{Q}_{2n}(g_1,g_2)\equiv \tilde{Q}_{2n}^{\ast}(h_{\ast n}(g_1),g_2)$ ranging over $(g_1,g_2)\in H_{20}^{\ast}(k)$, for all $n$ large enough. We can similarly define $\tilde{C}_L^{nk}(g_1)=C_L^{nk}(h_{\ast n}(g_1))$, $\tilde{C}_U^{nk}(g_1)=C_U^{nk}(h_{\ast n}(g_1))$, $\tilde{C}_R^{nk}=\tilde{C}_U^{nk}-\tilde{C}_L^{nk}$, and $\tilde{R}_n(k)=\{g_1\in\re^{p-1}:\;\tilde{C}_R^{nk}(g_1)>0,\;\|g_1\|\leq 0.9k_1^{-1}k\}$. We can now verify that
\[\hat{h}_{3n}(k)=\frac{\int_{\tilde{R}_n(k)}h_{\ast n}(g_1)\tilde{C}_R^{nk}(g_1)\left(\frac{d\nu_n(h_{\ast n}(g_1))}{d\mu(g_1)}\right)d\mu(g_1)}{\int_{\tilde{R}_n(k)}\tilde{C}_R^{nk}(g_1)\left(\frac{d\nu_n(h_{\ast n}(g_1))}{d\mu(g_1)}\right)d\mu(g_1)},\]
where, similar to what was done previously, $\mu$ is Lebesgue measure on $\re^{p-1}$ and $\nu_n$ is Lebesgue measure on the sphere $nS^{p-1}$ oriented so that a normal vector $e_p$ in the $p$'th dimension orthogonal to the hyperplane $\re^{p-1}$ passes through the center of the sphere and the sphere touches $\re^{p-1}$ at the point $e_p=0$. We can make this transformation because the orthonormal rotational transformation matrix from $\re^p$ to $\re_{\omega_0}\times\re_{\overline{\omega}_0}^{p-1}$, where $\re_{\omega_0}$ is just the linear span of $\omega_0$, which is the matrix $L_0=(\omega_0,\overline{\omega}_0)$, has determinant 1. The quantity $d\nu_n/d\mu$ is the Radon-Nikodym derivative of $\nu_n$ relative to $\mu$ defined via the transformation $g_1\mapsto h_{\ast n}(g_1)$. Because the curvature of the surface of $S^{p-1}$ is always greater than zero, we know that $d\nu_n/d\mu\geq 1$. We can also show that, for all $n$ large enough and $\|g_1\|\leq 0.9 k_1^{-1}k$, the infinitesimal curvature $h_{\ast n}(dg_1)/dg_1$ is bounded above by
\[\prod_{l=1}^{p-1}\sqrt{1+\frac{g_{1,l}^2}{n^2-\|g_1\|^2}}\;\rightarrow\; 1,\]
as $n\rightarrow\infty$, where $g_1=(g_{1,1},\ldots,g_{1,p-1})'$. Moreover, we can also readily verify that
\[\lim\sup_{n\rightarrow\infty}\;\max_{1\leq j\leq m_0^{-}(k)}\;\sup_{(g_1,g_2)\in H_{20}^{\ast}(k)}\;\left|h_{\ast n}'(g_1)X_{j-}-g_2-\left(g_1'\tilde{X}_j^{-}-g_2\right)\right|=0,\]
and 
\[\lim\sup_{n\rightarrow\infty}\;\max_{1\leq j\leq m_0^{+}(k)}\;\sup_{(g_1,g_2)\in H_{20}^{\ast}(k)}\;\left|h_{\ast n}'(g_1)X_{j+}-g_2-\left(g_1'\tilde{X}_j^{+}-g_2\right)\right|=0.\]

Recall $\tilde{V}_j^{-}(k)$ and $\tilde{V}_j^{+}(k)$, defined above for all $1\leq j<\infty$, and note that when $F_0(k)=1$, we have that $\tilde{V}_j^{-}(k)=-1$ for all $m_{10}^{-}(k)<j<\infty$ and $\tilde{V}_j^{+}(k)=1$ for all $m_{10}^{+}(k)<j<\infty$. Thus 
$\tilde{g}(k) \in \arg\min_{g\in H_{20}^\ast(k)}Q_{02}(g)=\arg\min_{g\in H_{20}^\ast(k)}Q_{02}^{\ast}(g)$, where
\begin{eqnarray*}
Q_{02}^{\ast}(g)&=&\sum_{j=1}^{m_{0}^{-}(k)}\ind\{-g_1'\tilde{X}_j^{-}+g_2<-\tilde{U}_j^{-}\leq 0\}\tilde{E}_j^{-}+\sum_{j=1}^{m_{0}^{+}(k)}\ind\{0<\tilde{U}_j^{+}\leq -g_1'\tilde{X}_j^{+}+g_2\leq 0\}\tilde{E}_j^{+}.
\end{eqnarray*}
Based on arguments up to this point, we can also verify that for all $n$ large enough $V_{j-}=\tilde{V}_j^{-}(k)=-1$ for all $m_{10}^{-}(k)<j\leq m_0^{-}(k)$ and $V_{j+}=\tilde{V}_j^{+}(k)=1$ for all $m_{10}^{+}<j\leq m_0^{+}(k)$. Thus the non-zero entries in both $\tilde{Q}_{2n}(g)$ and $Q_{02}^{\ast}(g)$, for all $g\in H_{20}^\ast(k)$, can only occur when $\ind\{-h_{\ast n}(g_1)'X_{j-}+g_2<-U_{j-}\leq 0\}=1$ or $\ind\{-g_1'\tilde{X}_j^{-}+g_2<-\tilde{U}_j^{-}\leq 0\}=1$, for $1\leq j\leq m_{10}^{-}(k)$; or when $\ind\{0<U_{j+}\leq -h_{\ast n}(g_1)'X_{j+}+g_2\}=1$ or $\ind\{0<\tilde{U}_j^{+}\leq -g_1'\tilde{X}_j^{+}+g_2\}=1$, for $1\leq j\leq m_{10}^{+}(k)$. Our next step will be verify that  $V_{j-}=\tilde{V}_j^{-}(k)$ for all $0\leq j\leq m_{10}^{-}(k)$ and $V_{j+}=\tilde{V}_j^{+}(k)$ for all $0\leq j\leq m_{10}^{+}$. This will allow us to link the level sets between the processes $Q_{02}^{\ast}$ and $\tilde{Q}_{2n}$.

Toward this end, the maximum number of possible level sets in $Q_{02}^{\ast}$ is at most $2$ to the power of ${m_{10}^{-}(k)+m_{10}^{+}(k)}$, as the number of unique indicator functions is $m_{10}^{-}(k)+m_{10}^{+}(k)$ and each one can be either 0 or 1 (at most). However, not all possibilities may be feasible for $g$ ranging over $H_{20}^\ast(k)$, so the actual number will generally be lower. For each of the realized level sets, the set of $g\in H_{20}^\ast(k)$ must include an open ball in $\re^{p-1}$ almost surely because the $\tilde{U}_j^{-}$ and $\tilde{U}_j^{+}$ are continuous random variables. Also, there will almost surely be no ties between level sets defined in terms of the indicator functions because all possible values of $\tilde{E}_j^{-}$ and $\tilde{E}_j^{+}$ and their finite sums will be unique due to Condition~C5 and the assumed continuity of the residual $\epsilon$. Because of the previously established uniform convergence of the components in ${\cal P}_n(k)$ to those in ${\cal P}_0$ combined with the compactness of $H_{20}^\ast(k)$, we obtain that all of the boundaries between the level sets in $Q_{02}^{\ast}$ and $\tilde{Q}_{2n}$ will converge uniformly and the magnitudes of the respective sums will converge. Recall that $\hat{g}(k)$ is in the interior of $\tilde{\Phi}_0(k)\cap H_{20}^\ast(k)$, the level set of the $\arg\max$ of $Q_{02}^{\ast}$. Thus, for all $n$ large enough $\hat{g}(k)$ will be contained in the interior of $\tilde{\Phi}_n(k)=\{g\in H_{2n}(k):\; h_{\ast n}(g)\in \Phi_n'(k)\}$, the level set of the $\arg\max$ of $\tilde{Q}_{2n}$. Hence $V_{j-}=\tilde{V}_j^{-}(k)$ for all $1\leq j\leq m_0^{-}(k)$ and $V_{j+}=\tilde{V}_j^{+}(k)$ for all $1\leq j\leq m_0^{+}(k)$ almost surely, for all $n$ large enough.

Let $\tilde{C}_L^k(g_1)$, $\tilde{C}_U^k(g_1)$, $\tilde{C}_R^k(g_1)$, $\tilde{R}_0(k)$, and $\hat{g}_1(k)$, for $(g_1,g_2)\in H_{20}^{\ast}(k)$, be as defined previously but applied to the fixed data in ${\cal P}_0(k)$; The forgoing convergence results taken together now yield that for each $k<\infty$, 
\begin{eqnarray*}
\lim_{n\rightarrow\infty}\sup_{g_1\in\re^{p-1}:\;\|g_1\|\leq 0.9k_1^{-1}k}\left|\tilde{C}_L^{nk}(g_1)-\tilde{C}_L^k(g_1)\right|\vee \left|\tilde{C}_U^{nk}(g_1)-\tilde{C}_U^k(g_1)\right|\vee\left|\tilde{C}_R^{nk}(g_1)-\tilde{C}_R^k(g_1)\right|&=&0.
\end{eqnarray*}
\FloatBarrier

\section{Proofs of Section \ref{sec:inference}} 
\label{asec:inference}
\begin{proof}[Proof of Theorem \ref{parametric.bootstrap}]
From previous results, it is easy to verify that $\max_{1\leq i\leq n}|\hat{U}_i-U_i|\leq e_n\equiv k_1\|\hat{\omega}_n-\omega_0\|+|\hat{\gamma}_n-\gamma_0| =O_P(n^{-1})$. Since $|U|\leq k_1+\gamma_0$ almost surely, we now have that $\hat{\tau}_n^2\rightarrow\mbox{var}(U)\equiv\tau_0^2$ in probability. Note also that for all $t\in\re$,
\[\ep_n\{U\leq t-e_n\}\leq \hat{F}_{n1}(t)\leq \ep_n\{U\leq t+e_n\},\]
and thus, by using standard empirical process results combined with Condition C2', we obtain that for
$F_0(t)\equiv\int_{-\infty}^tf(s)ds$,
\[\|\hat{F}_{n1}-F_0\|_{\re}\;\leq\; O_P(n^{-1/2})+\sup_{t\in\re}P\{t-e_n< U\leq t+e_n\}\;\leq\;O_P(n^{-1/2}).\]
By lemma~\ref{density.estimation} below, the proof of which will be given later in this section, we now have that $\hat{f}_{n0}\rightarrow f_0$ in probability, since $\tilde{\eta}_n=\hat{\eta}_{n1}$ readily satisfies the conditions of the lemma and since C2' ensures that $f$ is continuous at $u=0$ as needed for Part 1 of the lemma. We note that the $n^{1/5}$ power was selected because of its asymptotic optimality when the second derivative of the density being estimated is uniformly equicontinuous (see, e.g., \cite{Wahba1975}), but this stronger condition is not required since consistency is all we need and will still be achieved under the weaker assumptions we are using.

\begin{lemma}\label{density.estimation}
Let $\tilde{f}_0$ be a density on $\re$ that is uniformly bounded, and let $\tilde{F}_0(t)\equiv\int_{-\infty}^t\tilde{f}_0(s)ds$. Let $\tilde{F}_n$ be an estimator of $\tilde{F}_0$ such that $\|\tilde{F}_n-\tilde{F}_0\|_{\re}=O_P(n^{-1/2})$. Also, for $u\in\re$,
define
\[u\mapsto\tilde{f}_n(u)=\int_{\re}\frac{1}{\tilde{\eta}_n}\phi\left(\frac{u-t}{\tilde{\eta}_n}\right)d\tilde{F}_n(t),\]
where $\tilde{\eta}_n$ is a (possibly random) sequence satisfying $\tilde{\eta}_n=o_P(1)$ and $\tilde{\eta}_n^{-1}=o_P(n^{1/2})$, and where $\phi$ is the standard normal density. Then we have the following:
\begin{enumerate}
    \item Provided $\tilde{f}_0$ is continuous at $u$, then $\tilde{f}_n(u)\rightarrow\tilde{f}_0(u)$ in probability.
    \item Provided $\lim_{\eta\downarrow 0}\;\sup_{s,t\in\re:\; |s-t|\leq\eta} |\tilde{f}_0(s)-\tilde{f}_0(t)|=0$, then $\|\tilde{f}_n-\tilde{f}_0\|_{\re}=o_P(1)$.
\end{enumerate}
\end{lemma}

Continuing with the proof of Theorem~\ref{parametric.bootstrap}, we now verify consistency of $\tilde{M}_{n}$. As defined, $\tilde{M}_n$ satisfies $\tilde{M}_n\tilde{M}_n'=I-\hat{\omega}_n\hat{\omega}_n'\equiv \tilde{A}_n$ but is not in general unique since, for example, rearranging the columns of $\tilde{M}_n$ does not change $\tilde{A}_n$. This non-uniqueness will be addressed shortly. However, we first address consistency of the projection $\tilde{A}_n$ as an estimator of $\tilde{A}_0\equiv I-\omega_0\omega_0'$. Previous arguments yield that for any $t\in\re^p$, $t'(\tilde{A}_n-\tilde{A}_0)t=O_P(n^{-1})$. This implies that $\tilde{A}_n$ converges to $\tilde{A}_0$ in probability as a projection. Let $\hat{\re}_p$ be the subset of $\re^p$ equal to the range of the projection $\tilde{A}_n$. Note that $g_1'\tilde{X}_{\ast}(t)=g_1'\tilde{M}_n'X$, for $(X,Z)$ drawn from $\tilde{G}_n$. Thus minimizing over $g=(g_1,g_2)\mapsto \tilde{Q}_{02}^{\ast}(g)$, and then multiplying the resulting minimizer by $\tilde{M}_n$, is the same as minimizing $\tilde{Q}_{02}^{\ast}$ over $(h_3,g_2)$, after replacing $g_1'\tilde{X}_{\ast}(t)$ with $h_3'X$, for $h_3\in\hat{\re}_p$. The point of this is that the non-uniqueness of $\tilde{M}_n$ is not a problem provided $\tilde{M}_n\tilde{M}_n'=I-\hat{\omega}_n\hat{\omega}_n'$. 

Note that, by Condition C4,
\[|\hat{\epsilon}_i-\epsilon_i|\leq k_2(\|\hat{\beta}_n\|+\|\hat{\delta}_n\|)\hat{E}_i+k_2(\|\hat{\beta}_n-\beta_0\|+\|\hat{\delta}_n-\delta_0\|)\equiv\tilde{e}_n,\]
where $\hat{E}_i\equiv\ind\{\hat{U}_i\leq 0<U_i\}+\ind\{U_i\leq 0 < \hat{U}_i\}$, $i=1,\ldots,n$. Previous arguments verify that $n^{-1}\sum_{i=1}^n\hat{E}_i=O_P(n^{-1})$, and thus, after some derivation, $\hat{\sigma}_n^2=\sigma^2+O_P(n^{-1/2})$, $\tilde{\epsilon}_n=O_P(n^{-1/2})$, and also $\tilde{\Sigma}_{jn}=\Sigma_j+O_P(n^{-1/2})$ for $j=1,2$. Now, for all $t\in\re$,
\[\ep_n\{\epsilon\leq t+\tilde{\epsilon}_n-\tilde{e}_n\}\leq \hat{F}_{n2}(t)\leq \ep_n\{\epsilon\leq t+\tilde{\epsilon}_n+\tilde{e}_n\},\]
and thus, by using standard empirical process results combined with Condition C6, we obtain that for $F_0^{\ast}(t)=\int_{-\infty}^t\xi(s)ds$,
\begin{eqnarray*}
\|\hat{F}_{n2}-F_0^{\ast}\|_{\re}&\leq& O_P(n^{-1/2})+\sup_{t\in\re}P\{t-\tilde{e}_n<\epsilon\leq t+\tilde{e}_n\}\\
&&+\sup_{t\in\re}|P\{\epsilon\leq t+\tilde{\epsilon}_n\}-F_0^{\ast}(t)|\\
&=&O_P(n^{-1/2}),
\end{eqnarray*}
since, regarding the first supremum above, the supremum of a real function over $t\in\re$ is the same as the supremum over $t+\tilde{\epsilon}_n\in\re$. Now Condition C6 combined with Part 2 of Lemma~\ref{density.estimation} yield that $\|\hat{\xi}_n-\xi\|_{\re}=o_P(1)$.

Define $t_n=\sup\left\{t>0:\;\sum_{i=1}^n\ind\{|U_i|\leq t\}\leq r_n\right\}$, and note that, by definition of $e_n$, $t_n-e_n\leq \hat{t}_n\leq t_n+e_n$. Since it can also be shown that $t_n=r_n/(2f_0n)+O_P(n^{-1/2})$, we have that $e_n/t_n=o_P(1)$ as a consequence of the fact that $r_n n^{-1/2}\rightarrow\infty$. If we let $t_{n0}=r_n/(2f_0n)$, we now have that both $\hat{t}_n/t_n=1+o_P(1)$ and $\hat{t}_n/t_{n0}=1+o_P(1)$. Hence, for $\hat{J}_n\equiv \left\{i\in\{1,\ldots,n\}:\; |\hat{U}_i|\leq \hat{t}_n\right\}$, we have for any $\eta>0$ that $P\left\{\max_{i\in\hat{J}_n}\,|U_i|\leq t_n(1+\eta)\right\}\rightarrow 1$, as $n\rightarrow\infty$. Let $0<m<\infty$ be an integer, and let $(X_j^{\ast},Z_j^{\ast})$, $j=1,\ldots,m$ be an i.i.d. sample taken with replacement from $\tilde{G}_n$. Based on the above discussion, this is identical to drawing with replacement $m$ integers from $\hat{J}_n$, which integers we denote as $i(1),\ldots,i(m)$, and then letting $(X_j^{\ast},Z_j^{\ast})=(X_{i(j)},Z_{i(j)})$, for $j=1,\ldots,m$. Since $r_n\rightarrow\infty$, as $n\rightarrow\infty$, we have that the probability that $i(j)=i(k)$, for any $1\leq j<k\leq m$, is $m(m-1)/(2r_n)\rightarrow 0$, as $n\rightarrow\infty$. Thus, given $1\leq m<\infty$, we have, with probability going to one, that there are no ties in the indices drawn from a sample of size $m$, as $n\rightarrow\infty$. We now have, by combining previous results, along with reapplication of Lemma~\ref{lthm:cr0}, that for any $(s_1,u_1),\ldots,(s_m,u_m)\in\re^p\times\re^d$,
\begin{eqnarray*}
E\left[\left.\exp\left\{i\sum_{j=1}^m s_j'X_j^{\ast}+u_j'Z_j^{\ast}\right\}\right|\,{\cal X}_n\right]&=&E\left[\left.\left(\prod_{j=1}^m\exp\left\{i(s_j'X_j^{\ast}+u_j'Z_j^{\ast})\right\}\right)\right|\,{\cal X}_n\right]\\
&\rightarrow&\prod_{j=1}^m\chi(s_j,u_j),
\end{eqnarray*}
in probability, as $n\rightarrow\infty$, where $i$ in the above expression is $\sqrt{-1}$, ${\cal X}_n=\left\{(X_i,Z_i,Y_i),\,1\leq i\leq n\right\}$ is the observed data from a sample of size $n$, and $\chi$ is the characteristic function of $G$.

We now can rework the theory from the previous section to verify that the restrictions of the estimators $\hat{g}_{\ast}$ and $\check{g}_{\ast}$ as $\arg\min$s of $\tilde{Q}_{02}^{\ast}(g)$ over the compact set $H_{20}^{\ast}(k)$, for $0<k<\infty$, converge in distribution, conditional on ${\cal X}_n$, to $\hat{g}(k)$ and $\check{g}(k)$, respectively, in probability as $n\rightarrow\infty$. For the convergence of $\check{g}_{\ast}$ to $\check{g}$, the asymptotic impossibility of ties in the indices established in the previous paragraph ensures that, under Conditions C3', there are asymptotically no ties in the values of $X_1^{\ast}, X_2^{\ast},\ldots$ observed in $\tilde{Q}_{02}^{\ast}(g)$ for $g$ restricted to $H_{20}^{\ast}(k)$, for any fixed $0<k<\infty$. We can also establish the needed compactness of $\tilde{g}_{\ast}$, $\hat{g}_{\ast}$ and $\check{g}_{\ast}$, using arguments parallel to those used in the proof of Theorem~\ref{argmax.compactness}. By Theorem~\ref{convergence-modemid}, we can, under C3', replace $\hat{\phi}_n$ with $\check{\phi}_n$ in defining $\hat{U}_i$, without changing the conclusions. The desired conditional weak convergence for all components now follows via arguments similar to those used in the proofs of Theorems~\ref{convergence-meanmid} and~\ref{convergence-modemid}.
\end{proof}
\begin{proof} [Proof of Lemma~\ref{density.estimation}]
Using Riemann-Stieltjes integration by parts, we have for every $u\in\re$,
\[\int_{\re}\frac{1}{\tilde{\eta}_n}\phi\left(\frac{u-t}{\tilde{\eta}_n}\right)d\tilde{F}_n(t)-\int_{\re}\frac{1}{\tilde{\eta}_n}\phi\left(\frac{u-t}{\tilde{\eta}_n}\right)d\tilde{F}_0(t)
=-\int_{\re}\left[\tilde{F}_n(t)-\tilde{F}_0(t)\right]\frac{(u-t)}{\tilde{\eta}_n^3}\phi\left(\frac{u-t}{\tilde{\eta}_n}\right)dt\]
$\equiv E_1$, where 
\[|E_1|\leq \frac{\|\tilde{F}_n-\tilde{F}_0\|_{\re}}{\tilde{\eta}_n}\equiv E_2,\]
and where $E_2=O_P(\tilde{\eta}_n^{-1} n^{-1/2})=o_P(1)$ does not depend on $u$. Moreover, for every $\Delta>0$,
\begin{eqnarray*}
\int_{\re}\frac{1}{\tilde{\eta}_n}\phi\left(\frac{u-t}{\tilde{\eta}_n}\right)d\tilde{F}_0(t)&=&\int_{u-\Delta}^{u+\Delta}\frac{1}{\tilde{\eta}_n}\phi\left(\frac{u-t}{\tilde{\eta}_n}\right)\tilde{f}_0(u)dt\\
&&-\int_{u-\Delta}^{u+\Delta}\frac{1}{\tilde{\eta}}\phi\left(\frac{u-t}{\tilde{\eta}_n}\right)\left[\tilde{f}_0(u)-\tilde{f}_0(t)\right]dt\\
&&+\int_{(-\infty,u-\Delta]\cup(u+\Delta,\infty)}\frac{1}{\tilde{\eta}}\phi\left(\frac{u-t}{\tilde{\eta}_n}\right)\tilde{f}_0(t)dt\\
&=&\tilde{f}_0(u)E_2(u)+E_3(u)+E_4(u),
\end{eqnarray*}
where $E_2(u)=\left[\Phi(\tilde{\eta}_n^{-1}\Delta)-\Phi(-\tilde{\eta}_n^{-1}\Delta)\right]\equiv E_n^{\ast}\rightarrow 1$, in probability, as $n\rightarrow\infty$. Note that the dependency of $E_2(u)$ on $u$ vanishes via a change of variables in the integral. Next, it is easy to verify that
\[|E_4(u)|\leq [1-E_n^{\ast}]\times \|\tilde{f}_0\|_{\re}\rightarrow 0 \]
in probability, as $n\rightarrow\infty$, since $\|\tilde{f}_0\|_{\re}<\infty$ by assumption. Finally, we have that
\[|E_3(u)|\leq E_n^{\ast}\times\sup_{u-\Delta\leq t\leq u+\Delta} |\tilde{f}_0(t)-\tilde{f}_0(u)|\leq \sup_{u-\Delta\leq t\leq u+\Delta} |\tilde{f}_0(t)-\tilde{f}_0(u)|\equiv E_5(u,\Delta),\]
where, since $\Delta$ is arbitrary, we can allow it to get arbitrarily small. For Part 1 of the proof, we have for a fixed $u$, that continuity of $\tilde{f}_0$ at $u$ yields $\lim_{\Delta\downarrow 0}E_5(u,\Delta)=0$, and thus $\tilde{f}_n(u)\rightarrow \tilde{f}_0(u)$ in probability, as $n\rightarrow\infty$. Thus Part 1 follows. For Part 2, the uniform equicontinuity of $\tilde{f}_0$ yields that 
\[\lim_{\Delta\downarrow 0}\sup_{u\in\re}E_5(u,\Delta)=0,\]
and hence $\|\tilde{f}_n-\tilde{f}_0\|_{\re}\rightarrow 0$ in probability, as $n\rightarrow\infty$. Thus Part 2 also follows.
\end{proof}
\FloatBarrier

\section{Numerical estimation}
\subsection{Uniform search}
\label{asec:search}
In this section, we provide the algorithmic details of the numerical estimation procedures. We first state the algorithm, Algorithm \ref{al:uniform_search} below,  of the uniform search approach described in Section \ref{sec:simest_uniform}. This algorithm calls Algorithms 2--5 below and outputs $(\hat\omega, \hat\gamma)$ and $(\check\omega, \check\gamma)$.
\begin{algorithm}[h!]
\caption{The main algorithm. $SSR(\omega, \gamma)$ is defined as $\text{SSR}_n(\omega, \gamma)$ in Section~\ref{sec:simest_uniform}.}
\label{al:uniform_search}
Constants: $DECAY=0.8; MAXHITS=20; MAXITER=500, J=10^4$\\
\begin{algorithmic}[1]
\State $k\leftarrow 0$,  $h\leftarrow 0$, $C\leftarrow 1$; 
\State $\omega^* \leftarrow \omega^{(0)}\leftarrow~\texttt{Relax}()$;  \hfill (See Algorithm \ref{al:relax})
\State $\text{SSR}^{*} \leftarrow \text{SSR}^{(0)} \leftarrow \text{SSR}(\omega^{(0)})$;
\For{$k = 1, 2, ..., MAXITER$}
    \If {$h < MAXHITS$} 
        \State EXIT;
    \EndIf
    \If{$SSR^{(k)} > SSR^{*}$}
        \State $SSR^* \leftarrow SSR^{(k-1)}$; 
        \State $\omega^* \leftarrow \omega^{(k-1)}$, $\gamma^* \leftarrow \gamma^{(k-1)}$;
        \State $C \leftarrow C \cdot DECAY$;
        \State $h  \leftarrow 0$;
    \Else
        \State $h \leftarrow h + 1$;
    \EndIf
    
    \For{$j = 1,2, ...,J$}
        \State $\omega^{(k,j)} \leftarrow \texttt{UnifSphere}(\omega^*, C)$; \hfill (See Algorithm \ref{al:uniform})
        \State $\gamma^{(k,j)} \leftarrow \arg\max_{\gamma}SSR(\omega^{(k,j)}, \gamma)$; 
        \State $SSR(\omega^{(k,j)}) \leftarrow SSR(\omega^{(k,j)}, \gamma^{(k,j)})$;
    \EndFor
    \State $j^* \leftarrow {\arg\max}_{j=1,2,...,J} SSR(\omega^{(k,j)})$;
    \State $\omega^{(k)} \leftarrow \omega^{(k, j^*)}$; $\gamma^{(k)}  \leftarrow \gamma^{(k, j^*)}$;
    \State $k \leftarrow k + 1$
\EndFor  
\State $\hat \omega, \hat \gamma \leftarrow \texttt{MeanArg}(\omega^*, \gamma^*)$; \hfill (See Algorithm \ref{al:arg})
\State $\check \omega, \check \gamma \leftarrow \texttt{ModeArg}(\omega^*, \gamma^*)$;  \hfill (See Algorithm \ref{al:arg})
\State $\hat \beta  \leftarrow \check \beta \leftarrow  Z_+(\omega^*, \gamma^*)'\left\{ Z_+(\omega^*, \gamma^*)  Z_+(\omega^*, \gamma^*)'\right\}^-Y$;
\State $\hat \delta  \leftarrow \check \delta \leftarrow  Z_-(\omega^*, \gamma^*)'\left\{ Z_-(\omega^*, \gamma^*)  Z_-(\omega^*, \gamma^*)'\right\}^-Y$;
\end{algorithmic}
\end{algorithm}

\begin{algorithm}[h!]
\caption{$\texttt{Relax}()$ finds an initial guess for $\tilde \omega_n$ using a smooth surrogate loss function in the objective function.}
\label{al:relax}
Constants: $MAXITER=10^3, THRES = 10^{-8}, \boldsymbol e_1 = (1,\boldsymbol 0_{p-1}')'$;\\
Functions: \\
\hspace{1cm} $\mu_i(\beta, \delta) = (Y_i - Z_i'\beta)^2 - (Y_i - Z_i'\delta)^2, i=1,2,...,n$;\\
\hspace{1cm} $\bar Q(\omega; \beta, \delta) = \sum_{i=1}^n|\mu_i(\beta, \delta)| \log\Bigg\{1 + \exp\Big(\big[2\{\mu_i(\beta, \delta) > 0\} -1\big]\big[ X_i'\omega -\gamma\big] \Big)\Bigg\}$;\\
 \begin{algorithmic}[1]
\Function{\texttt{Relax}}{~}
\State
Initialize $\bar \beta^{(0)} \leftarrow {\arg\min}_{\beta}\sum_{i=1}^n(Y_i - Z_i'\beta)^2$ through the least squares;
\State
Initialize $\bar \delta^{(0)}\leftarrow \bar \beta^{(0)} + e^{(0)}$, where $e^{(0)}\sim N(\boldsymbol 0_{d},I_{d})$;
\State
Generate $\bar\omega^{(0)} \sim U_{S^{p-1}}(\boldsymbol e_1, 1)$, $\bar\gamma^{(0)} \sim U(l_0, u_0)$; \hfill (See Algorithm \ref{al:uniform})
\For{m=1,2,...,MAXITER}
    \If{$|\bar\omega^{(m)} - \bar\omega^{(m-1)}|< THRES$}
    \State EXIT;
    \EndIf
    \State $(\bar\omega^{(m)}, \bar\gamma^{(m)}) = {\arg\min}_{(\omega,\gamma)\in \mathcal S^{p-1} \times [l_0,u_0]} \bar Q(\omega; \bar\beta^{(m-1)}, \bar\delta^{(m-1)})$ \\
    \; ~ \hfill (e.g., through Newton Raphson);
    \State $\bar \beta^{(m)} \leftarrow {\arg\min}_{\beta}\sum_{i=1}^n(Y_i - Z_i'\beta)^21\{X_i'\bar\omega^{(m)} - \bar\gamma^{(m)}\le 0\}$;
    \State $\bar \delta^{(m)}\leftarrow {\arg\min}_{\delta}\sum_{i=1}^n(Y_i - Z_i'\delta)^21\{X_i'\bar\omega^{(m)} - \bar\gamma^{(m)} > 0\}$;
\EndFor\\
\Return $\bar\omega^{(m)}$
\EndFunction
\end{algorithmic}
\end{algorithm}

\begin{algorithm}[h!]
\caption{$\texttt{UnifSphere}(\omega_0, C; m=1)$ returns $m \ge 1$ uniform samples from a collection of points on a unit (semi-)hypersphere of which angle from $\omega_0$ is at most $\frac \pi 2 C$.}
\label{al:uniform}
 Constants: $N = 10 m$;\\
 \begin{algorithmic}[1]
\Function{$\texttt{UnifSphere}$}{$\omega_0, C, m=1$}
\If{$C=1$}
    \State Draw $\mathring\omega \sim N_p(0,1)$;
    \State $\mathring\omega\leftarrow \frac{\mathring\omega}{\|\mathring\omega\|_2}$;
\Else
    \State Draw $N$ \textit{i.i.d.} vectors $\alpha_l \equiv (\alpha_{l,1},...,\alpha_{l,p}) \sim U\{(-\frac \pi 2 C, \frac \pi 2 C]^p\}, ~ l = 1,2,...,N$;
    \State Calculate the weights as $w_l = \prod_{j=1}^{p-1}\cos^{p-j}(\alpha_{l,j}), ~ l = 1,2,...,N$;
    \State Generate $\{\alpha_{(k)}:k=1,2,...,m\}$ by sampling each element out of $\{\alpha_l:l=1,2,...,N\}$ with probability proportional to $\{w_l: l=1,2,...,N\}$;
    \For {k=1,2,...,m}
        \State $\varphi_{(k),0} \leftarrow \frac \pi 2$;
        \State $\varphi_{(k),j} \leftarrow {\arg\tan}\big\{\frac{\alpha_{(k), j+1}}{\alpha_{(k), j}}\sin(\varphi_{(k),j-1})\big\}, j= 1,2,...,p-1$;
        \State $\varphi_{(k)} \leftarrow  (\varphi_{(k),1}, ..., \varphi_{(k),p})'$;
        \State $\varphi_{(k)} \leftarrow \varphi_{(k)} + \texttt{toAngle}(\omega_0)$
        \State $\omega_{(k)} \leftarrow \texttt{toCoordinate}(\varphi_{(k)})$;
    \EndFor
\EndIf\\
\Return $\{\omega_{(k)}: k=1,2,...,m\}$;
\EndFunction
 \end{algorithmic}
\end{algorithm}

\begin{algorithm}[h!]
\caption{$\texttt{toAngle}(\omega)$ returns the corresponding angle of the point $\omega$ from the first basis vector $\boldsymbol e_1$, and $\texttt{toCoordinate}(\varphi)$ returns the Cartesian position that corresponds to the angle $\varphi$.}
\label{al:transform}
 \begin{algorithmic}[1]
\Function{$\texttt{toAngle}$}{$\omega$}
\State $\varphi_{0} \leftarrow \frac \pi 2$;
\State $\varphi_{j} = {\arg\tan}\{\frac{\omega_{j+1}}{\omega_{j}}\sin(\varphi_{j-1})\}, j= 1,2,...,p-1$;
\Return $\varphi = (\varphi_1, \varphi_2, ..., \varphi_{p-1})'$;
\EndFunction
\Function{$\texttt{toCoordinate}$}{$\varphi$}
\State $\varphi_{0} \leftarrow 0$;
\State $\omega_{j} = \sin(\varphi_{p-j+1})\prod_{j'=0}^{j-1}\cos(\varphi_{p-j'}) j= 1,2,...,p$;
\Return $\omega = (\omega_1, \omega_2, ..., \omega_{p})'$;
\EndFunction
 \end{algorithmic}
\end{algorithm}

\begin{algorithm}[h!]
\caption{\texttt{MeanArg}() and \texttt{ModeArg}() find $(\hat\omega, \hat\gamma)$ and $(\check\omega, \check\gamma)$, respectively, given $(\tilde\omega, \tilde\gamma)$.}
\label{al:arg}
Constants: $RESOLUTION = 200$;\\

\begin{algorithmic}[1]
\Function{$\texttt{MeanArg}$}{$\tilde\omega, 
\tilde \gamma$}\\
\Return $\texttt{MidArg}(\tilde\omega, 
\tilde \gamma, \texttt{TYPE} = \texttt{MEAN})$;
\EndFunction
\end{algorithmic}
~ \\

\begin{algorithmic}[1]
\Function{$\texttt{ModeArg}$}{$\tilde\omega, 
\tilde \gamma$}\\
\Return $\texttt{MidArg}(\tilde\omega, 
\tilde \gamma, \texttt{TYPE} = \texttt{MODE})$;
\EndFunction
\end{algorithmic}
~\\

 \begin{algorithmic}[1]
\Function{$\texttt{MidArg}$}{$\tilde\omega, 
\tilde \gamma$, \texttt{TYPE}}
\For {$i=1,2,...,n$}
    \State $V_i \leftarrow \text{sign}(X_i'\tilde\omega - \tilde\gamma)$;
\EndFor
\State $N\leftarrow 0; M \leftarrow 10^4; C\leftarrow 1$;
\While {$N < RESOLUTION$}
    \State Draw $\{\omega_1, ...,\omega_M\} \leftarrow U_{S^{p-1}}(\tilde\omega, C, m = M)$;
    \For {$m=1,2,...,M$}
        \State $U_m \leftarrow \underset{i:i=1,...,n, V_i=1}{\min} X_i'\omega_m$;
        \State $L_m \leftarrow \underset{i:i=1,...,n, V_i=-1}{\max}X_i'\omega_m$;
        \State $R_m \leftarrow U_m - L_m$;
    \EndFor
    \State $N\leftarrow \sum_{m=1}^M 1(R_m > 0)$;
    \State $C \leftarrow 0.7 C$;
    \State $M \leftarrow 1.1 M$;
\EndWhile
\If {\texttt{TYPE} == \texttt{MEAN}}
    \State $\omega^* \leftarrow \sum_{m=1}^M\omega_m R_m$;
    \State $\omega^* \leftarrow\frac{\omega^*}{\|\omega^*\|_2}$;
    \State $\gamma^* \leftarrow \left\{\underset{i:i=1,...,n, V_i=1}{\min} X_i'\omega^*\right\}/2 + \left\{\underset{i:i=1,...,n, V_i=-1}{\max} X_i'\omega^*\right\}/2$;
\Else
    \State $m^* \leftarrow {\arg\max}_m R_m$;
    \State $\omega^* \leftarrow \omega_{m^*}$;
    \State $\gamma^* \leftarrow \frac{U_{m^*} + L_{m^*}}{2}$;
\EndIf
\Return $(\omega^*, \gamma^*)$;
\EndFunction
 \end{algorithmic}
\end{algorithm}

\FloatBarrier

\subsection{Mixed integer programming (MIP) algorithm specification}
\label{asec:mip}
We reformulate the estimation problem using a mixed integer program set-up as described in Section \ref{sec:simest_mip}. We first provide a baseline framework that is intended to give the true optimum, which however, could be impractical for large $p$ and $n$ settings. Then, we provide a coordinate descent algorithm of which high-level idea was introduced in \cite{lee2021factor}, but our algorithm has a distinct algebraic formulation from theirs.

\subsection{A baseline MIP framework}
\label{asec:mip_base}

We describe the algorithm-level description of the baseline MIP framework.
Define our decision variable as $\bs k=(\omega', \gamma, \beta', \eta', \bs V_n', \bs q_n, \bs s_n).$ The dimension of $\bs k$ is $d_{\bs k}=3n + 2d + p + 1 (= p + 1 + d + d + n + n + n \text{ in the order})$, 
and the variable is binary-coded (0 or 1) for the elements corresponding to $\bs V_n$, $\bs q_n$ and $\bs s_n$, and is real-valued for the others. 
The elements of $\bs k$  that correspond to $(\omega, \phi)$ are bounded by a constant, e.g., 1, and the other continuous decision variables are not bounded.

To construct a quadratic objective function, consider a matrix $\mathbb A$ of dimension $d_{\bs k} \times d_{\bs k}$, which is made up of $6\times 6$ block matrices $\{(\mathbb A(j,k):k=1,2,...,6):j=1,2,...,6\}$. The six block-rows ($j$'s) have dimensions $p +1, d, d, n, n$, and $n$, and so do the six block columns  ($k$'s). 
Then the objective is to minimize 
$$\bs k' \mathbb A \bs k,$$
where $\mathbb A(j,k)=I_n$ for $(j,k)=(5,5), (6,6)$, and  a zero matrix, otherwise.

The constraint $V_i = \ind(\omega'X_i-\gamma\le 0)$, for $i=1,2,...,n$, is expressed as 
$$\bs k' \mathbb B_i \bs k + \bs b_i' \bs k + 10^{-8} \le 0,$$
where $\mathbb B_i$ is another matrix of exactly the same dimension as $\mathbb A$ such that the 4th block-row and the 1st block-column matrix is $\mathbb B(4,1) = \begin{pmatrix} \bs 0_{(i-1),(p+1)}\\ 2(X_i', -1) \\\bs 0_{(n-i), (p+1)}\end{pmatrix}$ and all the other block matrices are zero, $\bs 0_{m, l}$ is an $m\times l$ zero matrix, $\bs b_i = (-X_i', 0, \bs 0_{2d+3n}')'$. A small value $10^{-8}$ was added on the left-hand side of the constraint to avoid degenerate solutions (i.e., solutions at the boundary that are within a tolerance limit but are non-feasible).

Finally, the constraint for the unit norm, or $\omega'\omega = 1$, is expressed as 
$$\bs k' \mathbb M \bs k > 0.3,$$ where
$\mathbb M$ is a diagonal matrix with the diagonal entries $(\bs 1_{p}', \bs 0_{1+2d+3n}')'$. Then, $\tilde \phi$ is given by $\tilde \phi = (k_1,...,k_p, k_{p+1})/\sqrt{\sum_{j=1}^p k_j^2}$ for the found solution $\bs k$. The reason why the last constraint was not given as $\bs k' \mathbb M \bs k =1$ is to give the software degree of freedom. The small positive right-hand side is to discourage the software from picking $(k_1,...,k_p) = \bs 0_p$.

The estimate for $\zeta$, or $\tilde \beta$ and $\tilde \delta$, are $(k_{p+2}, ..., k_{p+d+1})$ and $(k_{p+d+2}, ..., k_{p+2d+1})$ for the found solution $\bs k$.

\subsection{An MIP framework with the coordinate-descent algorithm}
The baseline MIP algorithm often takes hours for large $p$ and $n$ settings. Even if it finds an incumbent solution that turns out to be sample-optimum, MIP performs further searches until it collects evidence that there is no better solution. To circumvent this issue, \cite{lee2021factor} came up with a heuristic where profile maximum likelihood is performed by optimizing $\omega$ and $\gamma$ one at a time fixing the latest solution of the other. Below we describe the approach in more detail.

Define $\texttt{MIP}(\texttt{start} = \omega_{\texttt{start}}; \texttt{bestobjstop} = \epsilon_{\texttt{bestobjstop}})$ as the baseline MIP algorithm described in Section \ref{asec:mip_base}, which starts with $\omega_{\texttt{start}}$ as the incumbent solution and returns the first incumbent solution of which objective function is larger than that of  $\omega_{\texttt{start}}$ by $\epsilon_{\texttt{bestobjstop}}$. If there is no solution found that satisfies the condition, it returns the current best solution.

We will postpone the definition of the MIP-based profile objective maximizer algorithm to the next section and denote it as  $\texttt{MIP-PROFILE}(\texttt{start}=(\beta_1, \delta_1))$, where $(\beta_1, \delta_1)$ is the best solution found for $(\beta, \delta)$ up to the previous iteration.

\begin{algorithm}[h!]
\caption{MIP-CDA.}
\label{al:mip_cda}
Constants: $\epsilon_{OBJ}=10^{-5}$, $MAXITER=20$\\
\begin{algorithmic}[1]
\State $\omega^{(0)}\leftarrow~\texttt{Relax}()$, $\text{OBJ}^{(0)}\leftarrow \infty$;  \hfill (See Algorithm \ref{al:relax})
\State $(\omega^{(1)}, \gamma^{(1)}, \beta^{(1)}, \delta^{(1)})\leftarrow~\texttt{MIP}(\texttt{start}=\omega^{(0)}; \texttt{bestobjstop}=10^{-2})$;
\State $\text{OBJ}^{(1)} \leftarrow M_n(\omega^{(1)}, \gamma^{(1)}, \beta^{(1)}, \delta^{(1)})$;
\State $(\omega^*, \gamma^*, \beta^*, \delta^*) \leftarrow (\omega^{(1)}, \gamma^{(1)}, \beta^{(1)}, \delta^{(1)})$
\For{$k = 2, 3,..., MAXITER$}
    \If {$|\text{OBJ}^{(k)} - \text{OBJ}^{(k-1)}| < \epsilon_{OBJ}$} 
        \State EXIT;
    \EndIf
    \State $(\omega^{(k)}, \gamma^{(k)})\leftarrow~\texttt{MIP-PROFILE}(\texttt{start}=(\beta^{(k-1)}, \delta^{(k-1)}))$;   \hfill (See Section \ref{asec:mip_profile})
    \State $\beta^{(k)} \leftarrow {\arg\min}_{\beta}\sum_{i=1}^n(Y_i - Z_i'\beta)^21\{X_i'\omega^{(k)} - \gamma^{(k)}\le 0\}$;  \hfill (Least squares)
    \State $\delta^{(k)}\leftarrow {\arg\min}_{\delta}\sum_{i=1}^n(Y_i - Z_i'\delta)^21\{X_i'\omega^{(k)} - \gamma^{(k)} > 0\}$; \hfill (Least squares)
    \State $\text{OBJ}^{(k)} \leftarrow M_n(\omega^{(k)}, \gamma^{(k)}, \beta^{(k)}, \delta^{(k)})$;
    \State $(\omega^*, \gamma^*, \beta^*, \delta^*) \leftarrow (\omega^{(k)}, \gamma^{(k)}, \beta^{(k)}, \delta^{(k)})$
\EndFor  
\State $\hat \omega, \hat \gamma \leftarrow \texttt{MeanArg}(\omega^*, \gamma^*)$; \hfill (See Algorithm \ref{al:arg})
\State $\check \omega, \check \gamma \leftarrow \texttt{ModeArg}(\omega^*, \gamma^*)$;  \hfill (See Algorithm \ref{al:arg})
\State $\hat \beta \leftarrow \check \beta \leftarrow \beta^*$;
\State $\hat \delta \leftarrow \check\delta \leftarrow \delta^*$;
\end{algorithmic}
\end{algorithm}

\subsection{A profile MIP framework}
\label{asec:mip_profile}

We describe the algorithm-level description of the profile MIP framework that, given fixed $\beta_1$ and $\delta_1$, finds $\tilde\phi(\beta_1, \delta_1) = \arg\min_{\omega\in\mathcal S^{p-1}, \gamma\in \mathbb R} M_n(\omega, \gamma, \beta_1, \delta_1)$.

It can be shown that the objective function can be set such that
$$\tilde\phi(\beta_1, \delta_1) = {\text{argmin}}_{\omega\in\mathcal S^{p-1}, \gamma\in \mathbb R} \sum_{i=1}^n \Lambda_i(\zeta_1) \ind\{\omega_0'X-\gamma_0\leq 0\},$$
where we recall the notation $\zeta_1 = (\beta_1, \eta_1)$, and 
$\Lambda_i(\zeta) = (Y_i - \beta'Z_i)^2- (Y_i - \delta'Z_i)^2$ for $i=1,2,...,n$.

Define our decision variable as $\bs k=(\omega', \gamma, \bs V_n')$ with dimension $d_{\bs k}=n + p + 1$. The variables are binary-coded (0 or 1) for the elements corresponding to $\bs V_n$ and are real-valued for the others. 
The elements of $\bs k$  that correspond to $(\omega, \phi)$ are bounded by constants, e.g., (-1, 1). Unlike the baseline MIP, this has a linear objective function, and we minimize $\bs k' \bs a$ with 
$$a=\left( \{\Lambda_i(\zeta_1):i=1,2,...,n\}', \bs 0_{p+1}'\right).$$

The constraints still involve quadratic computations. The constraint $V_i =\ind(\omega'X_i-\gamma\le 0)$, for $i=1,2,...,n$, is expressed as 
$$\bs k' \mathbb B_i \bs k + \bs b_i' \bs k + 10^{-8} \le 0,$$
where $\mathbb B_i$ is a $(n+p+1)\times(n+p+1)$ matrix of which $(n+i, l)$th entry is $X_{i,l}$ for $l=1,...,p$ and -1 for $l=p+1$, and zeros everywhere else. Specifically, all rows other than the $i$th row are zero vectors.
$\bs b_i = (X_i', -1, \bs 0_n')'$. Again, a small value $10^{-8}$ was added to the left-hand side of the constraint to avoid degenerate solutions.

The constraint for the unit norm, or $\omega'\omega = 1$, is expressed as 
$$\bs k' \mathbb M \bs k > 0.3,$$ where
$\mathbb M$ is a diagonal matrix with the diagonal entries $(\bs 1_{p}', \bs 0_{1+2d+3n}')'$. Then, $\tilde \phi$ is given by $\tilde \phi = (k_1,...,k_p, k_{p+1})/\sqrt{\sum_{j=1}^p k_j^2}$ for the found solution $\bs k$.
\FloatBarrier

\section{Simulation study}
\label{asec:sim}
\begin{figure}
    \centering
    \includegraphics[width = \textwidth] {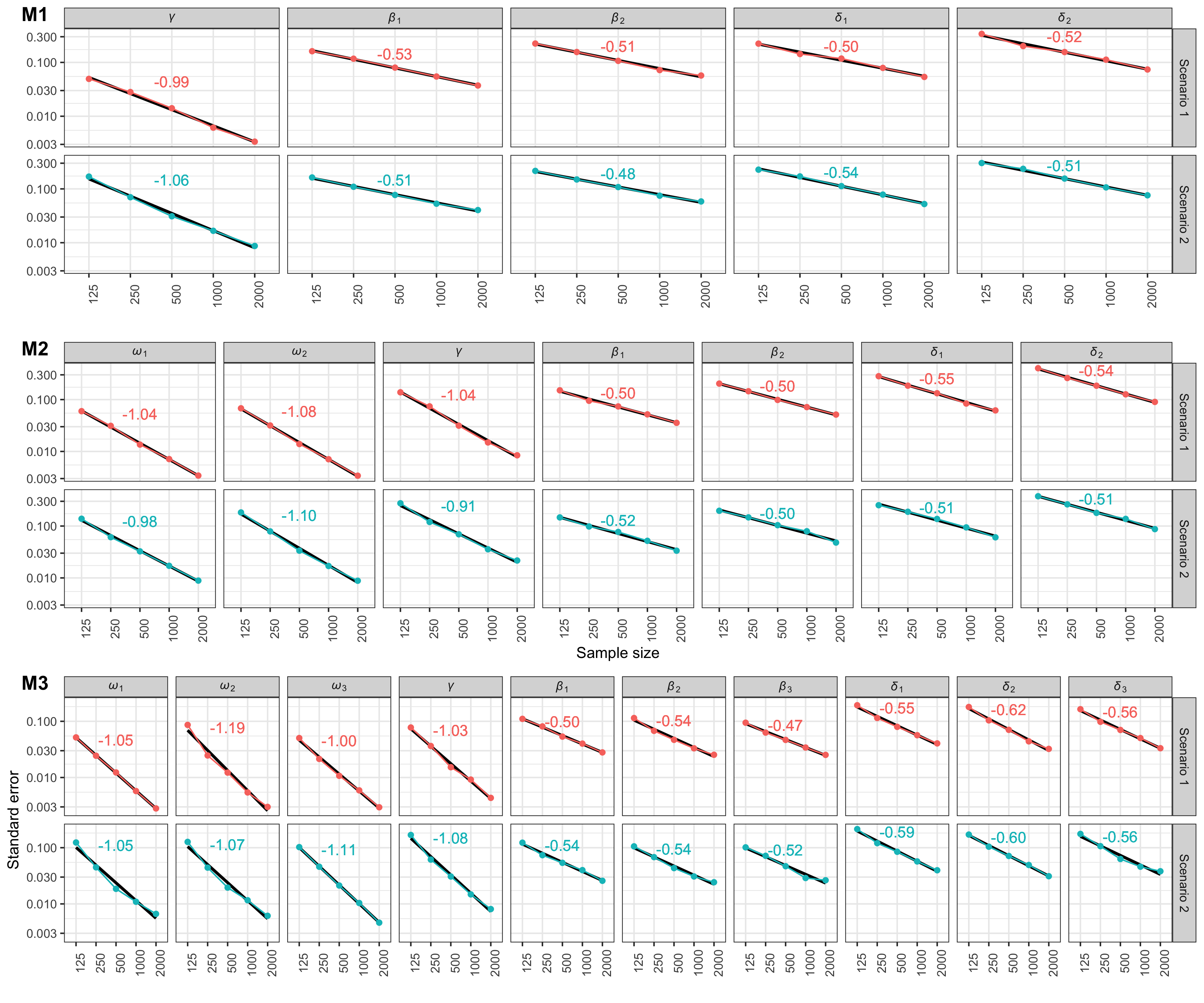}
    \caption{Rate of convergence results for mean-argmin estimators (based on the CDA-MIP algorithm). Both standard errors ($y$-axis) and sample sizes ($x$-axis) are presented on the log scale (with base 2). The lines and the annotated numbers are the least squares regression and the corresponding slope estimates representing the exponents of the rate of convergence. M1, Model 1; M2, Model 2; M3, Model 3}
    \label{fig:rate2}
\end{figure}

Figure \ref{fig:rate2} illustrates that the estimation error of all estimates, based on the CDA-MIP algorithm, becomes polynomially smaller as the sample size becomes larger for all models. The theoretical rate of convergence---$n^{-1}$ for $\hat\phi$ and $\check\phi$ and $n^{-1/2}$ for $\hat\zeta$ and $\check\zeta$---is well observed in the numerical results in most settings.

\begin{figure}
    \centering
    \includegraphics[width = 0.8\textwidth]{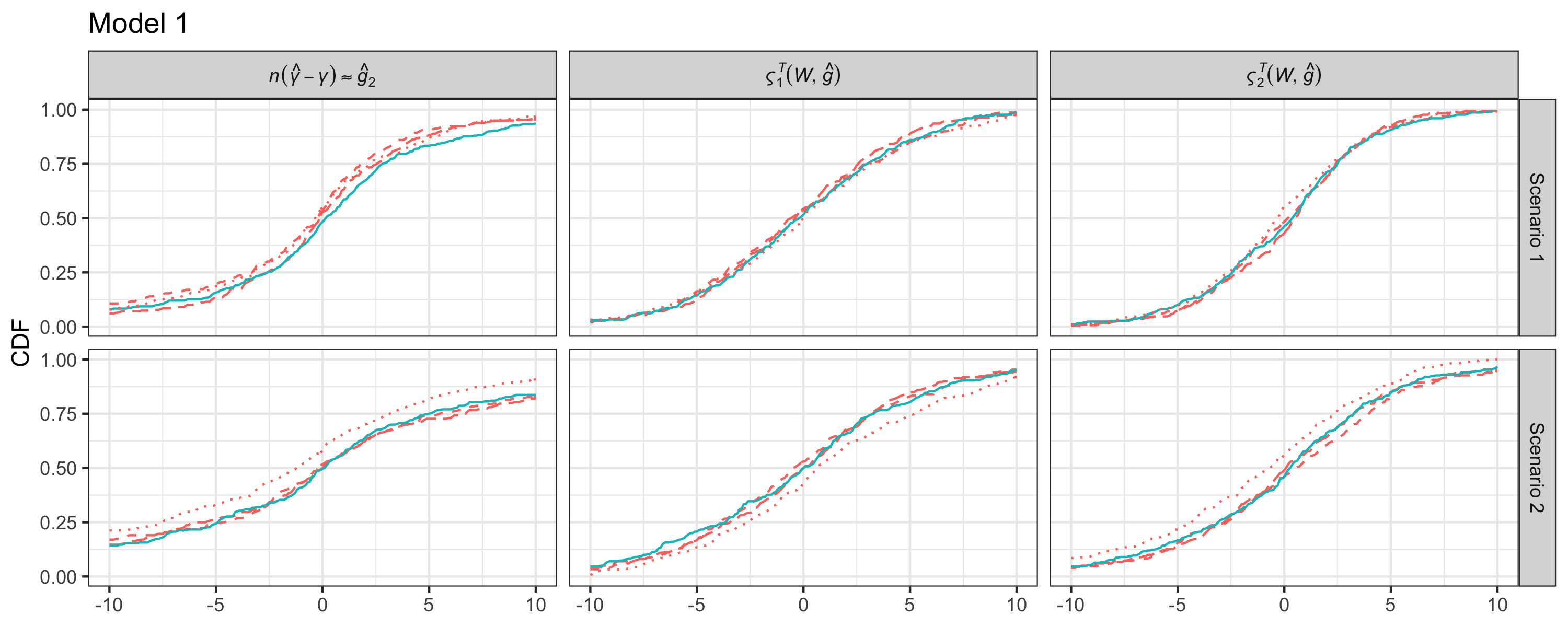}
    \includegraphics[width = 0.8\textwidth]{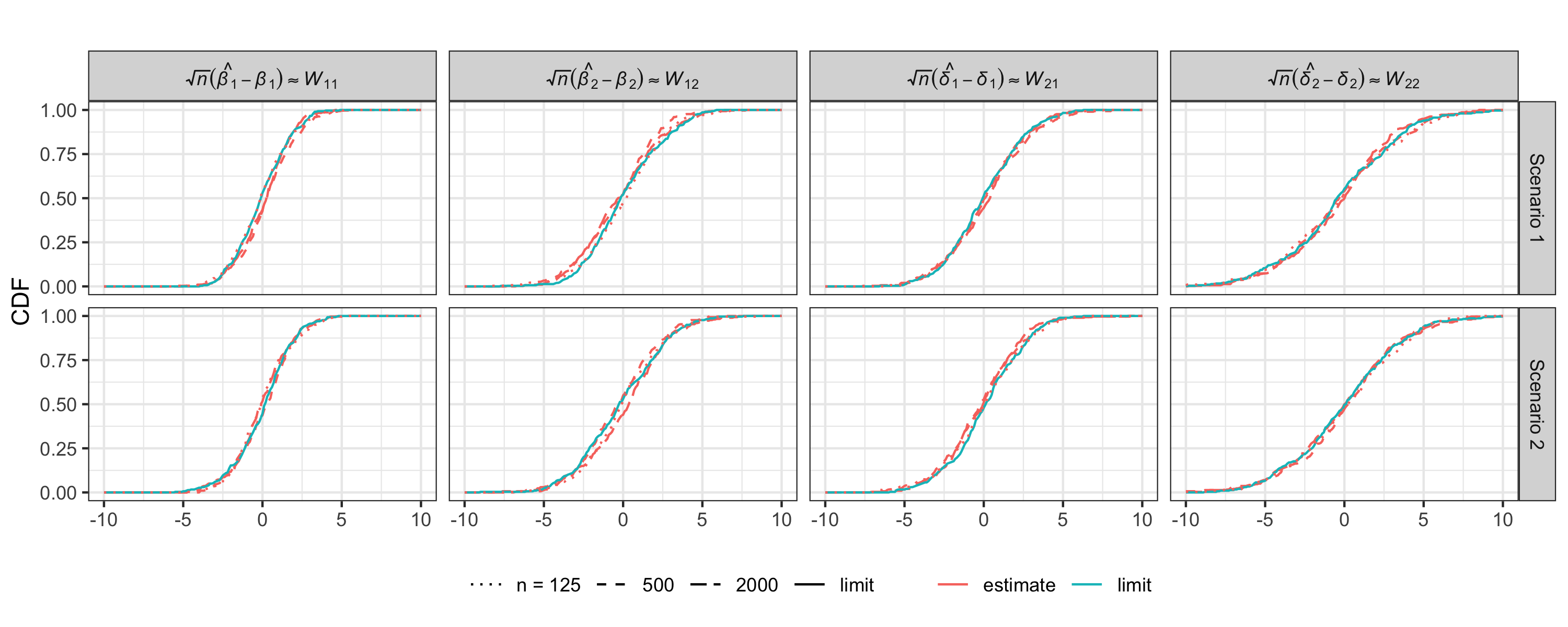}
    \includegraphics[width = 0.8\textwidth]{graph/ChangePlane_Model2_weakConv_mean_cda_phi.png}
    \includegraphics[width = 0.8\textwidth]{graph/ChangePlane_Model2_weakConv_mean_cda_zeta.png}
    \caption{The mean-argmin estimated and limiting CDFs for Models 1 and 2. The random coefficients are $\varsigma_1 = (-0.47, -0.26, 0.15, 0.82, -0.60, 0.80)'$, $\varsigma_2 = (0.89, 0.32, 0.26, -0.88, -0.59, -0.65)'$ for Model 1 and $\varsigma_1 = (-0.63, 0.40, 0.15, -0.66, 0.89, 0.89, -0.74)'$, $\varsigma_2 = ( 0.67, -0.06, 0.10, 0.11, -0.52, 0.52, -0.64)'$ for Model 2.}
    \label{fig:weak2}
\end{figure}

\begin{figure}
    \centering
    \includegraphics[width = 0.8\textwidth]{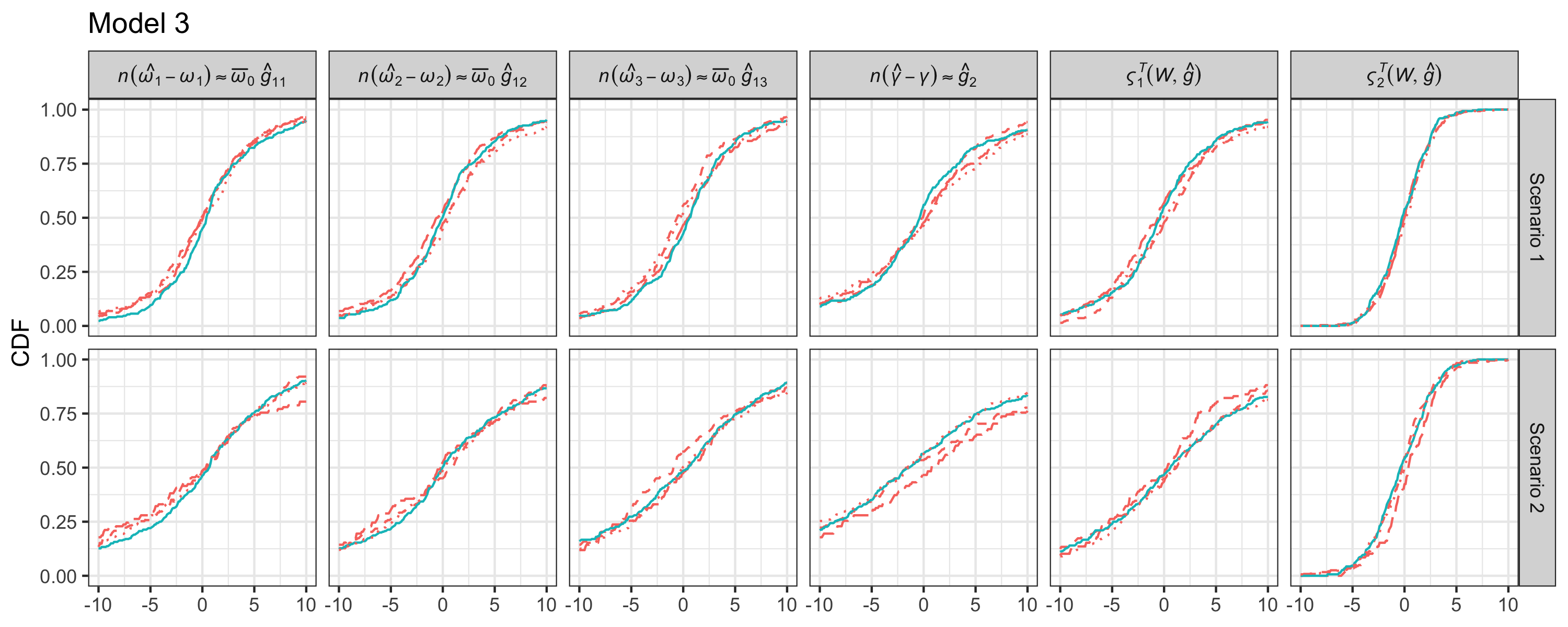}
    \includegraphics[width = 0.8\textwidth]{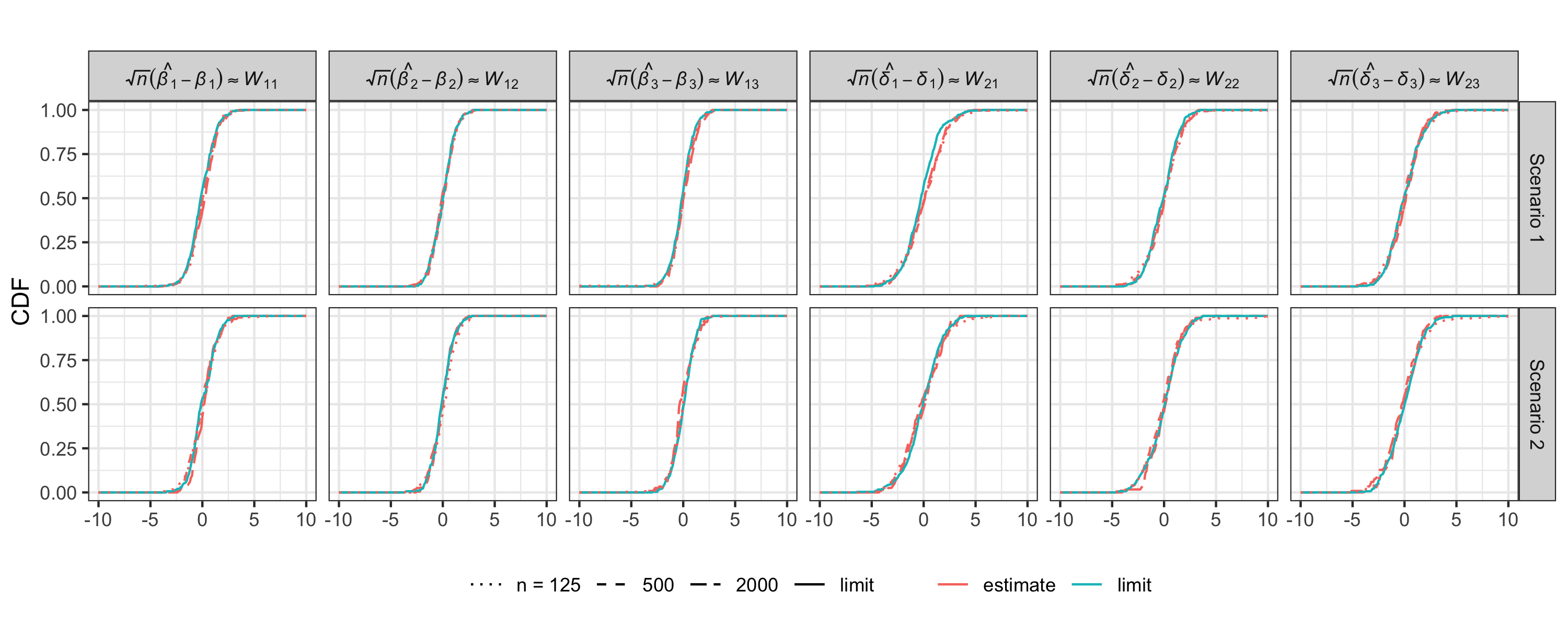}
    \includegraphics[width = 0.8\textwidth]{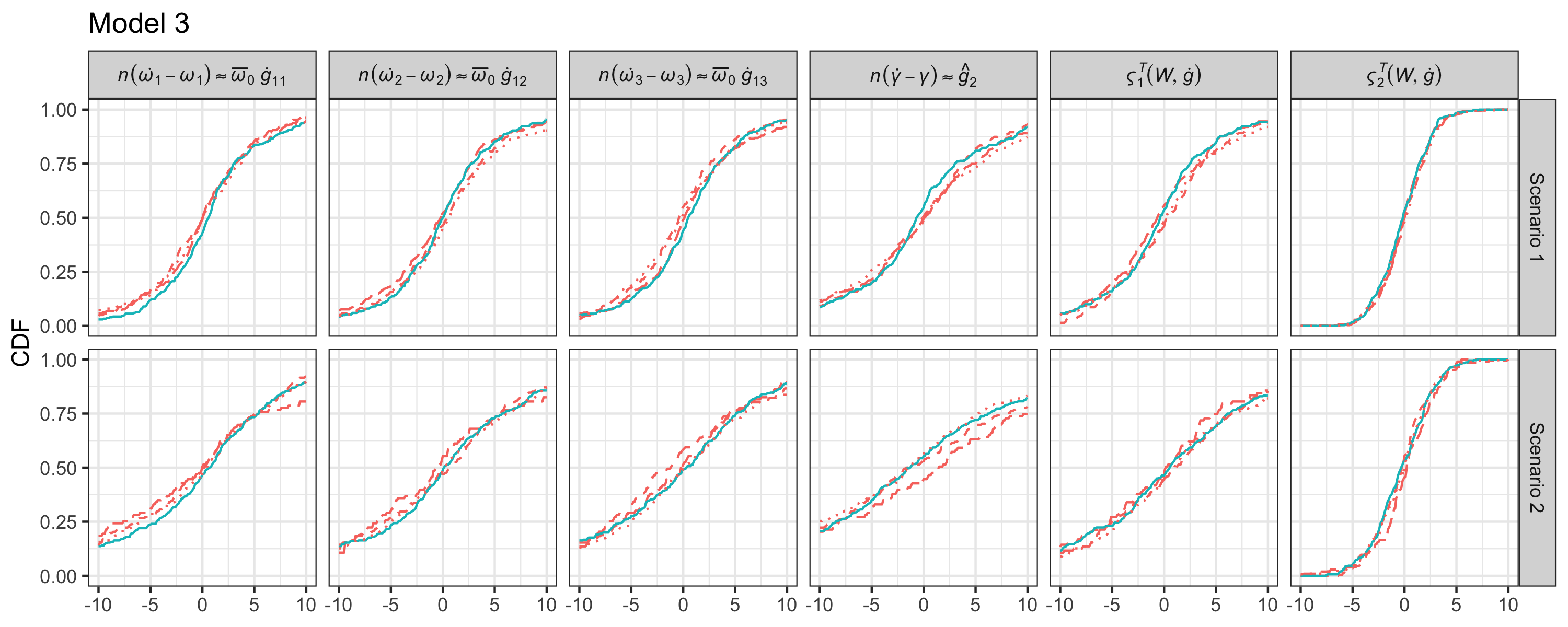}
    \includegraphics[width = 0.8\textwidth]{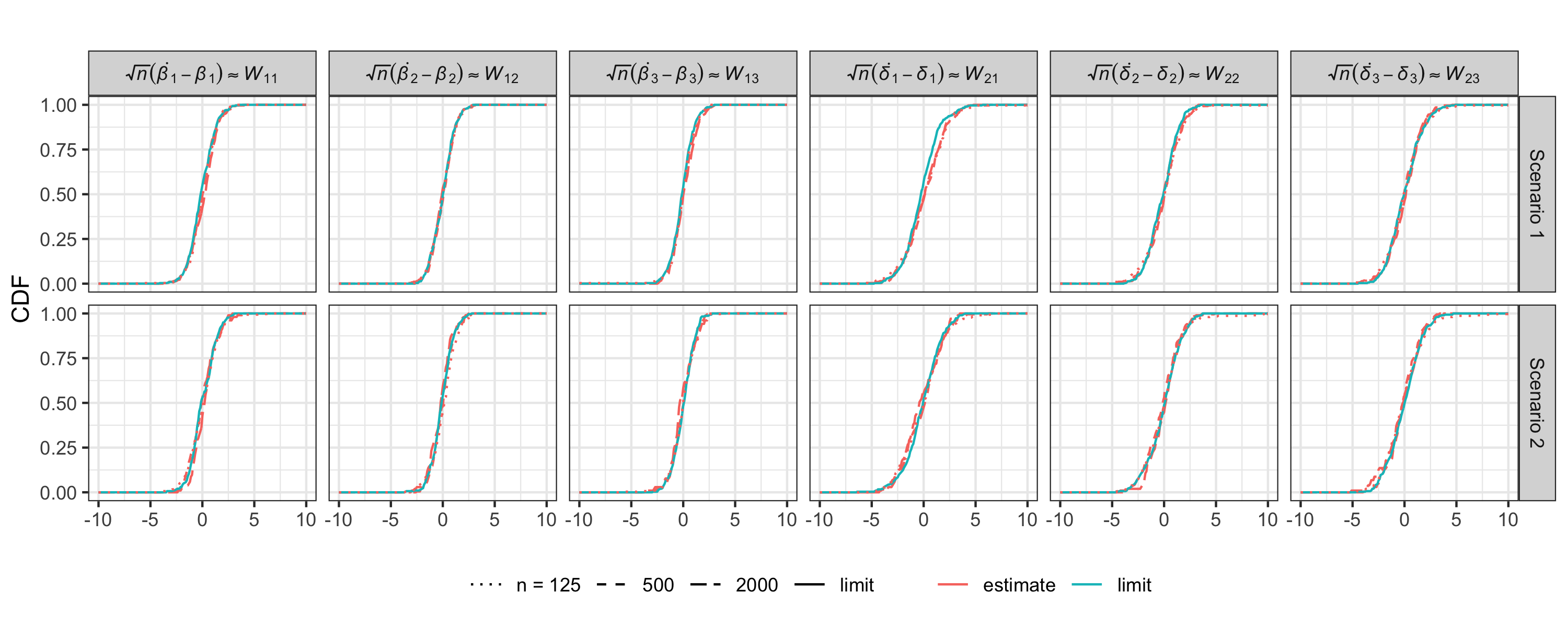}
    \caption{The mean-argmin (top two) and mode-argmin (bottom two) estimated and limiting CDFs for Model 3. The random coefficients $\varsigma_1 = (-0.66, 0.62, -0.23, -0.34, 0.20, 0.21, -0.75, -0.41, 0.16, 0.26)'$, $\varsigma_2 = (0.02, 0.01, 0.07, 0.11, 0.74, 0.66, -0.78, 0.41, 0.79, -0.44)'$ for Model 3.}
    \label{fig:weak3}
\end{figure}

In Figure \ref{fig:weak2}, the weak convergence simulation results for mean argmin estimators are presented for Models 1 and 2, and Figure \ref{fig:weak3} contains the results for Model 3 where both mean and mode argmin estimators are represented. Note that, for the change-point model (Model 1), mean- and mode-argmins are identical, and, for the binary change-plane covariate model (Model 2), $n\check\phi \not\leadsto (\bar \omega_0'\check g_1, \check g2)$, as expected, as Assumption $C3'$ is not satisfied. To see the joint convergence, per the Cram\'er-Wold device, the CDF of linear combinations of the estimates with random coefficients were compared to the corresponding limiting distribution. 
The numerical results support the weak convergence theory.

\begin{figure}
    \centering
    \includegraphics[width = 0.3\textwidth]{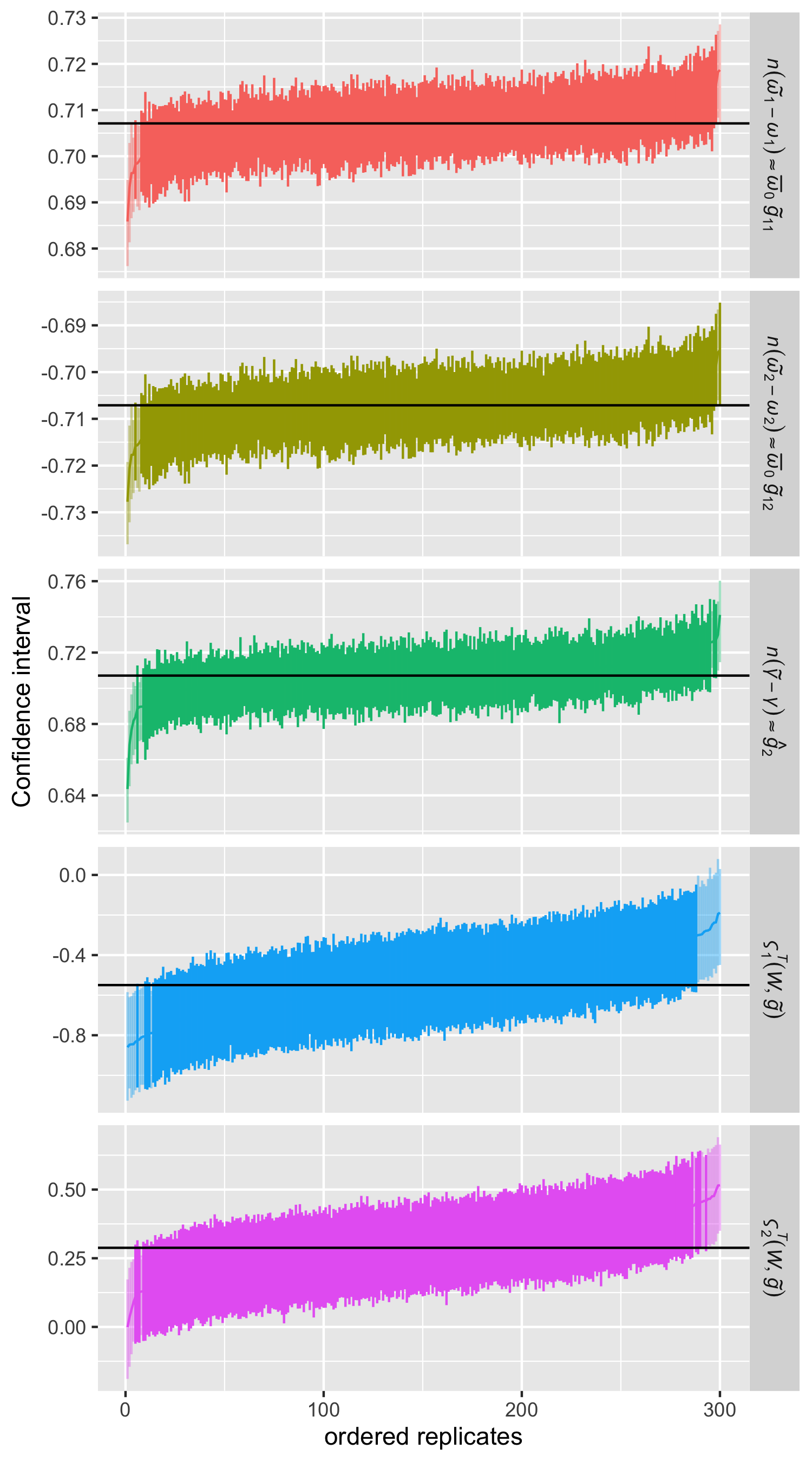}
    \includegraphics[width = 0.3\textwidth]{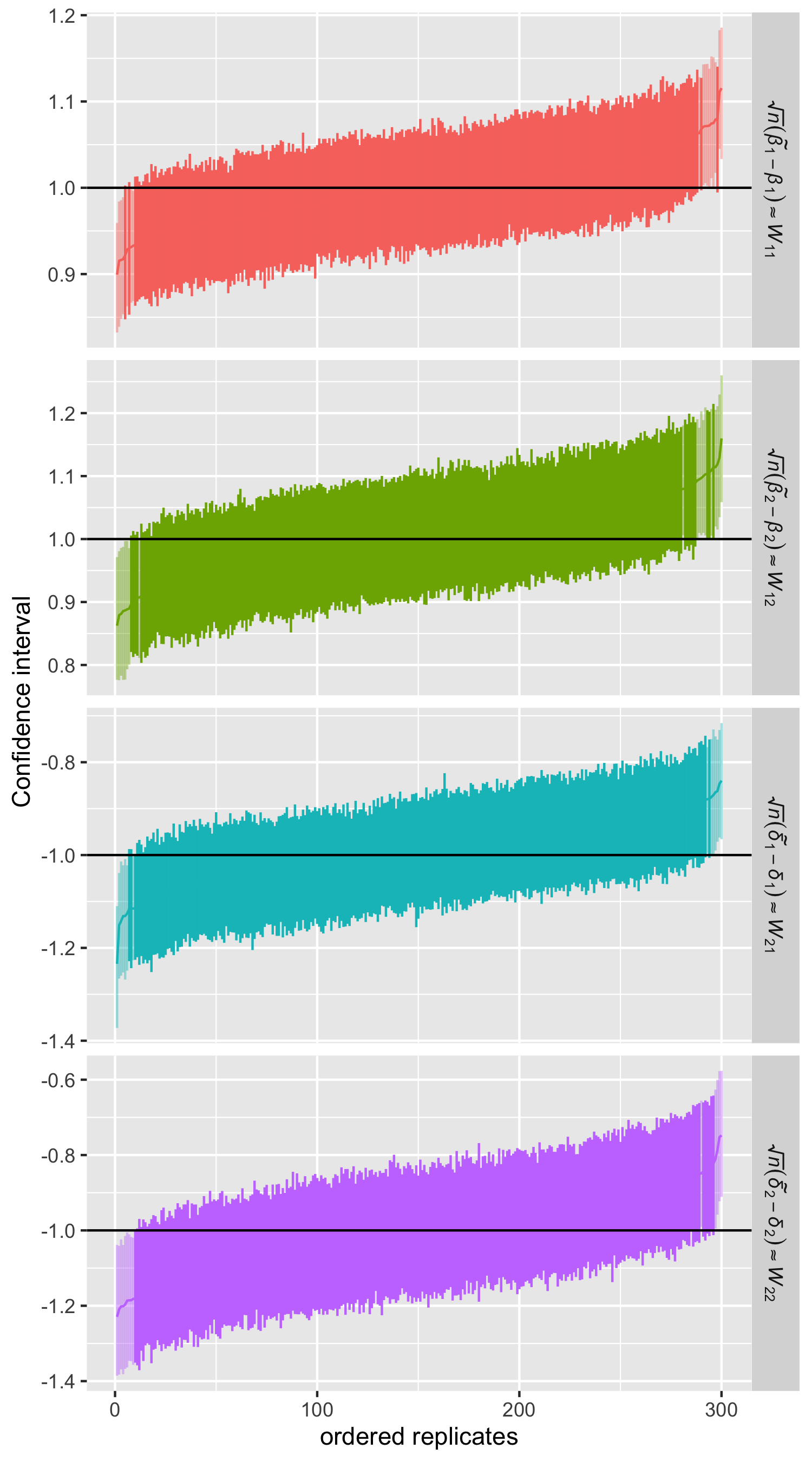}
    \caption{Vertical confidence intervals (CIs) sorted by their point estimates along the X-axis for Model 2 Scenario 1 with $n=2000$. The light (i.e., low opacity) CI bars represent no-coverage.}
    \label{fig:ci_plot}
\end{figure}

To examine the operating characteristics of the parametric bootstrap estimator, we provide, in Figure \ref{fig:ci_plot}, a confidence interval diagram for Model 2 Scenario 1 with $n=2000$.
Table \ref{tab:alternative_power} further shows the coverage probabilities for alternative values of the parameters. E.g., the estimated confidence intervals' coverage probabilities are presented for $\beta_1 - 0.04, \beta_1 - 0.02, \beta_1, \beta_1+0.02,$ and $\beta_1+0.04$. The same values corresponding to other parameters are also presented. With much faster rates than $\zeta$ estimators, the $\phi$ estimators have much higher specificity to shifts in the parameter values.

\begin{table}[ht]
\label{tab:alternative_power}
\centering
\begin{tabular}{cccccc}
  \hline
shift & -0.04 & -0.02 & 0 & +0.02 & +0.04\\ 
  \hline
   $\beta_1$    & 0.84 & 0.92 & 0.94 & 0.91 & 0.77 \\ 
   $\beta_2$    & 0.86 & 0.92 & 0.94 & 0.94 & 0.89 \\ 
   $\delta_1$   & 0.91 & 0.94 & 0.95 & 0.95 & 0.92 \\ 
   $\delta_2$   & 0.94 & 0.95 & 0.95 & 0.94 & 0.90 \\ 
   $\gamma$     & 0.02 & 0.49 & 0.97 & 0.42 & 0.02 \\ 
   $\omega_1$   & 0.00 & 0.01 & 0.97 & 0.01 & 0.00 \\ 
   $\omega_2$   & 0.00 & 0.01 & 0.98 & 0.01 & 0.00 \\ 
   $\varsigma_1$& 0.93 & 0.93 & 0.93 & 0.92 & 0.92 \\ 
   $\varsigma_2$& 0.92 & 0.94 & 0.95 & 0.96 & 0.92 \\ 
   \hline
\end{tabular}
\caption{The coverage probabilities for alternative values indicated by ``shift."}
\end{table}

\FloatBarrier

\subsection{Results for the uniform search algorithm}
We provide the results for the uniform search algorithm described in Section \ref{sec:simest_uniform} (see Algorithm \ref{al:uniform_search}).

\begin{figure}
    \centering
    \includegraphics[width = \textwidth] {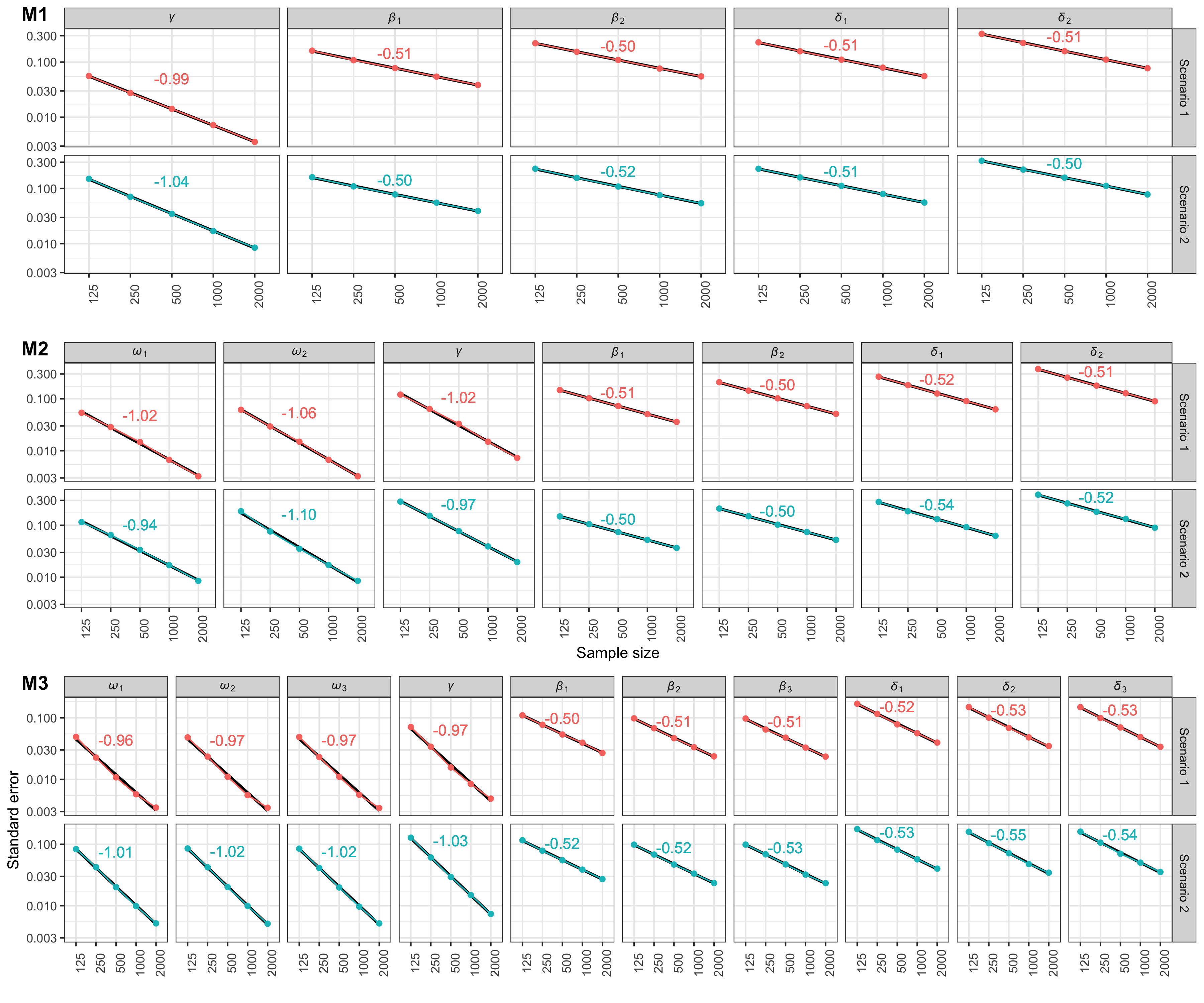}
    \caption{Rate of convergence results for mean-argmin estimators (based on the uniform search algorithm). Both standard errors ($y$-axis) and sample sizes ($x$-axis) are presented on the log scale (with base 2). The lines and the annotated numbers are the least squares regression and thse corresponding slope estimates represent the exponents of the rate of convergence. M1, Model 1; M2, Model 2; M3, Model 3}
    \label{fig:rate_uniform}
\end{figure}

Figure \ref{fig:rate_uniform} illustrates that the estimation error of all estimates, based on the uniform search algorithm, becomes polynomially smaller as the sample size becomes larger for all models. The theoretical rate of convergence---$n^{-1}$ for $\hat\phi$ and $\check\phi$ and $n^{-1/2}$ for $\hat\zeta$ and $\check\zeta$---is well observed in the numerical results in most settings.

\begin{figure}
    \centering
    \includegraphics[width = 0.4\textwidth]{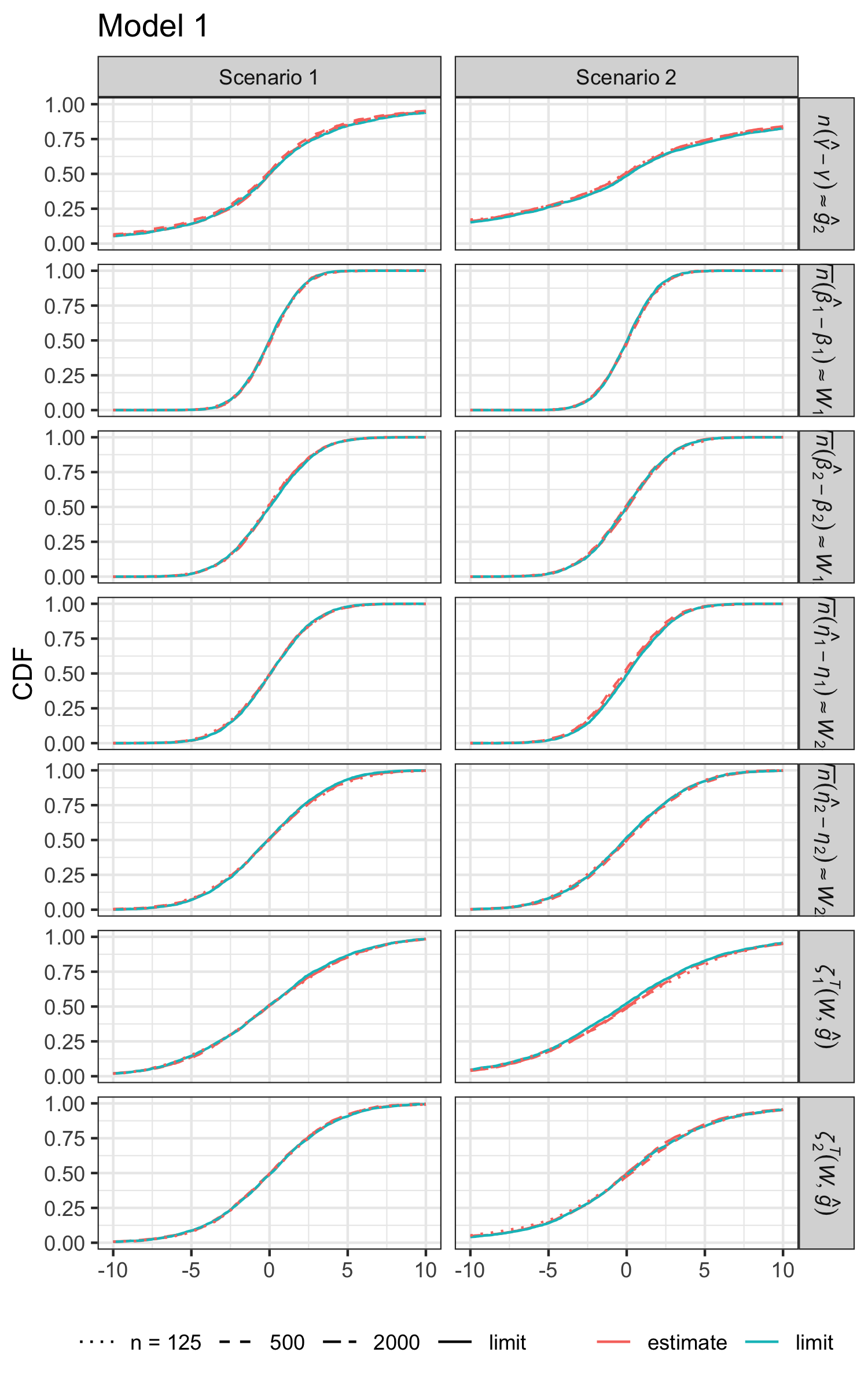}
    \includegraphics[width = 0.4\textwidth]{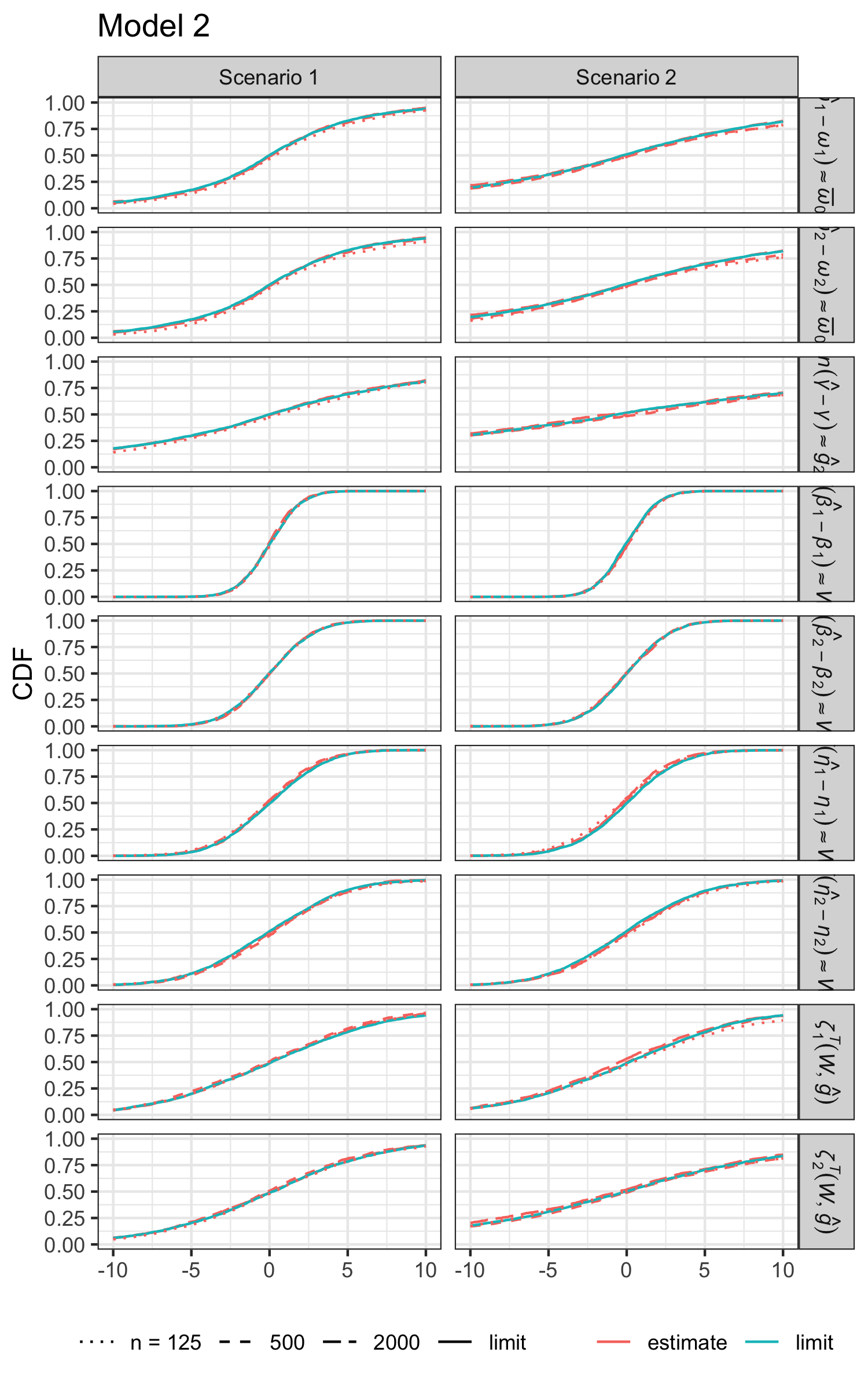}
    \caption{The mean-argmin estimated and limiting CDFs for Models 1 and 2---based on the uniform search algorithm. The random coefficients are $\varsigma_1 = (-0.47, -0.26, 0.15, 0.82, -0.60, 0.80)'$, $\varsigma_2 = (0.89, 0.32, 0.26, -0.88, -0.59, -0.65)'$ for Model 1 and $\varsigma_1 = (-0.63, 0.40, 0.15, -0.66, 0.89, 0.89, -0.74)'$, $\varsigma_2 = ( 0.67, -0.06, 0.10, 0.11, -0.52, 0.52, -0.64)'$ for Model 2.}
    \label{fig:weak2_uniform}
\end{figure}

We present the weak convergence results based on the uniform search algorithm. In Figure \ref{fig:weak2_uniform}, the weak convergence simulation results for mean argmin estimators are presented for Models 1 and 2, and Figure \ref{fig:weak3} contains the results for Model 3 where both mean and mode argmin estimators are represented. Note again that, for the change-point model (Model 1), mean- and mode-argmins are identical, and, for the binary change-plane covariate model (Model 2), $n\check\phi \not\leadsto (\bar \omega_0'\check g_1, \check g2)$, as expected, as Assumption $C3'$ is not satisfied. To see the joint convergence, per the Cram\'er-Wold device, the CDF of linear combinations of the estimates with random coefficients were compared to the corresponding limiting distribution. The numerical results based on the uniform search algorithm also supports the weak convergence theory.

\begin{figure}
    \centering
    \includegraphics[width = 0.4\textwidth]{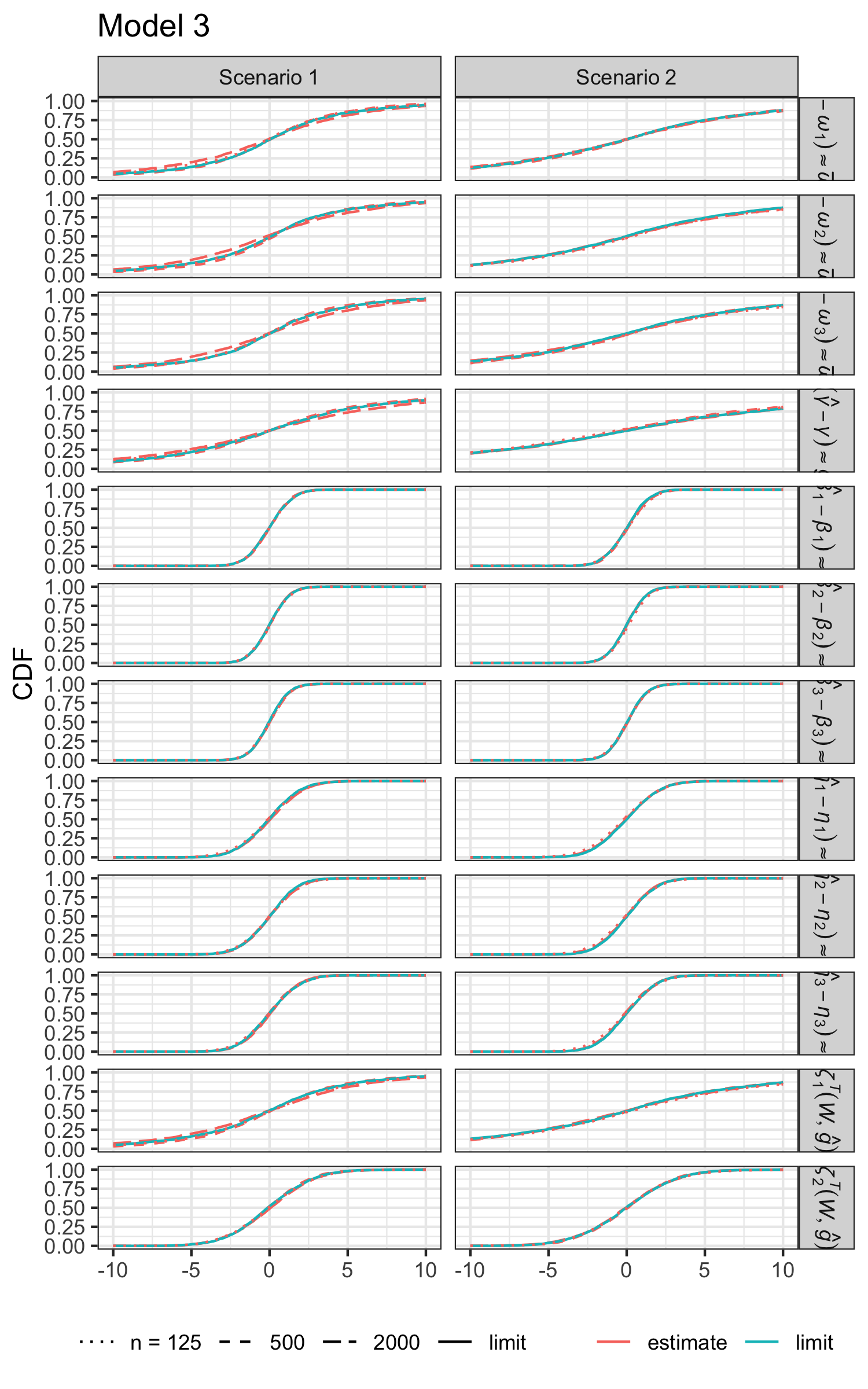}
    \includegraphics[width = 0.4\textwidth]{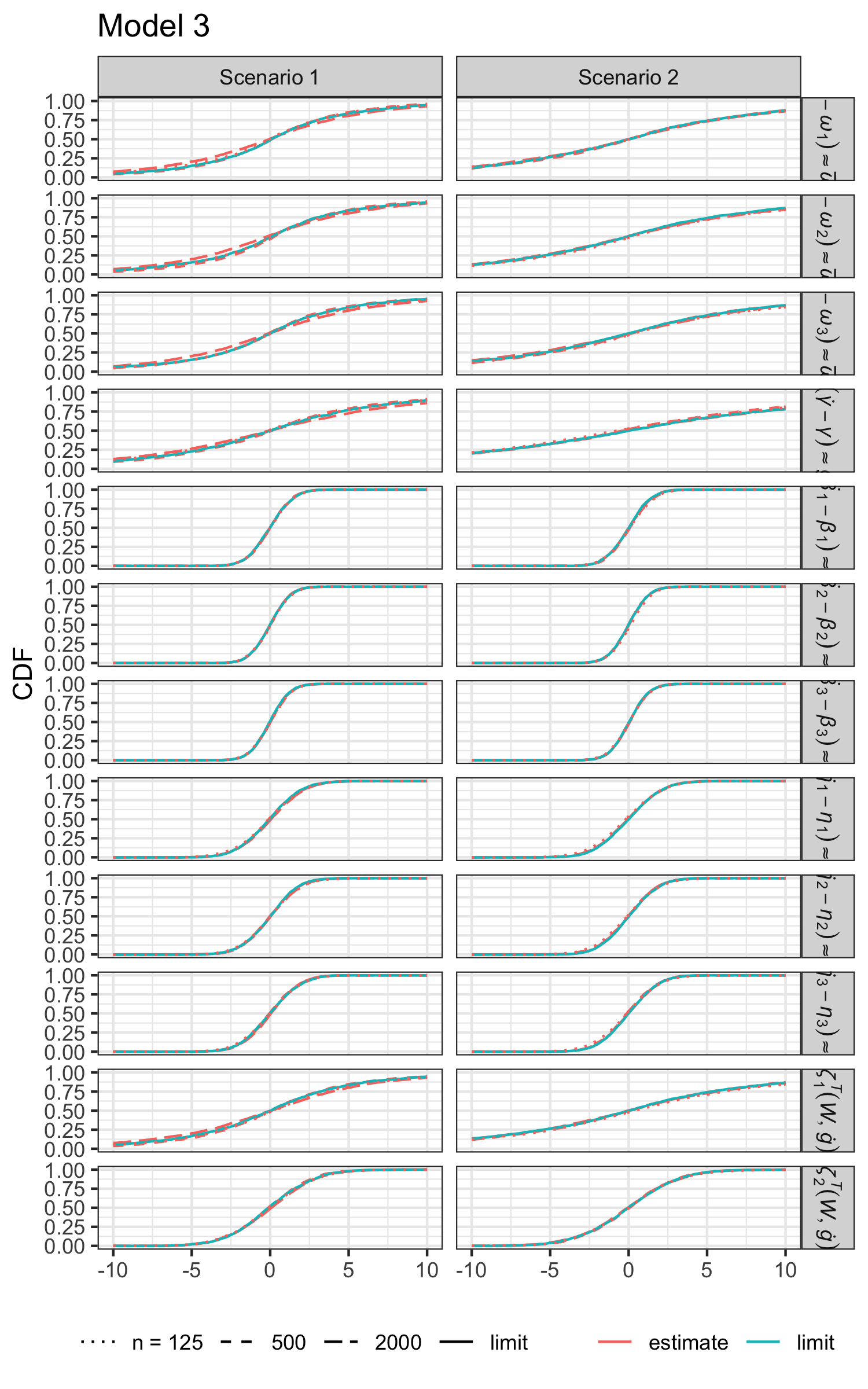}
    \caption{The mean-argmin (top two) and mode-argmin (bottom two) estimated and limiting CDFs for Model 3---based on the uniform search algorithm. The random coefficients $\varsigma_1 = (-0.66, 0.62, -0.23, -0.34, 0.20, 0.21, -0.75, -0.41, 0.16, 0.26)'$, $\varsigma_2 = (0.02, 0.01, 0.07, 0.11, 0.74, 0.66, -0.78, 0.41, 0.79, -0.44)'$ for Model 3.}
    \label{fig:weak3_uniform}
\end{figure}

\FloatBarrier

\subsection{Results for the group-balance settings}
\label{asec:balance}

To examine a finite sample estimation quality for imbalanced group sizes, we further consider two more scenarios, which are modifications of Model 2 Scenario 1. Namely, a perfectly balanced scenario (Scenario 3) with $(\omega,\gamma) =(1, -1, 0.5)/ \sqrt{2}$, and a highly imbalanced scenario (Scenario 4) with $(\omega,\gamma) = (1, -1, 1.5)/ \sqrt{2}$; all other parameters (the distribution of $X$ and $Z$ and $(\beta,\delta)$) are the same as Scenario 1 of Model 2, and the expected sample sizes per group are (250, 250) for Scenario 3 ($\rho = 0.5$) and (417, 83) for Scenario 4 ($\rho = 0.37$). The results are presented in Table \ref{tab:balance}. The overall estimation quality is worse for the unbalanced setting (Scenario 4) than the other (Scenario 3). The RMSE for $\zeta\equiv(\beta,\delta)$ is worse by a factor of 1.2--1.3 for the unbalanced setting, and the RMSE for $\phi$ is also worse by a factor of 1.0--1.8 for the unbalanced setting.

\begin{table}[ht]
\centering
\begin{tabular}{ccccccccccccccc}
  \hline
  $n$ & \multicolumn{3}{c}{125}& \multicolumn{3}{c}{500}& \multicolumn{3}{c}{2000}\\
 RMSE & $\hat\phi$  &$\check\phi$ & $\hat\zeta ~(=\hat\zeta)$ & $\hat\phi$  &$\check\phi$ & $\hat\zeta ~(=\hat\zeta)$ & $\hat\phi$  &$\check\phi$ & $\hat\zeta ~(=\hat\zeta)$\\ \hline
Model 2 Scenario 3 & 0.143 & 0.181 & 0.459 & 0.035 & 0.043 & 0.236 & 0.0070 & 0.0080 & 0.116\\
Model 2 Scenario 4 & 0.141 & 0.178 & 0.609 & 0.041 & 0.047 & 0.297 & 0.0123 & 0.0136 & 0.141\\
relative efficiency  & $\times 0.99$ & $\times 0.98$ &  $\times 1.33$ & $\times 1.17$ & $\times 1.10$ &  $\times 1.26$ & $\times 1.75$ & $\times 1.70$ &  $\times 1.22$ \\
  \hline
\end{tabular}
\label{tab:balance}
\caption{The root mean squared error (RMSE) of each estimator under balanced and imbalanced scenarios. The relative efficiency is defined as the RMSE of Scenario 4 divided by that of Scenario 3.}
\end{table}

\FloatBarrier

 \section{Application to the ACTG175 AIDS study}
 \label{asec:data_ACTG175}
We denote potential outcomes as $Y^{(1)}$ and $Y^{(0)}$, which are the outcomes that would have been observed if the patient 
were given Treatment $T=1$  or $T=0$, respectively. Using the standard counterfactual causal inference framework \citep{rubin2005causal, hernan2020},  
the treatment effect for the patient with individual characteristics $X=x$ is defined as the difference between the expected potential outcomes given $X=x$, that is,
\begin{equation*}
\Delta(x)=E\{Y^{*}(1)  -Y^{*}(0)|X=x\} =  E\{Y^{*}(1)|X=x\}-E\{Y^{*}(0)|X=x\}.
\end{equation*}
$\Delta(x)$  is also known as the conditional average treatment effect (CATE). Under the assumptions mentioned in Section \ref{sec:data_ACTG175} (SUTVA, no unmeasured confounders, positivity assumptions), $E\{Y^{*}(t)|X=x\} =E\{Y^{*}(t)|T=t, X=x\}=E\{Y^{*}(t)|T=1-t, X=x\} = E\{Y|T=t, x\}$ for $t=0, 1$ \citep{rubin2005causal, hernan2020}. The first two equations are true under the no unmeasured confounders (exchangeability) and positivity assumptions. The last equation is true under the consistency assumption, which is included in SUTVA. That is, the potential outcomes are independent of the actual treatment given $X=x$. Therefore, the conditional expectation of the outcome in the group treated with $T=t$ is the same as the conditional expectation of the outcome if every individual in the population of interest were treated with $T=t$, given $X=x$. We note that the independence between the potential outcomes $Y^{*}(t)$ and the actual treatment $T$ does not imply the independence of the observed outcomes $Y$ and treatment $T$. We refer to Sections 2.1, 2.2, and 3.3 in \cite{hernan2020} for further details. 

By plugging the change-plane regression model (12) into  $E(Y|T, X)$, we can express $\Delta(x)$ as follows: 

\begin{eqnarray*}
\Delta(x) &=& E\{Y^\ast(1)|X=x\}-E\{Y^\ast(0)|X=x\} \\
	&=& E(Y|T=1,X=x)-E(Y|T=0,X=x)\\
          &=& \ind\{\omega'x-\gamma \leq 0\} \{E(Y|T=1, x) -E(Y|T=0, x)  \} \\
           &&  +\ind\{\omega'x-\gamma > 0\}\{E(Y|T=1,x) -E(Y|T=0,x)  \}\\ 
            & =& \ind\{\omega'x-\gamma \leq 0\} \{(\beta_0 + \beta_T \times 1 + \beta_A \times A +\beta_H \times H )-(\beta_0 +  \beta_T \times 0 +  \beta_A \times A + \beta_H \times H) \} \\
             &&+ \ind\{\omega'x-\gamma > 0\} \{(\eta_0 + \eta_T \times 1 +\eta_A \times A +\eta_H \times H)-(\eta_0 +\eta_T \times 0 + \eta_A \times A+\eta_H \times H ) \} \\
              &= &\ind\{\omega'x-\gamma \leq 0\} \beta_T + \ind\{\omega'x-\gamma > 0\}\eta_T\\
              &=& \left[ 1-\ind\{\omega'x-\gamma > 0\} \right] \beta_T + \ind\{\omega'x-\gamma > 0\}\eta_T\\
              &=&\beta_T +\ind\{\omega'x-\gamma > 0\}(\eta_T -\beta_T) .
\end{eqnarray*}

This implies that, for the patients with $\omega'x-\gamma =\omega_A A+\omega_H H -\gamma \leq 0$, the treatment effect is explained by $\beta_T$; for the patients with $\omega'x-\gamma =\omega_A A +\omega_H H -\gamma > 0$, the treatment effect is explained by $\eta_T$. Hence, the treatment effect differs depending on A (age) and H (homosexual activity), and the heterogeneous subgroups are identified through the change-plane,  $\omega_A A+\omega_H H -\gamma=0$.  The difference in treatment effects between the two identified subgroups is expressed as $\eta_T -\beta_T$.

Note that, in our application, $\Delta(A, H) \geq 0$ indicates a greater expected CD4 cell count (a better immune system) under $T=1$ than that under $T=0$ for the patients with $(A, H)$. For 571 (54.6\%) patients in the data who satisfy $0.077 A - 0.997 H -1.889 > 0$, $T=1$ (ZDV $+$ ddI) is recommended over $T=0$ (ZDV $+$ zal) to maximize the CD4 cell count at 20 weeks, with an estimated treatment effect of $\hat{\Delta}=\hat{\eta}_T-\hat{\beta}_T=69.94$. For the other subgroup of 475 patients, $T=0$ is recommended, as $T=1$ compared to $T=0$ has a negative estimated treatment effect, $\hat{\Delta}=\hat{\beta}_T=-6.50$. The difference of treatment effects in two identified subgroups, $\hat\eta_T - \hat\beta_T = 69.94$, is relatively large to the degree that the confidence interval, $(37.7, 103.1)$, does not contain zero ($p < 0.002$).


\bibliographystyle{imsart-nameyear} 
\bibliography{mcp2.bib}       

\end{document}